\documentclass[12pt]{article}
\usepackage{graphicx}
\usepackage{graphics}
\usepackage{comment}
\usepackage{epsfig}
\usepackage{amsmath}
\usepackage{amssymb}
\usepackage{epsf}
\usepackage{authblk}
\usepackage{fullpage}

\usepackage{xcolor} 

\def\pist#1#2{\noindent\hangindent 2em\hangafter1\hbox to 2em{#1\hfil~~}#2}

\title{\vspace{-.6in}Global existence and singularity formation for the generalized Constantin-Lax-Majda equation with dissipation: The real line vs. periodic domains}

\author[1]{David M. Ambrose\thanks{dma68@drexel.edu}}
\author[2]{Pavel M. Lushnikov\thanks{plushnik@unm.edu}}
\author[3]{Michael Siegel\thanks{misieg@njit.edu (corresponding author)}}
\author[4]{Denis A. Silantyev \thanks{dsilanty@uccs.edu}}

\affil[1]{Department of Mathematics, Drexel University, Philadelphia, PA
19104, USA$~~~~~~~$} \affil[2]{Department of Mathematics and Statistics,
University of New Mexico,     Albuquerque, MSC01 1115,  NM, 87131, USA}
\affil[3]{Department of Mathematical Sciences and Center for Applied
Mathematics and$~~~~~$ Statistics, New Jersey Institute of Technology,
Newark, NJ 07102, USA} \affil[4]{Department of Mathematics, University of
Colorado, Colorado Springs, CO 80918,$~$ USA }

\newcommand{\be}{\begin{equation}}
\newcommand{\ee}{\end{equation}}
\newcommand{\bx}{\mathbf{x}}
\newcommand{\by}{\mathbf{y}}

\newcommand{\bu}{\mathbf{u}}
\newcommand{\dH}{\dot{H}^s}
\newcommand{\dHsig}{\dot{H}^\eta}
\newcommand{\dHm}{\dot{H}^m}
\newcommand{\dHr}{\dot{H}^r}
\newcommand{\dHH}{\dot{H}^1}

\newcommand{\emL}{e^{- t \mathcal{L}}}
\newcommand{\dX}{X_\infty}
\newcommand{\dXs}{X_\infty^\eta}
\newcommand{\omt}{{\omega_{-2}}}
\newcommand\bom{\boldsymbol{\omega}}
\newcommand{\tomega}{{\tilde{\omega}}}

\newcommand{\DETAILS}[1]{}
\newcommand{\I}{\mathrm{i}}
\newcommand\D{\mbox{d}}

\newcommand{\R}{\mathbb R}

\newcommand{\e}{\eqref}

\newcommand{\revc}[1]{{\color{red}\sout{}}}  
\newcommand{\revmike}[1]{{\color{cyan}}} 

\newtheorem{theorem}{Theorem}[section]
\newtheorem{lemma}[theorem]{Lemma}

\begin{document}

\maketitle

\begin{abstract}
The  question of global existence versus finite-time singularity formation
is considered for the generalized Constantin-Lax-Majda equation with
dissipation $-\Lambda^\sigma$, where $\widehat
{{\Lambda}^\sigma}=|k|^\sigma$,  both for the problem on the circle  $x
\in [-\pi,\pi]$ and the real line.   In the periodic geometry,   two
complementary approaches are used to prove global-in-time existence of
solutions  for  $\sigma \geq 1$ and all real values of an advection
parameter $a$ when the data is small.
We also derive new analytical  solutions in both geometries when $a=0$,
and on the real line when  $a=1/2$, for various values of $ \sigma$. These
solutions exhibit  self-similar finite-time  singularity formation, and
the similarity exponents and  conditions for singularity formation are
fully characterized.   We revisit an analytical solution on the real line
due to Schochet for $a=0$ and $\sigma=2$, and reinterpret it terms of
self-similar finite-time collapse.  The analytical solutions on the real
line allow finite-time singularity formation for arbitrarily small data,
even  for  values of $\sigma$ that are greater than or equal to one,
thereby illustrating a critical difference between the problems on the
real line and the circle. The analysis is complemented by accurate
numerical simulations, which are able to track the formation and motion of
singularities in the complex plane. The computations validate and extend
the analytical theory.

\end{abstract}

\section{Introduction}

In this paper we investigate global well-posedness and singularity
formation for the generalized Constantin-Lax-Majda (gCLM) model with
dissipation,
\begin{equation}
\begin{aligned}\label{CLM}
&\tilde{\omega}_{t}  =-au\tilde{\omega}_{x}+\tilde{\omega}
{\cal{H}}(\tilde{\omega})- \nu  {\Lambda}^{\sigma}(\tilde{\omega}),
\qquad \tilde{\omega} \in \mathbb{R}, x \in \mathbb{S} \mbox{~or~} \mathbb{R}, t>0, \\
&u_x =\cal{H} \tilde{\omega}, \\
&\tilde\omega(x,t) \rightarrow 0 \mbox{~for~} x \rightarrow \pm \infty \mbox{~when~} x \in \mathbb{R},\\
&\tilde\omega(x,0)=\tilde\omega_{0}(x).
\end{aligned}
\end{equation}
The equation is considered on both the circle $\mathbb{S}$ for $x \in
[-\pi, \pi]$  and the real line $\mathbb{R}$. Here $\cal{H}$ is the usual
Hilbert transform, which in the periodic case takes the form
\be \nonumber
{\cal H} f (x)= \frac{1}{2 \pi} PV \int_{-\pi}^{\pi}  f(x') \cot
\left(\frac{x-x'}{2} \right) \  dx', \ee while for the problem on the real
line \be \nonumber {\cal H} f (x)= \frac{1}{\pi} PV
\int_{-\infty}^{\infty} \frac{f(x')} {x - x'} \  dx'. \ee The operator
$\Lambda$ is given by ${\cal H} \partial_{x}$. The Hilbert transform has
Fourier symbol
\begin{equation}\nonumber
\hat{\cal{H}}=-i\mathrm{sgn}(k),
\end{equation}
so that the symbols of $\Lambda$ and $\Lambda^{\sigma}$ are
\begin{equation}\nonumber
\hat{\Lambda}(k)=|k|,\qquad \widehat{\Lambda^{\sigma}}(k)=|k|^{\sigma}.
\end{equation}
Note that $-\Lambda^2$ gives the usual diffusion operator $\partial_{xx}$,
and $-\Lambda^\sigma$ represents a generalized dissipation. The equation
$u_{x}={\cal H}(\tilde{\omega})$ defines $u$ up to its mean, and we take
the mean of $u$ to equal zero. The parameters $a$, $\sigma$ and $\nu$
satisfy $a\in\mathbb{R}$, $\sigma>0$ and $\nu>0$.

Constantin et al.\ \cite{ConstantinLaxMajda} first introduced (\ref{CLM})
with $a=\nu=0$ as a simple 1D model to study finite-time singularity
formation in the 3D incompressible Euler equations. It was later
generalized by DeGregorio to include an advection term $u
\tilde{\omega}_x$. Okamoto et al.\ \cite{Okamoto2008} introduced the
generalized advection term $a u \tilde{\omega}_x$, with real parameter
$a$, to investigate different relative weights of advection  and vortex
stretching,  $\tilde{\omega} {\cal{H}}(\tilde{\omega})$. This generalized
advection  is  motivated by recent studies of potential singularity
formation in Euler and Navier-Stokes systems, which show that advection
can have an unexpected smoothing effect  \cite{hou2018potential,
hou2014finite,  hou2012singularity, Lei2019, OkamotoOkhitani}. We will
refer to the Okamoto et al.\ model as the generalized Constantin-Lax-Majda
(gCLM) equation.   A diffusion term (can be also called by a viscosity
term) $-\Lambda^2 \tilde{\omega}=\partial_x^2 \tilde{\omega}$ was first
introduced into the Constantin et al.\ model (with $a=0$) by Schochet
\cite{Schochet}. When $a=-1$ the gCLM equation with generalized
dissipation is equivalent to the Cordoba-Cordoba-Fontelos equation
\cite{CCF},  which has been extensively studied. For $\sigma>2$ one can
interpret the term $- \nu  {\Lambda}^{\sigma}(\tilde{\omega})$ in
(\ref{CLM}) as a hyperviscosity which is widely used in many applications,
see e.g. Ref. \cite{SilantyevLushnikovRosePartIPhysPlasm2017} for the
hypervicosity in high temperature plasmas.

The dissipative gCLM system (\ref{CLM}) with $\sigma=2$ can be considered
as a 1D model of the incompressible Navier-Stokes equations, which are
written in terms of the vorticity $\boldsymbol{\omega}=\nabla \times
\boldsymbol{u}$ as
\begin{align}
\partial_t \bom &+ \bu \cdot \nabla \bom = \bom \cdot \nabla \bu + \nu \nabla^2 \boldsymbol{\omega},~~~ \bx \in \mathbb{R}^3 \mbox{~or~} \mathbb{S}^3  ,~t>0, \label{eq:vorticity}\\
\bu&= \nabla \times  (- \Delta)^{-1} \bom. \label{eq:bs}
\end{align}
The second equation above is the Biot-Savart law, which in free-space has
an equivalent representation as a convolution integral
\begin{equation} \label{eq:bsint}
\bu(\bx,t)= \frac{1}{4 \pi} \int_{\mathbb{R}^3} \frac{(\bx-\by) \times
\bom(\by,t)}{|\bx-\by|^3} \ d\by.
\end{equation}
The term $\bom \cdot \nabla \bu$ on the right-hand side  of
(\ref{eq:vorticity}) is known as the vortex stretching term, and $\nabla
\bu$ can be represented via (\ref{eq:bsint}) as a matrix of singular
integrals, which we denote by $S(\bom)$. The dissipative  gCLM equation
with $\sigma=2$  is obtained from (\ref{eq:vorticity})-(\ref{eq:bsint})
by replacing the advection term $\mathbf{u} \cdot \nabla
\boldsymbol{\omega}$ with $a u \tomega_x$, the vortex stretching term
$S(\boldsymbol{\omega}) \boldsymbol{\omega} $ by its 1D analogue $ {\cal
H} (\tomega) \tomega$, and the diffusion term by $\tomega_{xx}$.  The
Hilbert transform is the unique linear singular integral operator in 1D
that, like $S(\boldsymbol{\omega})$,  commutes with translations and
dilations \cite{ConstantinLaxMajda}. This motivates the  replacement of
$S(\boldsymbol{\omega})$ from the 3D equations with ${\cal H} (\tomega)$
in the 1D model.

Singularities to (\ref{CLM}), when they occur, are generally found to be locally self-similar with the form%
\begin{equation}\label{self-similar1}
\tilde{\omega}=\frac{1}{\tau^\beta}f\left ( \xi \right ), \
\xi=\frac{x-x_0}{\tau^\alpha}, \ \tau=t_c-t,
\end{equation}
in a space-time neighborhood of $(x_0,t_c)$, where $t_c>0$\ is the
singularity time, $x_0 \in \mathbb{R}$ is its location, and $\alpha$,
$\beta$ are real similarity parameters. There are a number of results on
finite-time singularity formation in the inviscid problem for (\ref{CLM})
with $\nu=0$, which we now briefly describe; see
\cite{Lushnikov_Silantyev_Siegel} for a more complete review.  In this
case of $\nu=0$, one has that $\beta=1$ while $\alpha$ depends on $a$.
Constantin, Lax and Majda  \cite{ConstantinLaxMajda} present a closed-form
exact solution to  the initial value problem
for (\ref{CLM}) with $a=0$. Their solution develops a singularity of the
local form (\ref{self-similar1}) with $\alpha=\beta=1$ for a  class of
analytic initial data.    Castro and Cordoba \cite{CastroCordoba} prove
finite-time blow-up for  $a<0$  using a Lyapunov-type argument.
For  $\epsilon-$small values of $a>0$,  Elgindi and Jeong
\cite{Elgindi2020} and Chen et al. \cite{Chen2021} prove the existence of
singularities of the form (\ref{self-similar1}) with $\beta=1$ and
$\alpha$ approaching $1$ in the limit $a\to 0^+$.

More recently,   \cite{Lushnikov_Silantyev_Siegel} and
\cite{Chen2020Singularity} independently find an exact self-similar
solution to the inviscid problem as a superposition of double-pole
singularities for $a=1/2$ with $\alpha=1/3$ and $\beta=1$
(\cite{Lushnikov_Silantyev_Siegel} further show that, beyond the
particular cases $a=0$ and $a=1/2$, no exact solutions as a superposition
of pole singularities exist).    Lushnikov et al.\
\cite{Lushnikov_Silantyev_Siegel} also perform numerical simulations over
a wide range of $a$ and find the existence of a critical value $a_c=0.6890
\ldots$ for which the  self-similar blow up of solutions
changes character. More precisely, they find  self-similar collapse with
$\alpha>0$ when $a<a_c$ for both $x \in \mathbb{S}$ and $\mathbb{R}$,
expanding self-similar blow up with $\alpha<0$ when $a_c<a \leq 1$ and $x
\in \mathbb{R}$, and `neither expanding nor collapsing' blow-up with
$\alpha=0$ when $a_c<a\le 0.95$ and $x \in \mathbb{S}$ (with the
expectation that the latter behavior occurs for $a$ going all the way up
to, but not including, $a=1$).  Here the terminology  ``collapse'' or
``wave collapse'' was first introduced in  \cite{ZakharovJETP1972a} in
analogy with  gravitational collapse and has been widely used ever since
to mean that the solution shrinks in $x$ as $t\to t_c$ while its amplitude
diverges in that limit; see
\cite{ZakharovJETP1972a,ChPe1981,SulemSulem1999,brenner1999diffusion,KuznetsovZakharov2007,DyachenkoLushnikovVladimirovaKellerSegelNonlinearity2013,LushnikovDyachenkoVladimirovaNLSloglogPRA2013}
for a more general description.  Existence of the expanding similarity
solution for $x \in \mathbb{R}$ and the  `neither expanding nor
collapsing'  similarity solution for $x \in \mathbb{S}$  are proven in
\cite{Chen2020Singularity}, \cite{Chen2021}, when $a$ is near $1^-$.
Analytical \cite{JiaStewartSverak} and numerical
\cite{Lushnikov_Silantyev_Siegel} evidence is consistent with  global
well-posedness when $a \geq 1$ in the periodic problem, and $a >1$ in the
problem on the real line. However,  at present there is  no proof of this
for general analytic or $C^\infty$ initial data.

Much less is known about solutions to (\ref{CLM}) when there is nonzero
dissipation.  Schochet \cite{Schochet} constructs an explicit solution  on
the real line for $a=0$ and $\sigma=2$,    which blows up in finite time.
When $a=-1$, so that (\ref{CLM}) is the  Cordoba-Cordoba-Fontelos
equation, finite time blow up can occur for $\sigma<1/2$  \cite{Kiselev,
LiRodrigo, SilvestreVicol}, although there is global well-posedness for
sufficiently small data \cite{Dong2008}.   Global well-posedness  of the
CCF equation for   $\sigma \geq 1$ is shown in \cite{CCF, Dong2008,
Kiselev}.   When $a \leq -2$ is even and $\sigma=1$, global well-posedness
for small data in the periodic setting is shown in \cite{Wunsch}. More
recently, Chen \cite{Chen2020Singularity} shows that for the problem on
the real line, there exists self-similar blow up when $a$ is close to
$1/2$ and $\sigma=2$, and global well-posedness for  $\sigma \in
[|a|^{-1}, 2]$ with $a<-1$.  We note that for $a>-1$, there is no known
coercive conserved quantity for general initial data, which complicates
attempts to prove global well-posedness.

The focus of this paper is to further investigate conditions under which
(\ref{CLM}) is well-posed globally in time, for different values  of the
parameters $a$ and $\sigma$.
We find a surprising dependence of the global well-posedness on the domain
of $x$, i.e., whether it is $\mathbb{S}$ or $\mathbb{R}$.   In particular,
we prove that the initial value problem (\ref{CLM}) with $\sigma \geq 1$
has global-in-time solutions for all sufficiently small data and all $a
\in \mathbb{R}$, when the problem is considered on the periodic domain $x
\in \mathbb{S}$.  These solutions are analytic for $t>0$.

We present examples in the form of  exact analytical solutions and
numerical simulations which show this result does not hold on the real
line $x \in \mathbb{R}$.
The analytical results include new  `pole dynamics' solutions to
(\ref{CLM}) -- there are numerous examples of exact pole dynamics
solutions in  both Hamiltonian and dissipative systems, see e.g.,
\cite{CalogeroBook2001,Senouf1CaflischErcolaniNonlinearity1996,LushnikovPhysLettA2004,LushnikovZubarevPRL2018}.
Our exact solutions for the problem on the real line form finite-time
singularities of the type \eqref{self-similar1} for arbitrarily  small
initial data  in $L^2$ and usually also in $L^\infty$. They include: (1) a
solution for $a=1/2$ and $\sigma=1$  expressed as the sum of a complex
conjugate (c.c.) pair of  second order poles in $\omega$,  (2) solutions
for $a=0$ and $\sigma=1$ expressed as the sum of one or two c.c.\ pairs of
first order poles in $\omega$, and (3) a solution for $a=0$ and $\sigma=0$
expressed as the sum of a c.c.\ pair of first order poles in $\omega$.
We also revisit and slightly correct a previous example due to Schochet
for $a=0$ and $\sigma=2$, which forms singularities in finite-time from
arbitrarily small data, and reinterpret it as self-similar blow up.
Overall, the exact solutions display different similarity exponents
$\alpha$ and $\beta$, depending on the location and `strength' (i.e.,
power or exponent)  of their poles in the complex plane, and whether they
impinge on the real line with a nonzero or zero velocity.

Additionally, we find a new pole dynamics solution to the periodic problem
for $a=0$ with `marginal' dissipation $\sigma=0$. This solution consists
of a c.c.\ pair of simple poles and can form a finite-time singularity of
the form (\ref{self-similar1}) for data which is arbitrarily small in
$L^2$, but not necessarily small in  $L^\infty$. This supplies a lower
bound in $\sigma$ for which a global existence theory in $L^2$ can apply,
when $\nu$ is nonzero.

The analysis is complemented by accurate numerical simulations which
confirm and build upon the analytical results in  the periodic and real
line problems. As part of the numerics, the formation and motion of
singularities is tracked in the complex plane.  When a singularity reaches
the real line (at time $t_c$) a finite-time singularity of the form
(\ref{self-similar1}) occurs. We make use of two methods to trace
singularities in the complex plane. One is based on the asymptotic decay
of Fourier amplitudes, which gives precise (quantitative) information on
the singularity that is closest to the real line. The other method, known
as the AAA algorithm \cite{TrefethenAAA}, utilizes rational function
approximation to obtain information on  singularities beyond the one
closest to the real line.

 Our analysis of the periodic problem makes use of two  complementary approaches. We first prove that when $\sigma\geq1,$ the solution exists globally in time  for small initial data in the periodic Wiener algebra, which describes the set of functions with Fourier coefficients in $l^1$. A consequence of the proof is that solutions are analytic at all positive times in a strip in the complex plane that contains the real line, with the width of the strip growing linearly in time.  The proof employs the method of Duchon and Robert \cite{DuchonRobert}, who developed it to show the existence of global vortex sheet solutions for certain types of small data. Other applications of this method to show global existence are \cite{AmbroseRadius, AmbroseMazzucato}.

 We also prove  global-in-time existence of mild solutions with small initial data in $L^2$, when $\sigma>1$.  A particular challenge in the proof is to obtain an exponential decay estimate for the solution operator when $t\gg 1$.  We are able to do this, but the result relies in an essential way on the periodicity of the geometry.
The proof guarantees that
the solution at any time $t>0$ exists in $H^\gamma$ for all
$1/2<\gamma<\min [1,\sigma-1/2]$.
 We further expect that solutions become analytic for $t>0$, even starting from rough $L^2$ data. This can be shown using the approach of Grujic and Kukavica \cite{GrujicKukavica}, which has been used in several  related problems to show analyticity of solutions on a strip which grows initially like $t^{1/\sigma}$ (see e.g., \cite{AmbroseMazzucato}). We do not provide details, and instead refer the interested reader to the relevant work.

The rest of this paper is organized as follows. After some mathematical
preliminaries in \S \ref{sec:Preliminaries}, a solution operator is
written in  \S \ref{sec:Solutionoperator} using the Duhamel
representation. Section \ref{sec:globalsolutionsWieneralgebra}   proves
global existence for small periodic initial data with  $\sigma \geq 1$
as a fixed point of the Duhamel representation by  using a Wiener algebra
approach. Section \ref{sec:Mildsolutionapproach} proves  global existence
for small periodic initial data in $L^2$ with $\sigma>1$ using a mild
solution approach. Section \ref{sec:real_line} focuses on the derivation
of exact solutions on the real line and their relation to the self-similar
form (\ref{self-similar1}).
Section \ref{sec:periodic} derives an exact solution to the periodic
problem for $a=0$ and $\sigma=0$ which can develop a  finite-time
singularity for arbitrarily small data in $L^2$. Section
\ref{sec:numerical_results} presents numerical results, with the numerical
method described in \S \ref{sec:numerical_method}, numerical results for
the  periodic problem given in \S \ref{sec:periodic_numerics}, and
numerical results for  the problem on the real line discussed in \S
\ref{sec:real_line_numerics}. Concluding remarks are given in \S
\ref{sec:conclusion}. An appendix, \S \ref{sec:appendix}, provides a proof
of inequality (\ref{desiredBound})  used in the Wiener space analysis, and
Lemmas \ref{lemma:int_est} and \ref{lemma:sumest} used in the mild
solution analysis.

\section{Preliminaries}
\label{sec:Preliminaries}

By rescaling each of $t$ and $\tilde{\omega},$ we can eliminate $\nu$ from
the problem.  We therefore set $\nu=1$ without loss of generality, unless
otherwise noted.

Notice that for any periodic function $f,$ we have $\int_{\mathbb{S}} f
{\cal H}(f)\ dx=0.$ Also, $u  \tilde{\omega}_x = (u \tilde{\omega})_x -
\tilde{\omega} \mathcal{H} \tilde{\omega}$ has zero mean. Thus the mean of
$\tilde{\omega}$   is preserved under the evolution \eqref{CLM} on the
circle. In the periodic problem, we make the decomposition
$\tilde{\omega}=\omega+{\omega}_{av},$ where ${\omega}_{av}$ is the mean
of $\tilde{\omega}$ and $\omega$ has zero mean. Substituting this
decomposition in \eqref{CLM} yields
\begin{equation}\label{mainEquationTemporary}
(\omega+{\omega}_{av})_{t}+au(\omega+{\omega}_{av})_{
x}=(\omega+{\omega}_{av})  \mathcal{H}(\omega+{\omega}_{av})
-\Lambda^{\sigma}(\omega+{\omega}_{av}) \mbox{~for~} x \in \mathbb{S},
\end{equation}
with $u$ now being defined through $u_{x}=\mathcal{H}(\omega);$ this is
the same as the previous formula since the periodic Hilbert transform of a
constant function is equal to zero. Since
$({\omega}_{av})_{t}=({\omega_{av}})_{x}=\mathcal{H}({\omega}_{av})=\Lambda^{\sigma}({\omega}_{av})=0,$
we can rewrite \eqref{mainEquationTemporary} as
\begin{equation} \label{mainEquation2}
\omega_{t}+au\omega_{x}=\omega \mathcal{H}(\omega)+{\omega}_{av}
\mathcal{H}(\omega)-\Lambda^{\sigma}\omega  \mbox{~for~} x \in \mathbb{S},
\end{equation}
with initial data $\omega(x,0)=\omega_0(x)$, which are used instead of the
first and last equation in  (\ref{CLM}) for the periodic problem.  We
continue to use  (\ref{CLM}) for the problem on the real line, but omit
the tilde from $\omega$, with the understanding that when $x \in
\mathbb{R}$ the function $\omega$ is allowed to  have a nonzero mean.

Notice the Hilbert transform also has the representation
\begin{equation} \label{Hilbert}
\mathcal{H} \omega=-i(\omega_+ - \omega_-),
\end{equation}
where $\omega=\omega_+ + \omega_-$ with $\omega_+$  analytic in the upper
complex half-plane $\mathbb{C}^+$,  and $\omega_-$ is analytic in the
lower  complex half-plane $\mathbb{C}^-$. In the periodic problem,
$f_+=\sum_{k>0} \widehat{f}_k e^{i k x}$ and $f_-=\sum_{k<0} \widehat{f}_k
e^{i k x}$ are the projections onto the upper and lower analytic
components of $f$, respectively.

\subsection{Solution operator in the periodic case}
\label{sec:Solutionoperator}

The solution to (\ref{mainEquation2}) can be written using the Duhamel
representation \be \label{Duhamel} \omega(\cdot, t)=e^{- t \mathcal{L}}
\omega_0 + \int_0^t e^{- (t-\tau) \mathcal{L}} \left(-a u \omega_x +
\omega u_x \right)(\cdot, \tau) \ d \tau, \ee in which  the operator $e^{-
t \mathcal{L}}$ is defined by \be \label{emLdef} e^{- t \mathcal{L}} f =
\mathcal{F}^{-1}  \left(e^{- t |k|^\sigma-it {\omega_{av}}\mathrm{sgn}(k)}
\widehat{f}(k) \right) \ee where $\mathcal{F}$ is the Fourier transform
operator.
It is helpful to rewrite \eqref{Duhamel} slightly; we do so by first
rewriting \eqref{mainEquation2} using
\begin{equation}\nonumber
u\omega_{x}=(u\omega)_{x}-u_{x}\omega=(u\omega)_{x}-\omega
\mathcal{H}(\omega),
\end{equation}
leading to
\begin{equation}\nonumber
\omega_{t}=(1+a)\omega
\mathcal{H}(\omega)-a(u\omega)_{x}+{\omega_{av}}\mathcal{H}(\omega)-\Lambda^{\sigma}\omega.
\end{equation}

We again rewrite this using Duhamel's principle, finding
\begin{equation}\label{almostOurDuhamel}
\omega=e^{-t\mathcal{L}}\omega_{0} +\int_{0}^{t}e^{-(t-\tau)\mathcal{L}}
\left[(1+a)\omega \mathcal{H}(\omega)-a(u\omega)_{x}\right](\cdot,\tau)\
d\tau.
\end{equation}
We use again the fact that for any periodic function $f,$ the integral
$\int_{\mathbb{S}} f \mathcal{H}(f)\ dx=0;$  introducing the operator
$\mathbb{P}_{0}$ to be the projection which zeroes out the mean of a
periodic function, we have
\begin{equation}\nonumber
\mathbb{P}_{0}\left[(1+a)\omega
\mathcal{H}(\omega)-a(u\omega)_{x}\right]=\left[(1+a)\omega
\mathcal{H}(\omega)-a(u\omega)_{x}\right].
\end{equation}
We then use this with \eqref{almostOurDuhamel} as the basis for
introducing an operator $\mathcal{T},$
\begin{equation}\nonumber
\mathcal{T}(\omega)=e^{-t\mathcal{L}}\omega_{0}
+\int_{0}^{t}e^{-(t-\tau)\mathcal{L}}\mathbb{P}_{0} \left[(1+a)\omega
\mathcal{H}(\omega)-a(u\omega)_{x}\right](\cdot,\tau)\ d\tau.
\end{equation}
We will obtain solutions of the gCLM equation with dissipation by finding
a fixed point of  $\mathcal{T}$. As we have said above, we will do this
twice, once in function spaces based on the Wiener algebra, and once in
$L^{2}$-based Sobolev spaces.

\section{Small global solutions in spaces based on the Wiener algebra}
\label{sec:globalsolutionsWieneralgebra}

In this section we will prove global existence of small solutions when the
initial data is taken from the Wiener algebra.  This uses an adaptation of
the argument of Duchon and Robert used to prove existence of small global
vortex sheets \cite{DuchonRobert}.  The unregularized vortex sheet is an
elliptic problem in space-time, but the method has also been applied to
parabolic problems in \cite{AmbroseRadius, AmbroseMazzucato}.

\subsection{Function spaces and operators}

We denote the periodic Wiener algebra as $B_{0};$ this is the set of
functions $f:\mathbb{S}\rightarrow\mathbb{R}$ such that the norm
\begin{equation}\nonumber
\|f\|_{B_{0}}=\sum_{k\in\mathbb{Z}}|\hat{f}(k)|
\end{equation}
is finite.

For $\varpi>0$ and $\theta\geq0,$ we define the function space
$\mathcal{B}_{\varpi}^{\theta}$ to be the set of periodic functions
continuous in time with values in $B_{0},$ such that the norm
\begin{equation}\nonumber
\|h\|_{\varpi,\theta}=\sum_{k\in\mathbb{Z}}(1+|k|^{\theta})\sup_{t\in[0,\infty)}e^{\varpi
t|k|}|\hat{h}(k,t)|
\end{equation}
is finite.  We will demonstrate that this is a Banach algebra.  First,
note that for all $k\in\mathbb{Z},$ for all $j\in\mathbb{Z},$ we have
\begin{equation}\label{desiredBound}
|k|^{\theta}\leq
\max\{1,2^{\theta-1}\}\left(|k-j|^{\theta}+|j|^{\theta}\right).
\end{equation}
(We prove this inequality in Appendix \ref{inequalityAppendix}.)  We
denote $C=\max\{1,2^{\theta-1}\}.$ We compute the norm of $fg,$ for
$f\in\mathcal{B}_{\varpi}^{\theta}$ and
$g\in\mathcal{B}_{\varpi}^{\theta}:$
\begin{multline}\nonumber
\|fg\|_{\varpi,\theta}=\sum_{k\in\mathbb{Z}}(1+|k|^{\theta})\sup_{t\in[0,\infty)}e^{\varpi
t|k|}|\widehat{(fg)}(k,t)|
\\
\leq C
\sum_{(k,j)\in\mathbb{Z}^{2}}\left(1+|k-j|^{\theta}+|j|^{\theta}\right)
\left[\sup_{t\in[0,\infty)}e^{\varpi t |k-j|}|\hat{f}(k-j,t)|\right]
\left[\sup_{t\in[0,\infty)}e^{\varpi t|j|}|\hat{g}(j,t)|\right].
\end{multline}
We sum first in $k$ and then in $j,$ finding
\begin{equation}\nonumber
\|fg\|_{\varpi,\theta}\leq C\|g\|_{\varpi,0}\|f\|_{\varpi,\theta}
+C\|f\|_{\varpi,0}\|g\|_{\varpi,\theta}\leq
2C\|f\|_{\varpi,\theta}\|g\|_{\varpi,\theta}.
\end{equation}

For any $\varpi>0$ and $\theta\geq0,$ we let
$\mathcal{B}_{\varpi,0}^{\theta}$ be the subspace of
$\mathcal{B}_{\varpi}^{\theta}$ of functions with zero mean. We then
define the integral operator
$I^{+}:\mathcal{B}_{\varpi}^{\theta}\rightarrow\mathcal{B}_{\varpi,0}^{\theta+\sigma}$
by
\begin{equation}\nonumber
(I^{+}h)(\cdot,t)=\int_{0}^{t}e^{-(t-\tau)\mathcal{L}}\mathbb{P}_{0}h(\cdot,\tau)\
d\tau.
\end{equation}

We compute the operator norm of $I^{+}.$ The norm for
$\mathcal{B}_{\varpi,0}^{\theta+\sigma}$ is the same as for
$\mathcal{B}_{\varpi}^{\theta+\sigma},$ except that the $k=0$ mode is
excluded from the summation.  Therefore we have
\begin{equation}\nonumber
\|I^{+}h\|_{\varpi,\theta+\sigma}=\sum_{k\in\mathbb{Z}\setminus\{0\}}(1+|k|^{\theta+\sigma})
\sup_{t\in[0,\infty)}e^{\varpi t|k|}\left|\int_{0}^{t}e^{-
|k|^{\sigma}(t-\tau)-i{\omega_{av}}\mathrm{sgn}(k)(t-\tau)}
\hat{h}(\cdot,\tau)\ d\tau\right|.
\end{equation}
We use the triangle inequality and rearrange the exponentials, finding
\begin{equation}\nonumber
\|I^{+}h\|_{\varpi,\theta+\sigma}=\sum_{k\in\mathbb{Z}\setminus\{0\}}(1+|k|^{\theta+\sigma})\sup_{t\in[0,\infty)}
e^{(\varpi|k|-|k|^{\sigma})t}\int_{0}^{t}e^{
|k|^{\sigma}\tau}|\hat{h}(k,\tau)|\ d\tau.
\end{equation}
We adjust factors of the weights, arriving at
\begin{multline}\nonumber
\|I^{+}h\|_{\varpi,\theta+\sigma}=\sum_{k\in\mathbb{Z}\setminus\{0\}}
\left(\frac{1+|k|^{\theta+\sigma}}{1+|k|^{\theta}}\right)
\Bigg(\sup_{t\in[0,\infty)} e^{(\varpi|k|-|k|^{\sigma})t}\cdot
\\
\cdot
\int_{0}^{t}e^{(|k|^{\sigma}-\varpi|k|)\tau}\left[(1+|k|^{\theta})e^{\varpi|k|\tau}|\hat{h}(k,\tau)|\right]\
d\tau\Bigg).
\end{multline}
We estimate this by taking the supremum two more times, once with respect
to $\tau$ and once with respect to $k,$ and then rearranging:
\begin{multline}\nonumber
\|I^{+}h\|_{\varpi,\theta+\sigma}\leq\left(\sum_{k\in\mathbb{Z}\setminus\{0\}}(1+|k|^{\theta})\sup_{\tau\in[0,\infty)}
e^{\varpi|k|\tau}|\hat{h}(k,\tau)|\right)
\\
\left(\sup_{k\in\mathbb{Z}\setminus\{0\}}\left(\frac{1+|k|^{\theta+\sigma}}{1+|k|^{\theta}}\right)
\sup_{t\in[0,\infty)}e^{(\varpi|k|-|k|^{\sigma})t}\int_{0}^{t}e^{(|k|^{\sigma}-\varpi|k|)\tau}\
d\tau \right).
\end{multline}
We identify the first factor on the right-hand side as simply being
$\|h\|_{\varpi,\theta},$ and we evaluate the last integral and simplify.
These considerations yield the following:
\begin{equation}\nonumber
\|I^{+}h\|_{\varpi,\theta+\sigma}\leq\|h\|_{\varpi,\theta}
\left(\sup_{k\in\mathbb{Z}\setminus\{0\}}\left(\frac{1+|k|^{\theta+\sigma}}{1+|k|^{\theta}}\right)
\sup_{t\in[0,\infty)}
\frac{1-e^{(\varpi|k|-|k|^{\sigma})t}}{|k|^{\sigma}-\varpi|k|} \right).
\end{equation}
The last denominator on the right-hand side is positive as long as
$\sigma\geq1$ and $\varpi<1.$ With these conditions, we may then ignore
the negative term in the numerator, arriving at
\begin{equation}\nonumber
\|I^{+}h\|_{\varpi,\theta+\sigma}\leq\|h\|_{\varpi,\theta}\left(\sup_{k\in\mathbb{Z}\setminus\{0\}}
\frac{1+|k|^{\theta+\sigma}}{(1+|k|^{\theta})(|k|^{\sigma}-\varpi|k|)}\right).
\end{equation}
We then estimate this using $1\leq |k|^{\sigma}$ and simplifying, arriving
at
\begin{equation}\label{I+Estimate}
\|I^{+}h\|_{\varpi,\theta+\sigma}\leq\frac{\|h\|_{\varpi,\theta}}{1-\varpi}.
\end{equation}

We need an entirely analagous bound for the composition operator
$I^{+}\partial_{x}.$  The above proof of the estimate \eqref{I+Estimate}
works just the same to show that $I^{+}\partial_{x}$ maps
$\mathcal{B}_{\varpi}^{\theta}$ to
$\mathcal{B}_{\varpi}^{\theta+\sigma-1},$ with the estimate
\begin{equation}\label{I+dxEstimate}
\|I^{+}\partial_{x}h\|_{\varpi,\theta+\sigma-1}\leq\frac{\|h\|_{\varpi,\theta}}{1-\varpi}.
\end{equation}

We also need to demonstrate the boundedness of the semigroup, acting on
$B_{0}.$ Letting $h\in B_{0},$ we consider the norm of
$e^{-t\mathcal{L}}h:$
\begin{equation}\nonumber
\|e^{-t\mathcal{L}}h\|_{\varpi,0}
=2\sum_{k\in\mathbb{Z}}\sup_{t\in[0,\infty)} e^{\varpi
t|k|}e^{-|k|^{\sigma}t}|\hat{h}(k)| \leq
2\|h\|_{B_{0}}\sup_{k\in\mathbb{Z}}\sup_{t\in[0,\infty)} e^{(\varpi
|k|-|k|^{\sigma})t}.
\end{equation}
With $\sigma\geq1$ and $\varpi<1,$ we may estimate this as
\begin{equation}\nonumber
\|e^{-t\mathcal{L}}h\|_{\varpi,0}\leq 2\|h\|_{B_{0}}.
\end{equation}

\subsection{Existence of a solution}

In the current notation, our operator $\mathcal{T}$ may be expressed as
\begin{equation}\label{TDefinition}
\mathcal{T}\omega=e^{-t\mathcal{L}}\omega_{0}
+\left(I^{+}\left[(1+a)\omega \mathcal{H}(\omega)-a(u\omega)_{x}\right]
\right)(t).
\end{equation}
A fixed point of \eqref{TDefinition} is a solution of the initial value
problem for \eqref{mainEquation2}.

We see that if $\omega_{0}\in B_{0}$ and if $\varpi<1$ and $\sigma\geq1,$
then $\mathcal{T}$ maps $\mathcal{B}_{\varpi}^{0}$ to itself. We want to
show that there exists $X\subseteq\mathcal{B}_{\varpi}^{0}$ such that
$\mathcal{T}$ is a contraction on $X.$  We let $X$ be the ball of radius
$r_{0}$ centered at $e^{-t\mathcal{L}}\omega_{0},$ and we denote
$r_{1}=\|\omega_{0}\|_{B_{0}}.$  We will show that $\mathcal{T}$ is a
contraction on $X$ for an appropriate choice of $r_{0}$ and $r_{1}.$  Note
that for any $\omega\in X,$ we have $\|\omega\|_{\varpi,0}\leq
r_{0}+r_{1}.$

We have two properties to establish: that $\mathcal{T}:X\rightarrow X,$
and that there exists $\lambda\in(0,1)$ such that for any $\omega_{1}\in
X$ and for any $\omega_{2}\in X,$
\begin{equation}\label{contractionCondition}
\|\mathcal{T}(\omega_{1}-\omega_{2})\|_{\varpi,0}\leq
\lambda\|\omega_{1}-\omega_{2}\|_{\varpi,0}.
\end{equation}

To show that $\mathcal{T}$ maps $X$ to $X,$ we let $\omega\in X$ be given,
and we need to establish that
\begin{equation}\nonumber
\left\|I^{+}\left[(1+a)\omega
\mathcal{H}(\omega)-a(u\omega)_{x}\right]\right\|_{\varpi,0}\leq r_{0}.
\end{equation}
We immediately have
\begin{equation}\nonumber
\left\|I^{+}\left[(1+a)\omega
\mathcal{H}(\omega)-a(u\omega)_{x}\right]\right\|_{\varpi,0} \leq
\frac{|1+a|}{1-\varpi}(r_{0}+r_{1})^{2}+|a|\|I^{+}(u\omega)\|_{\varpi,1}.
\end{equation}
We then bound this as
\begin{equation}\nonumber
\left\|I^{+}\left[(1+a)\omega
\mathcal{H}(\omega)-a(u\omega)_{x}\right]\right\|_{\varpi,0} \leq
\frac{|1+a|}{1-\varpi}(r_{0}+r_{1})^{2}
+\frac{|a|}{1-\varpi}(r_{0}+r_{1})\|u\|_{\varpi,0}.
\end{equation}
We also have $\|u\|_{\varpi,0}\leq\|\omega\|_{\varpi,0},$ so that
\begin{equation}\nonumber
\left\|I^{+}\left[(1+a)\omega
\mathcal{H}(\omega)-a(u\omega)_{x}\right]\right\|_{\varpi,0} \leq
\frac{|1+a|+|a|}{1-\varpi}(r_{0}+r_{1})^{2}.
\end{equation}
Our first condition that $r_{0}$ and $r_{1}$ must satisfy, then, is
\begin{equation}\label{firstCondition}
\frac{|1+a|+|a|}{1-\varpi}(r_{0}+r_{1})^{2}\leq r_{0}.
\end{equation}

Next, we work on establishing \eqref{contractionCondition}.  To begin, we
express the difference $\mathcal{T}(\omega_{1}-\omega_{2}),$ doing some
adding and subtracting:
\begin{multline}\nonumber
\mathcal{T}(\omega_{1}-\omega_{2})=\int_{0}^{t}e^{-(t-\tau)\mathcal{L}}\left[
(1+a)\omega_{1}\mathcal{H}(\omega_{1})-a(u_{1}\omega_{1})_{x}\right](\cdot,\tau)\
d\tau
\\
-\int_{0}^{t}e^{-(t-\tau)\mathcal{L}}[(1+a)\omega_{2}\mathcal{H}(\omega_{2})-a(u_{2}\omega_{2})_{x}](\cdot,\tau)\
d\tau
\\
=A_{1}+A_{2}+A_{3}+A_{4},
\end{multline}
where the $A_{i}$ are given by
\begin{equation}\nonumber
A_{1}=\int_{0}^{t}e^{-(t-\tau)\mathcal{L}}[(1+a)(\omega_{1}-\omega_{2})\mathcal{H}(\omega_{1})](\tau)\
d\tau,
\end{equation}
\begin{equation}\nonumber
A_{2}=\int_{0}^{t}e^{-(t-\tau)\mathcal{L}}[(1+a)\omega_{2}(\mathcal{H}(\omega_{1})-\mathcal{H}(\omega_{2})](\cdot,\tau)\
d\tau,
\end{equation}
\begin{equation}\nonumber
A_{3}=-\int_{0}^{t}e^{-(t-\tau)\mathcal{L}}a((u_{1}-u_{2})\omega_{1})_{x}(\cdot,\tau)\
d\tau,
\end{equation}
\begin{equation}\nonumber
A_{4}=-\int_{0}^{t}e^{-(t-\tau)\mathcal{L}}a(u_{2}(\omega_{1}-\omega_{2})_{x})(\cdot,\tau)\
d\tau.
\end{equation}
We may estimate these as follows:
\begin{equation}\nonumber
\|A_{1}\|_{\varpi,0}\leq\frac{|1+a|}{1-\varpi}(r_{0}+r_{1})\|\omega_{1}-\omega_{2}\|_{\varpi,0},
\end{equation}
\begin{equation}\nonumber
\|A_{2}\|_{\varpi,0}\leq\frac{|1+a|}{1-\varpi}(r_{0}+r_{1})\|\omega_{1}-\omega_{2}\|_{\varpi,0},
\end{equation}
\begin{equation}\nonumber
\|A_{3}\|_{\varpi,0}\leq\frac{|a|}{1-\varpi}(r_{0}+r_{1})\|\omega_{1}-\omega_{2}\|_{\varpi,0},
\end{equation}
\begin{equation}\nonumber
\|A_{4}\|_{\varpi,0}\leq\frac{|a|}{1-\varpi}(r_{0}+r_{1})\|\omega_{1}-\omega_{2}\|_{\varpi,0}.
\end{equation}
We combine these estimates to find
\begin{equation}\nonumber
\|\mathcal{T}(\omega_{1}-\omega_{2})\|_{\varpi,0}\leq\frac{2(|1+a|+|a|)}{1-\varpi}(r_{0}+r_{1})
\|\omega_{1}-\omega_{2}\|_{\varpi,0}.
\end{equation}
Thus, our second condition which $r_{0}$ and $r_{1}$ must satisfy is
\begin{equation}\label{secondCondition}
\frac{2(|1+a|+|a|)}{1-\varpi}(r_{0}+r_{1})<1.
\end{equation}

To demonstrate that \eqref{firstCondition} and \eqref{secondCondition} may
be satisfied, we take $r_{1}=r_{0},$ and we will choose $r_{1}.$  In this
case, \eqref{firstCondition} becomes
\begin{equation}\label{firstConditionAgain}
r_{1}\leq\frac{1-\varpi}{4(|1+a|+|a|)},
\end{equation}
while \eqref{secondCondition} becomes
\begin{equation}\label{secondConditionAgain}
r_{1}<\frac{1-\varpi}{4(|1+a|+|a|)}.
\end{equation}
Of course \eqref{secondConditionAgain} implies
\eqref{firstConditionAgain}.

We have proved the following theorem:

\begin{theorem} \label{thm_ambrose}
Let $a\in\mathbb{R},$ ${\omega_{av}}\in\mathbb{R},$ and $\sigma\geq1$ be
given. Let $\omega_{0}\in B_{0}$ be given, such that $\omega_{0}$ has zero
mean and such that $\|\omega_{0}\|_{B_{0}}<\frac{1}{4(|1+a|+|a|)}.$ Let
$\varpi\in(0,1)$ be given such that
$\|\omega_{0}\|_{B_{0}}<\frac{1-\varpi}{4(|1+a|+|a|)}.$ Then the initial
value problem for \eqref{mainEquation2} with initial data $\omega_{0}$ has
a unique solution $\omega\in\mathcal{B}_{\varpi}^{0}.$
\end{theorem}

We make a few remarks on Theorem \ref{thm_ambrose}. Since the solution is
in $\mathcal{B}_{\varpi}^{0}$ with $\varpi>0,$ we know automatically that
the solution exists for all $t\in[0,\infty),$ and that the solution is
analytic at all positive times with radius of analyticity at least $\varpi
t.$
Next, we notice that the value of $a$ does not matter as far as whether we
can get global existence of a solution, except that it does affect the
maximum allowable size of the data; specifically, for larger $|a|,$ we
need to take the data smaller. As noted above, $\varpi$ is the rate at
which analyticity is gained; if we want this to be larger, the data must
be taken smaller.  Finally we note that the value of $\omega_{av}$ does
not affect the allowable size of the data or the rate at which analyticity
is gained.

\section{Mild solutions with data in $L^{2}$}
\label{sec:Mildsolutionapproach}

In this section we complement Theorem \ref{thm_ambrose} with another
theorem on existence of small global solutions, now taking initial data in
$L^{2}.$ In this approach, we will need more detailed mapping properties
for the semigroup associated to the diffusive term than in the Wiener
algebra case; we establish these properties in Section
\ref{semigroupSection} below.

\subsection{Function spaces and preliminary Lemmas}

Throughout we use the notation $L^2, H^s$ etc$.$ to denote the spaces
$L^2[-\pi, \pi]$, $H^s[-\pi,\pi]$ (with periodic boundary conditions) and
so forth. We consider data and solutions with finite $L^2$ norm, i.e.,
finite energy. Hence, it will be convenient to work with the norm in
homogeneous Sobolev spaces $\dot{H}^s$, defined by \be \nonumber \| f
\|^2_{\dH}=\sum_{k=-\infty}^\infty |k|^{2s} |\widehat{f}(k)|^2,~~~s \in
\mathbb{R}. \ee Note that if $f \in L^2$ then $f \in \dH$ if and only if
$f \in H^s$. We denote the subspace of functions in $L^2$ with zero mean
as
\[
\mathring{L}^2=\left\{ f \in L^2 \mid \int_{-\pi}^\pi f  \ dx = 0
\right\}.
\]
If a function $f$ has zero mean, then by Poincare's inequality $\| f
\|_{L^2} \leq c \| f_x \|_{L^2}$ so that $\dHH \subset \mathring{L}^2$. In
particular,
 $\| \mathcal{H} \omega \|_{L^2}   = \| \omega_+  + \omega_{-} \|_{L^2} \leq c \| \omega_x \|_{L^2}$, but note that if a function  $f$ has nonzero mean its $L^2$ norm cannot in general  be bounded by the $L^2$ norm of  its derivative.

In the fixed point analysis,   we make use of the adapted space \be
\nonumber \dX^{\eta}=\left\{ \omega: \mathbb{S} \times [0,\infty)
\rightarrow \mathbb{R} \mid \omega \in L^\infty ([0, \infty); L^2),
\sup_{0<t<\infty} t^{\eta/\sigma} \| \omega \|_{\dHsig}  < \infty
\right\}, \ee where $\sigma>0$ and $\eta>0$, with norm \be \nonumber \|
\omega \|_{\dXs} = \mbox{max} \left( \sup_{0<t<\infty} \| \omega \|_{L^2},
\sup_{0<t<\infty}   t^{\eta/\sigma} \| \omega \|_{\dHsig} \right). \ee The
factor of $t^{\eta/\sigma}$ is motivated by the estimate in Lemma
\ref{lemma:HsHr} below with $s=\eta$ and $r=0$.


%
%
%

We will make use of the following elementary  result, which is proven in
the appendix:

\begin{lemma} \label{lemma:int_est}

Let  $\hat{q}>0$, $0 \leq \hat{\alpha} < 1$, and let ${\hat{\beta}}$ and
${\hat{\delta}}$ be nonnegative numbers with $0 \leq
\hat{\alpha}+{\hat{\beta}} \leq 1$ and  $0 \leq
{\hat{\beta}}+{\hat{\delta}}<1$. Then there exists a positive constant $C$
such that \be \label{time_int} \int_0^t
\frac{e^{-{\hat{q}}(t-\tau)}}{(t-\tau)^{\hat{\alpha}}}
\frac{t^{\hat{\delta}}}{\tau^{{\hat{\beta}}+{\hat{\delta}}}} \ d \tau<C,
\ee where $C$  is independent of $t \in [0,\infty).$

\end{lemma}

\subsection{Operator estimates}\label{semigroupSection}

We estimate the smoothing properties of the semigroup $e^{- t
\mathcal{L}}$ for $t>0$. First, it is clear that \be \label{HsHs} \| \emL
f \|_{\dH}  =\left( \sum_{k \in \mathbb{Z}}  |k|^{2s} e^{-2 t |k|^\sigma}
| \widehat{f} |^2 \right)^{1/2} \leq \| f \|_{\dH}. \ee
 Let $s,r \in \mathbb{R}$ with $0 \leq r < s$.
We next estimate $\| \emL f \|_{\dH} $ in terms of the $\dHr$ norm of $f$:
\begin{lemma} \label{lemma:HsHr}
Let $f \in L^2$, $t>0$ and define the positive number $p=(s-r)/\sigma$ for
$\sigma>0$. Then \be \label{HsHr} \| \emL f \|_{\dH} \leq C e^{-t/2}
\left( 1+   t^{-p}  \right) \| f \|_{\dHr}, \ee where $C$ is a positive
constant that depends only on $p$.
\end{lemma}

\noindent {\it Proof.}   Using (\ref{emLdef}), we write (after multiplying
and dividing by $|k|^{2r}$),
\begin{align}
\| \emL f \|_{\dH}^2 &= \sum_{k= -\infty}^\infty |k|^{2(s-r)} |k|^{2r}
e^{-2  t |k|^\sigma} |\widehat{f}_k|^2 \nonumber
\\
& \leq \|  |k|^{s-r}   e^{-  t |k|^\sigma} \|_{l^\infty}^2 \|  f
\|_{\dHr}^2. \label{l_infty_Hr}
\end{align}
The first factor above is now estimated.  Define $g(\kappa)=\kappa^{p}
e^{- t^{1/2} \kappa}$ for $\kappa \geq 0$ and let $C_{\sigma, s,r}=\left[
(s-r)/\sigma e \right]^{(s-r)/ \sigma}$. The maximum of $g$ occurs at
$\kappa=\kappa_m= \frac{p}{ t^{1/2}}$, at which point
$g(\kappa_m)=C_{\sigma, s,r} t^{-p/2}$.  Set $\kappa=|k|^\sigma t^{1/2}$
and substitute into the definition of $g$ to find $ |k|^{s-r}  e^{-  t
|k|^\sigma} \leq C_{\sigma, s,r} t^{-p}$, which when used in
(\ref{l_infty_Hr}) and taking square roots gives \be \label{prelim_est} \|
\emL f \|_{\dH} \leq C_{\sigma, s,r} \ t^{-p} \| f \|_{\dHr}
~~\mbox{for}~~t>0. \ee If $t>p$, the estimate above can be improved.  In
this case, the wavenumber $k_m$ at which the maximum of $g$ occurs is less
than one. Since the minimum (nonzero) wavenumber in our periodic problem
is $k=1$, for this range of $t$ the maximum of $g$ occurs at $k_m=1$ or
$\kappa_m=t^{1/2}$, at which point $g(\kappa_m)=t^{p/2}e^{-t}$. Since
$t^{p/2}e^{-t} \leq  C^{1/2}_{\sigma, s,r} e^{-t/2}$, it follows that \be
\label{improved_est} \| \emL f \|_{\dH}  \leq C^{1/2}_{\sigma, s,r}
e^{-t/2} \| f \|_{\dHr}~~\mbox{for}~~ t>p. \ee The estimate (\ref{HsHr})
follows from combining (\ref{prelim_est}) and (\ref{improved_est}).
$\blacksquare$

We also need to estimate $\| \emL f \|_{\dH}$ in terms of $\| f \|_{L^1}$
to bound some of the nonlinear terms.  We start by deriving a bound on $\|
\emL f \|_{L^2}$ in terms of $\| f \|_{L^1}$. From Plancherel's theorem
and the Young-Haussdorf inequality,
\begin{align}
\| \emL f   \|_{L^2} &\leq \max_{k} | \widehat{f}_k | \| e^{- t
\rho(\cdot)} \|_{l_2} \nonumber
\\
 &\leq \| f \|_{L^1}   \| e^{- t \rho(\cdot)} \|_{l_2} ,
\label{exp_norm}
\end{align}
where $\rho(k)=|k|^\sigma$. Note that if $f$ has zero mean then the $k=0$
term can be omitted from the $l^2$ norm in (\ref{exp_norm}). An elementary
estimate of this $l^2$  norm is proven in the appendix:
\begin{lemma} \label{lemma:sumest}
Assume $\sigma > 0$. Then there exists a constant $C>0$ that is
independent of $t$ (but which may depend on $\sigma$)  such that
\begin{align} \label{l2_est}
 \| e^{- t \rho(\cdot) } \|_{l^2}^2  \leq 1 + C e^{-t} \left( 1+ t^{-1/\sigma} \right),
\end{align}
\begin{align} \label{l2_estTake2}
 \| e^{- t \rho(\cdot) }  \|_{l^{2}_{0}}^2  \leq C e^{-t} \left( 1+ t^{-1/\sigma} \right),
\end{align}
for $t>0$.
\end{lemma}
Note that here we have introduced the set of sequences $\ell^{2}_{0},$
where a sequence $\{a_{k}\}_{k=-\infty}^{\infty}$ is in $\ell^{2}_{0}$ if
it is in $\ell^{2}$ and if also $a_{0}=0;$ in \eqref{l2_estTake2}, we use
this to mean that we simply exclude the $k=0$ term when calculating the
norm.

The above lemma applied to (\ref{exp_norm}) immediately yields  the
estimate
\begin{lemma} \label{lemma:L2L1}
Let $f \in L^2$ and  $\sigma>0$. Then $f$ is in $L^1$ and \be \label{L2L1}
 \| \emL f   \|_{L^2} \leq \left[ 1 + C e^{-t} \left( 1+ t^{-1/\sigma} \right) \right]^{1/2}   \| f \|_{L^1},
\ee for $t>0$, where $C$ is a constant that is independent of $t$. If $f$
has zero mean, then the first  $1$ in (\ref{L2L1}) can be omitted, per the
comments following (\ref{exp_norm}).
\end{lemma}

We now use the above result to  estimate $\| \emL f \|_{\dH}$ in terms of
$\| f \|_{L^1}$.  We first apply (\ref{HsHr}) to the function $e^{- t
\mathcal{L}/2} f$ with operator $e^{- t \mathcal{L}/2}$ and $r=0$, $s>0$
to find
\[
\| \emL f \|_{\dH} \leq C e^{-t/4} \left(1 + \ ( t/2)^{-s/\sigma} \right)
\|  e^{- t \mathcal{L}/2} f \|_{L^2}.
\]
We next use  (\ref{L2L1}) to bound the $L^2$-norm above in terms of the
$L^1$-norm to  obtain the following lemma.
\begin{lemma}  \label{lemma:HsL1}
Let $s>0$. Under the same conditions as in Lemma \ref{lemma:L2L1} we have
\be \nonumber
 \| \emL f   \|_{\dH} \leq   C e^{-t/4}  \left( 1+t^{-(2s+1)/2 \sigma} \right)  \| f \|_{L^1} ,
\ee where $C$ is a constant independent of $t$.
\end{lemma}

In our global existence proof for small $L^2$ data in Section
\ref{mildExistenceSection} below, we will make use of the following
estimate on the Sobolev norm of a product of two functions, which is a
straightforward generalization of an exercise in \cite{Folland}.
\begin{lemma}  \label{Banach_norm}
Let $s>1/2$  and $m \in [0,s]$.  Let $f \in \dHm$ and $g \in \dH$ be
given. Then $fg \in \dHm$ and $\| fg \|_{\dHm} \leq c \| f \|_{\dHm} \| g
\|_{\dH}$.
\end{lemma}

\subsection{Global existence for small data in $L^2$}
\label{mildExistenceSection}

We construct solutions of the initial value problem  for
(\ref{mainEquation2}) by demonstrating the existence of a fixed point of
the operator $\mathcal{T}$  in (\ref{almostOurDuhamel}). The main result
is
\begin{theorem} \label{thm_main}
Let $ \omega_0  \in L^2$ and $\sigma>1$.  Let $\eta_{m}= \min(1,
\sigma-1/2).$  There exists $\epsilon>0$ small enough such that if $\|
\omega_0 \|_{L^2} < \epsilon$, then the initial value problem for
(\ref{mainEquation2}) with initial data $\omega_0$  has a unique solution
$\omega$ in $\dX^{\eta}$ for $1/2<\eta<\eta_{m}$.
\end{theorem}

\noindent{\em Remark.} Theorem \ref{thm_main} gives solutions in
$H^{\eta}$ at positive times, with $\eta>0,$ starting from $L^{2}$ initial
data.  As is usually the case for parabolic evolutions, this gain of
regularity can be bootstrapped to find that solutions are actually
$C^{\infty}$ at positive times.  We expect more than this, as we expect
solutions to in fact be analytic at positive times, as was demonstrated
for the solutions of Theorem \ref{thm_ambrose}.  We do not include a proof
of analyticity the solutions of Theorem \ref{thm_main}, but we expect that
the corresponding argument from \cite{AmbroseMazzucato}, which itself
followed the argument
of \cite{GrujicKukavica}, would be effective.\\

\noindent {\it Proof of Theorem \ref{thm_main}.}  We first show that
$\mathcal{T}: \dXs  \rightarrow \dXs$. Througout the proof, we employ the
notation $\lesssim$ to denote  $\leq C$ with $C>0$ independent of $\omega,
~\omega_0,$ and $t$.

Decompose the map $\mathcal{T}$ in (\ref{almostOurDuhamel})  into its
linear part $\emL \omega_0$, which
is called the `trend,' and the nonlinear part  $ \int_0^t e^{- (t-\tau)
\mathcal{L}} \left[ - a (u \omega)_x +(1+a)  \omega \mathcal{H}(\omega)
\right] (\cdot,\tau) \ d \tau$, which is called the `fluctuation.' The
trend is bounded as follows. First use (\ref{HsHs}) to see that $\| \emL
\omega_0 \|_{L^2} \lesssim  \| \omega_0 \|_{L^2}$, then apply
(\ref{prelim_est})   with $s=\eta$ and $r=0$  to find  $t^{\eta/\sigma} \|
\emL \omega_0 \|_{\dHsig}  \lesssim  \| \omega_0 \|_{L^2}$.  It
immediately follows that $\emL \omega_0\in \dXs$ with $ \| \emL \omega_0
\|_{\dXs} \lesssim \| \omega_0 \|_{L^2}$.

We next bound the norm of the  fluctuation, with the terms $-a(u
\omega)_x$ and $(1+a) \omega \mathcal{H} \omega$  in
(\ref{almostOurDuhamel}) treated separately. First consider the $L^2$ norm
of the contribution from $-a(u \omega)_x$. Since derivatives and
$\mathcal{L}$ commute as Fourier multipliers on the circle, we have \be
\nonumber \left\|  \int_0^t e^{- (t-\tau) \mathcal{L}} \left[-a (u
\omega)_x \right]  \ d \tau \right\|_{L^2}
 \lesssim \int_0^t \| e^{- (t-\tau) \mathcal{L}}   (u \omega) \|_{\dot{H}^1 } \ d \tau.
\ee The right hand side is bounded  by applying Lemma \ref{lemma:HsHr}
with $s=1$ and $r=\eta$  (which requires $0 \leq \eta<1$) followed by
Lemma \ref{Banach_norm} (which further requires $\eta>1/2$)
 to obtain
\begin{align}
\int_0^t \| e^{- (t-\tau) \mathcal{L}}   u \omega \|_{\dot{H}^1 } \ d
\tau.
&\lesssim \int_0^t   e^{- (t-\tau)/2} \left[1+  \left( t-\tau \right)  ^{-\frac{1-\eta}{\sigma}} \right]  \| u \omega \|_{\dHsig} \ d \tau  \nonumber \\
 &\lesssim \int_0^t   e^{- (t-\tau)/2} \left[1+  \left( t-\tau \right) ^{-\frac{1-\eta}{\sigma}} \right]  \| u \|_{\dHsig} \|  \omega \|_{\dHsig} \ d \tau  \nonumber \\
&\lesssim \int_0^t   e^{- (t-\tau)/2} \left[1+  \left( t-\tau \right) ^{-\frac{1-\eta}{\sigma}} \right]   \| \omega \|_{L^2} \frac{ \tau^\frac{\eta}{\sigma} \|  \omega \|_{\dHsig}}{\tau^\frac{\eta}{\sigma} }\ d \tau  \nonumber \\
&\lesssim \left(   \int_0^t  e^{-(t - \tau)/2}  \frac{ \left[ 1+ \left(
t-\tau \right)^{-\frac{1-\eta}{\sigma}} \right]
}{\tau^\frac{\eta}{\sigma}}
 \ d \tau \right) \| \omega \|_{\dXs}^2 \nonumber \\
& \lesssim \| \omega \|_{\dXs}^2. \label{fluctuation_L2}
\end{align}
In the above estimate  we have used   $\| u \|_{\dHsig} \lesssim \| \omega
\|_{L^2}$, which holds for $\eta<1$.   The integral in the second-to-last
inequality is bounded for $\sigma  \geq 1$ and the assumed range of $\eta$
by applying Lemma \ref{lemma:int_est} with ${\hat{q}}=1/2$,
$\hat{\alpha}=(1-\eta)/\sigma$, ${\hat{\beta}}=\eta/\sigma$, and
${\hat{\delta}}=0$.

We next use Lemma \ref{lemma:HsHr} (with $s=\eta+1$ and $r=\eta$) and
Lemma \ref{Banach_norm} to estimate the $\dHsig$ norm for $1/2<\eta<1$:

\begin{align}
\left\|  \int_0^t e^{- (t-\tau) \mathcal{L}} \left[-a (u \omega)_x \right]
\ d \tau \right\|_{\dHsig}
 &\lesssim \int_0^t \| e^{- (t-\tau) \mathcal{L}}   (u \omega) \|_{\dot{H}^{\eta+1 }} \ d \tau \nonumber \\
&\lesssim \int_0^t   e^{- (t-\tau)/2} \left[1+  \left( t-\tau \right)  ^{-\frac{1}{\sigma}} \right]  \| u \omega \|_{\dHsig} \ d \tau  \nonumber \\
 &\lesssim \int_0^t   e^{- (t-\tau)/2} \left[1+  \left( t-\tau \right) ^{-\frac{1}{\sigma}} \right]  \| u \|_{\dHsig} \|  \omega \|_{\dHsig} \ d \tau  \nonumber \\
&\lesssim \frac{1}{t^\frac{\eta}{\sigma}} \int_0^t   e^{- (t-\tau)/2} \left[1+  \left( t-\tau \right) ^{-\frac{1}{\sigma}} \right]   \| \omega \|_{L^2} \frac{ (t \tau)^\frac{\eta}{\sigma} \|  \omega \|_{\dHsig}}{\tau^\frac{\eta}{\sigma} }\ d \tau  \nonumber \\
&\lesssim \frac{1}{t^\frac{\eta}{\sigma}}  \left(   \int_0^t  e^{-(t -
\tau)/2}  \frac{ \left[ 1+ \left( t-\tau \right)^{-\frac{1}{\sigma}}
\right] t^\frac{\eta}{\sigma} }{\tau^\frac{\eta}{\sigma}}
 \ d \tau \right) \| \omega \|_{\dXs}^2 \nonumber \\
& \lesssim \frac{1}{t^\frac{\eta}{\sigma}}  \| \omega \|_{\dXs}^2.
\end{align}
In the above estimate,   we have again used   $\| u \|_{\dHsig} \lesssim
\| \omega \|_{L^2}$.   The integral in the second-to-last inequality is
bounded for $\sigma  > 1$ and the assumed range of $\eta$   by applying
Lemma \ref{lemma:int_est} with ${\hat{q}}=1/2$,  $\hat{\alpha}=1/\sigma$,
${\hat{\beta}}=0$, and ${\hat{\delta}}=\eta/\sigma$.

A different method is required to  bound  the  fluctuation associated with
the term $(1+a) \omega \mathcal{H} \omega$  in (\ref{almostOurDuhamel}).
We first bound the $L^2$ norm of this fluctuation by applying Lemma
\ref{lemma:L2L1} with the first 1 in (\ref{L2L1}) omitted (since the
integrand has zero mean) to obtain
\begin{align}
\left\| \int_0^t e^{- (t-\tau) \mathcal{L}} (1+a) \omega \mathcal{H}
\omega  \ d \tau \right\|_{L^2}
&\lesssim \int_0^t   e^{-\frac{t - \tau}{2}} \left[1+  \left( t-\tau \right) ^{-\frac{1}{2 \sigma}} \right]    \| \omega \mathcal{H} \omega \|_{L^1}  \ d \tau  \nonumber \\
&\lesssim   \int_0^t  e^{-\frac{t - \tau}{2}}  \left[ 1 +   \left( t-\tau \right) ^{-\frac{1}{2 \sigma}} \right]   \| \omega \|_{L^2}  \| \mathcal{H} \omega \|_{L^2} \ d \tau  \nonumber \\
&\lesssim \int_0^t  e^{-\frac{t - \tau}{2}}  \left[1+  \left( t-\tau \right)^{-\frac{1}{2 \sigma}} \right]  \| \omega \|^2_{L^2}\ d \tau  \nonumber \\
& \lesssim  \| \omega \|^2_{\dXs},
\end{align}
where we have used H\"{o}lder's inequality, $\| \mathcal{H} \omega
\|_{L^2} = \| \omega \|_{L^2}$  and   Lemma \ref{lemma:int_est} with
$\hat{\alpha}=1/2 \sigma$ and ${\hat{\beta}}={\hat{\delta}}=0$ to bound
the integral.

We next use Lemma \ref{lemma:HsL1} to similarly bound the $\dHsig$ norm of
this fluctuation:

\begin{align}
\left\| \int_0^t e^{- (t-\tau) \mathcal{L}} (1+a) \omega \mathcal{H}
\omega  \ d \tau \right\|_{\dHsig}
&\lesssim \int_0^t   e^{-\frac{t - \tau}{4}} \left[1+  \left( t-\tau \right) ^{-\frac{2 \eta +1}{2 \sigma}} \right]    \| \omega \mathcal{H} \omega \|_{L^1}  \ d \tau  \nonumber \\
&\lesssim   \int_0^t  e^{-\frac{t - \tau}{4}}  \left[ 1 +   \left( t-\tau \right) ^{-\frac{2 \eta + 1}{2 \sigma}} \right]   \| \omega \|_{L^2}  \| \mathcal{H} \omega \|_{L^2} \ d \tau  \nonumber \\
&\lesssim \frac{1}{t^\frac{\eta}{\sigma}}  \int_0^t  e^{-\frac{t - \tau}{4}}  \left[1+  \left( t-\tau \right)^{-\frac{2 \eta + 1}{2 \sigma}} \right]  \| \omega \|_{L^2}  \frac{(t \tau)^{\eta/\sigma}  \| \omega \|_{\dHsig}}{\tau^{\eta/\sigma}} \ d \tau  \nonumber \\
&\lesssim \frac{1}{t^\frac{\eta}{\sigma}}   \left(   \int_0^t  e^{-\frac{t
- \tau}{4}}  \frac{ \left[ 1+ \left( t-\tau \right)^{-\frac{2 \eta +1}{2
\sigma}} \right] }{\tau^{\frac{\eta}{\sigma}}}  t^\frac{\eta}{\sigma}
 \ d \tau \right) \| \omega \|_{\dXs}^2 \nonumber \\
& \lesssim  \frac{1}{t^\frac{\eta}{\sigma}} \| \omega \|_{\dXs}^2.
\label{fluctuation_Hsig_2}
\end{align}
In the above estimate,   we have  used   $\| \mathcal{H} \omega \|_{L^2} =
\| \omega \|_{L^2} \lesssim \| \omega \|_{\dHsig}$ (since $\omega$ has
zero mean).   The integral in the second-to-last inequality is bounded for
$\sigma>1$ and $\eta<\sigma-1/2$   by applying Lemma \ref{lemma:int_est}
with $\hat{\alpha}=(2 \eta + 1)/2 \sigma$, ${\hat{\beta}}=0$, and
${\hat{\delta}}=\eta/\sigma$.

Combining the bound on the trend with   (\ref{fluctuation_L2}) and
(\ref{fluctuation_Hsig_2}) yields: \be \label{T_est} \|
\mathcal{T}(\omega) \|_{\dXs} \leq A \left( \| \omega_0  \|_{L^2} + \|
\omega \|^2_{\dXs} \right), \ee for some constant $A>0$. We similarly may
establish a Lipschitz estimate on $\mathcal{T}$ in $\dXs$:
\begin{multline}
\| \mathcal{T}(\omega_1)-\mathcal{T}(\omega_2) \|_{\dXs}  =
\\
 \left\| \int_0^t e^{-(t - \tau) \mathcal{L}} \left\{ -a \left( u_1 \omega_{1} - u_2 \omega_{2} \right)_x +(1+a) \left( \omega_1 \mathcal{H} \omega_1- \omega_2 \mathcal{H} \omega_2 \right) \right\} \ d \tau \right\|_{\dXs}
\\
 =\bigg\|
\int_0^t e^{-(t - \tau) \mathcal{L}} \bigg\{ -a \left[ u_1
(\omega_{1}-\omega_{2} ) +  (u_1- u_2) \omega_{2} \right]_x
\\
 +   (1+a) \left[ \omega_1 (\mathcal{H} \omega_1-\mathcal{H} \omega_2 ) + (\omega_1-\omega_2) \mathcal{H} \omega_2 \right] \bigg\} \ d \tau
 {\bigg \|}_{\dXs}
 \\
\leq A \left( \| \omega_1 \|_{\dXs} +  \| \omega_2 \|_{\dXs}  \right) \|
\omega_1 - \omega_2  \|_{\dXs}  \label{Lip_est},\nonumber
\end{multline}
where we have repeated the analysis leading to
(\ref{fluctuation_L2})-(\ref{fluctuation_Hsig_2}) to obtain the last
inequality, and $A$ is the constant in (\ref{T_est}).

Let $\mathcal{B}_M$ denote the ball $\left\{ \omega: \| \omega \|_{\dXs} <
M \right\}$, and set $\widetilde{M} = \| \omega_0 \|_{L^2}$. We will
determine $M$ and $\widetilde{M}$ so that  $\mathcal{T}$ is a contraction
on $\mathcal{B}_M$.   From (\ref{T_est}), $\mathcal{T}$ will be a mapping
from $\mathcal{B}_M$ into $\mathcal{B}_M$ if $A \widetilde{M} + A M^2<M$,
which can be arranged by choosing  $M<1/(2A)$ and $\widetilde{M}< M/(2A)$.
$\mathcal{T}$ is  automatically a contraction on $\mathcal{B}_M$ under
these conditions on ${M}, \widetilde{M},$ since
\[
\| \mathcal{T}(\omega_1)-\mathcal{T}(\omega_2) \|_{\dXs} \leq 2 A M \|
\omega_1 - \omega_2  \|_{\dXs},~~\omega_1, \ \omega_2 \in \mathcal{B}_M.
\]
Thus by the Contraction Mapping Theorem,  there is a unique fixed point
$\omega$ of the map $\mathcal{T}$ in $\mathcal{B}_M$.  By a standard
continuation argument, the solution is unique in $\dXs$.  $\blacksquare$

\section{Exact solutions for the problem on the real line} \label{sec:real_line}

We now consider the Constantin-Lax-Majda problem (\ref{CLM})  on the real
line $x \in \mathbb{R}$. We derive several new analytical solutions and
revisit the exact solution of Schochet \cite{Schochet}. These solutions
exhibit self-similar finite-time singularity formation from arbitrarily
small data, in contrast to the periodic problem. In this section the
viscosity parameter $\nu$ is mostly retained so that we may compare
analytical solutions for $\nu>0$ with inviscid solutions derived in
\cite{Lushnikov_Silantyev_Siegel}.

\subsection{Schochet's solution for $a=0$ and $\sigma=2$}
\label{sec:Schochet} Schochet \cite{Schochet} constructs a solution to
(\ref{CLM})  in the case $a=0$ and $\sigma=2$  by the method of pole
dynamics.

 To describe his solution, introduce the  operator $\mathbb{P}_+$ which projects onto upper analytic function space, i.e., $\mathbb{P}_+ f=f_+$.   Apply $\mathbb{P}_+$ to  (\ref{CLM}) with $a=0$ to obtain
\be \omega_{+ t} = -i \omega_+^2 + \nu \omega_{+xx},  \label{CLM_+} \ee
where $x$ is now considered  complex. Since $\omega$ is real for $x \in
\mathbb{R}$,  its lower analytic component satisfies
$\omega_-(x,t)=\overline{\omega_{+}(\bar{x},t)}$ for $x \in \mathbb{C}^-$.
Note that for $x \in \mathbb{R}$, $\omega=\omega_+ + \omega_-=
2Re[\omega_+]$. If an upper analytic function $\omega_+$ satisfies
(\ref{CLM_+}) and vanishes at infinity, then $2Re[\omega_+]$ satisfies
(\ref{CLM}).

Schochet looks for solutions  of the form (using his notation) \be
\label{Schochet_soln} \omega_+(x,t)=\frac{1}{2}
\left\{\frac{A(t)}{x-x_1(t)} + \frac{B(t)}{[x-x_1(t)]^2} +
\frac{C(t)}{x-x_2(t)} + \frac{D(t)}{[x-x_2(t)]^2} \right\}. \ee
Substituting  into (\ref{CLM_+}) and equating like-power poles yields
\begin{align}
A(t)&=-K_\pm \nu i/([x_1(0) - x_2(0)]^2 -\frac{5}{3} K_\pm \nu t)^{1/2} \label{A_eqn} \\
B(t)&=-12 \nu i,~~C(t)=-A(t),~~D(t)=B(t), \label{B_eqn}
\end{align}
in which $K_\pm=24(3 \pm \sqrt{6})$  (correcting the value of  $K_\pm=12(6
\pm \sqrt{6})$ given in \cite{Schochet})  and
\begin{align}
x_1(t)&=\frac{1}{2} \left[x_1(0)+x_2(0) +\left( \left[x_1(0)-x_2(0) \right]^2- \frac{5}{3} K_\pm \nu t \right)^{1/2} \right], \label{x1_eqn}\\
x_2(t)&=\frac{1}{2} \left[x_1(0)+x_2(0) -\left( \left[x_1(0)-x_2(0)
\right]^2- \frac{5}{3} K_\pm \nu t \right)^{1/2} \right]. \label{x2_eqn}
\end{align}
Here $x_1(0), x_2(0) \in \mathbb{C}^-$ and the sign of $K_\pm$ can be
chosen arbitrariy.

As long as $x_1(t)$ and $x_2(t)$ both remain in the lower half-plane, the
real part of (\ref{Schochet_soln}) yields a smooth (analytic) solution to
(\ref{CLM}) for $x \in \mathbb{R}$. This smooth solution has finite
kinetic energy
\be \nonumber
E_K = \int_{-\infty}^{\infty} u^2 (x) \ dx. \ee However, Schochet shows
that  for all $x_1(0)$ and $x_2(0)$ in the lower half-plane and either
choice of sign in $K_\pm$, the solution blows up in finite time. His
argument is based on adding and subtracting (\ref{x1_eqn}) and
(\ref{x2_eqn}) to obtain
\begin{align}
x_1(t)+x_2(t)&=x_1(0)+x_2(0), \label{add}\\
[x_1(t)-x_2(t)]^2 &= [x_1(0)+x_2(0)]^2 - \frac{5}{3} K_\pm \nu t.
\label{subtract}
\end{align}
Let $x_j(t)=\xi_j(t)+i \eta_j(t)$ for $j=1,2$.  Then by (\ref{add}),
$\eta_1(t)+\eta_2(t)=constant$, and the real part of (\ref{subtract})
implies that $|\eta_1(t)-\eta_2(t)| \rightarrow \infty$ as $t \rightarrow
\infty$. It follows that one of $\eta_1(t)$ or $\eta_2(t)$ must cross zero
in finite time, at which point the solution blows up.

The solution (\ref{Schochet_soln}) can be made to have arbitrarily small
initial data in either the $L^2$ or $L^\infty$ norm by taking Im$~x_j(0)
\ll 0$ for $j=1,2$. Therefore, it providess an example of finite-time blow
up starting from arbitrarily small data for the problem on the real line.

\subsubsection{Self-similar form of Schochet's solution}

\begin{figure}[ht!]
\centering
\includegraphics[scale=0.75]{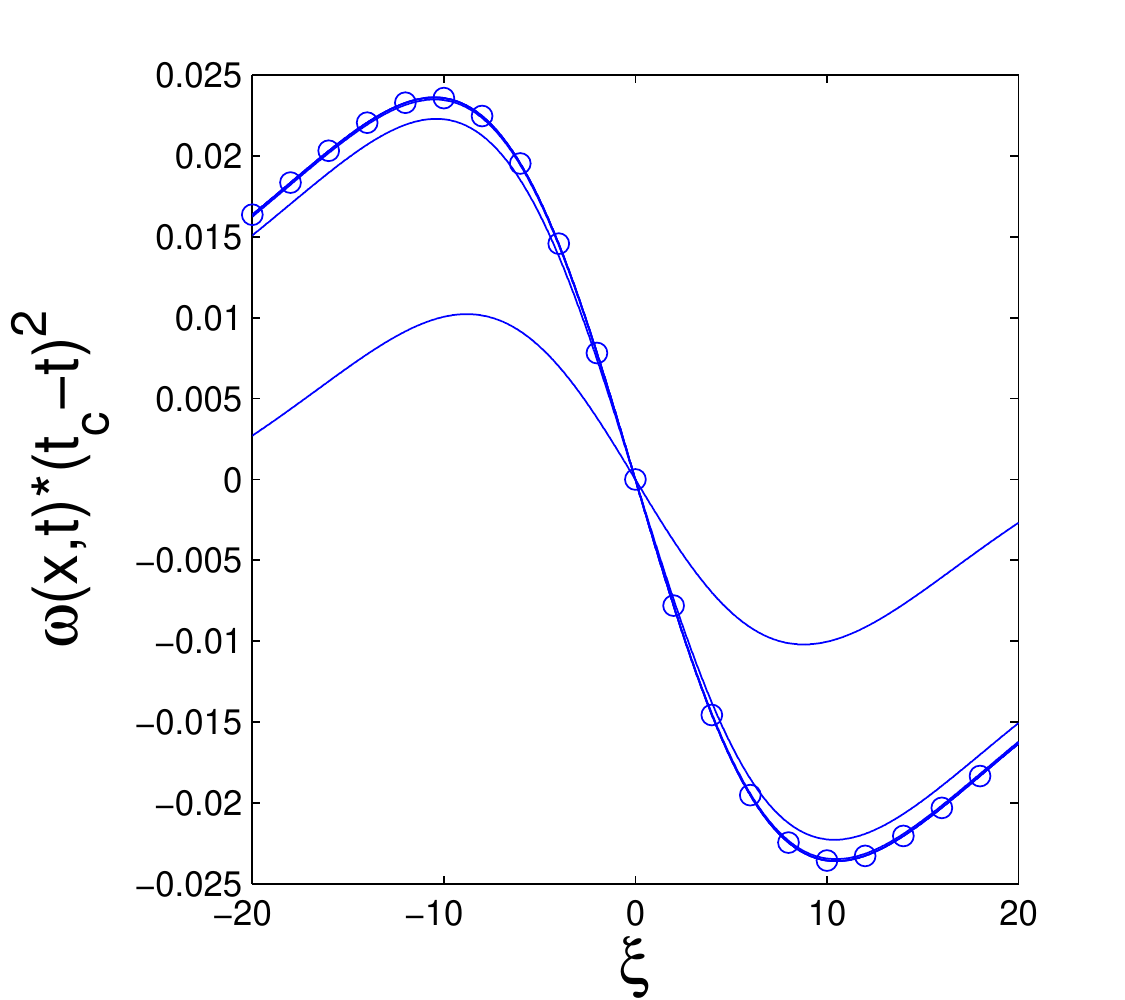}
\caption{Exact solution of Schochet for $x_1(0)=-i$, $x_2(0)=-2i$,
$\nu=1$, and $K_+=24(3+\sqrt{6})$ plotted using similarity variables
$\omega(x,t)*(t_c-t)^2$ versus $\xi$, for $t_c-t=10^{-k}, k=2, ..., 7$
(solid cuves). The asymptotic similarity solution (\ref{eqn:selfsim}) is
shown by open circles.    \label{figschochet} }
\end{figure}

Schochet's exact solution gives self-similar blow up for any initial data.
For example,  consider his solution with initial data $x_1(0)$ and
$x_2(0)$ on the negative imaginary axis, $0>Im[x_1(0)]>Im[x_2(0)]$. Then
$\omega$ is odd about $x=0$. It is easy to see that the solution for
$\omega(x,t)$ blows up at time
\begin{equation}\label{TcSchochet}
  t_c=-\frac{12}{5} \frac{x_1(0) x_2(0)}{K_\pm \nu}
\end{equation}
and that the blow up  is asymptotically self-similar in a space-time
neighborhood of $t=t_c$ and $x=0$, i.e., as $x \rightarrow 0$ and $t
\rightarrow t_c$
\begin{equation} \label{eqn:selfsim}
\omega(x,t) \simeq - \frac{24 \tilde{v}}{(t_c-t)^2}
\frac{\xi}{(\xi^2+\tilde{v}^2)^2} + O(t_c-t)^{-1},
\end{equation}
where the similarity variable $\xi$ and $\tilde{v}$ are given by
\[
\xi=\frac{x}{t_c-t},~~~\tilde{v}=-\frac{5 i}{12}  \frac{K_\pm
\nu}{x_1(0)+x_2(0)}.
\]
Figure \ref{figschochet}  shows the exact time-dependent solution for
$\omega$  with initial singularity positions $x_1(0)=-i$,$x_2(0)=-2i$ and
$\nu=1$. The solution is plotted  using similarity variables
$\omega(x,t)*(t_c-t)^2$ versus $\xi=x/(t_c-t)$  at the six times
$t_c-t=10^{-k}, k=2, ..., 7$, and  is found to approach a single universal
profile. Indeed only three separate profiles are distinguishable, with the
four curves  for $k=4$ to $7$ that are closest to $t_c$  all plotting on
top of each other. The open circles show the asymptotic self-similar
profile (\ref{eqn:selfsim}) which is approached by  the time-dependent
solution as $t \rightarrow t_c$.



\subsection{Exact  solution  for $a=1/2$ and $\sigma=1$} \label{a=1/2_semi_analytical}
When $a=1/2$ and $\sigma=1$ a new solution to (\ref{CLM})  is found using
the method of  pole dynamics. Following the analysis of
\cite{Lushnikov_Silantyev_Siegel} in the inviscid case $\nu=0$,  we look
for a solution of the form \be \label{double_pole} \omega(x,t)=i \ \omt(t)
\left( \frac{1}{[x - x_0-i v_c(t)]^2} - \frac{1}{[x - x_0+i v_c(t)]^2}
\right), \ee for which
\begin{align}
u(x,t) &=  \omt(t) \left( \frac{1}{x - x_0-i v_c(t)} +  \frac{1}{x - x_0+i v_c(t)} \right), \\
\mbox{and}~~\mathcal{H} \omega(x, t) &=  \ -\omt(t) \left( \frac{1}{[x -
x_0-i v_c(t)]^2} +  \frac{1}{[x - x_0+i v_c(t)]^2} \right).
\end{align}
Here $v_c$, $\omt$ and $x_0$ are real with $v_c>0$ and $\omt \neq 0$. The
vorticity (\ref{double_pole}) is analytic in a strip  $|Im~ x| < v_c(t)$
in the complex plane and has  double poles at $x-x_0 = \pm i v_c(t)$. The
pure imaginary amplitude $i \ \omt(t)$ implies that $\omega$ is real and
odd for $x \in \mathbb{R}$.

We substitute the ansatz (\ref{double_pole}) into (\ref{CLM}) (or
equivalently the upper analytic component $\omega_+$ of
(\ref{double_pole}) into the analog of (\ref{CLM_+}) for $a=1/2$ and
$\sigma=1$) and equate like-power poles.  Note that the leading order
$1/[x-x_0 \pm i v_c(t)]^4$ poles cancel out when $a=1/2$, which motivates
that choice for $a$ (other choices of $a$ are not consistent with a pole
dynamics solution of the form (\ref{double_pole})).  After multiplication
by $(x-x_0 +iv_c) \cdot (x-x_0  - iv_c)$ we obtain an equation which has
only single and double poles with spatially independent coefficients.
Setting the coefficients of the double poles  to zero gives \be
\label{v_c_ODE} \frac{d v_c(t)}{dt} = -\left(\frac{\omt(t)}{4 v_c(t)} -
\nu \right). \ee with initial data $v_c(0)>0$. Setting the coefficients of
the single poles to zero and using (\ref{v_c_ODE}) to eliminate $v'_c(t)$
gives \be \label{omtil_ODE} \frac{d \omt(t)}{dt}  =\frac{\omt^2(t)}{4
v_c^2(t)}. \ee with initial data $\omt(0)$. Note that (\ref{v_c_ODE}),
(\ref{omtil_ODE}) reduce to the equations derived in
\cite{Lushnikov_Silantyev_Siegel} for the inviscid case when $\nu=0$. It
is easily verified that (\ref{double_pole})-(\ref{omtil_ODE}) provide an
exact solution of the problem (\ref{CLM}) on the real line. In the
following we set $\nu=1$, which as noted earlier is equivalent to
rescaling $\omega$ and $t$.

It is instructive to define $\Omega=\omt/v_c$ and rewrite the system
 (\ref{v_c_ODE})-(\ref{omtil_ODE}) as
\begin{align}
\frac{d\Omega}{dt}&=\frac{\Omega}{v_c} \left( \frac{\Omega}{2} - 1\right), \label{AODE1} \\
\frac{d v_c}{dt} &= 1 - \frac{\Omega}{4}. \label{AODE2}
\end{align}
  Clearly, $\Omega(t)=2$ is an unstable equilibrium solution
to (\ref{AODE1}), for which $v_c(t)=(t+c)/2$ is the corresponding solution
to (\ref{AODE2}), where $c \in \mathbb{R}$ is an arbitrary constant. In
terms of the original variables, this solution is \be v_c(t)=\frac{1}{2}
(t+c), ~\omt(t)= (t+c), ~c \in \mathbb{R}~ \mbox{is a constant}.
\label{exact1} \ee A second  equilibrium solution to (\ref{AODE1}) is
$\Omega(t)=0$, for which $v_c(t)=t+c$  is the corresponding solution to
(\ref{AODE2}). This equilibrium is stable and an attractor for all
solutions with data $\Omega(0)<2$.

The above discussion implies that blow up of (\ref{double_pole}) is
determined solely by the sign of the data $\Omega(0)-2$.  More precisely,
(1)  there is  finite-time blow up with $v_c(t) \rightarrow 0$ when
$\Omega(0)>2$, and (2) the solution is analytic  and  $v_c(t)$ is
increasing for all $t>0$ when $\Omega(0) \leq 2$.

Straightforward calculations  from (\ref{double_pole}) show that \be
\nonumber \| \omega(x,t) \|_{L^\infty}= \frac{3 \sqrt{3}}{4}
\frac{\Omega(t)}{v_c(t)}, \qquad \| \omega(x,t) \|_{L^2} = \sqrt{\pi}
\frac{\Omega(t)}{\sqrt{v_c(t)}}. \ee It is therefore possible to obtain
finite-time blow up starting from arbitrarily small data, as measured by
either the $L^2$ norm or $L^\infty$ norm, by taking $\Omega(0)>2$ and
$v_c(0)\gg 1$.

\begin{figure}[ht!] \label{fig_a=1/2_sigma=1_exact}
\begin{center}
\includegraphics[scale=0.60]{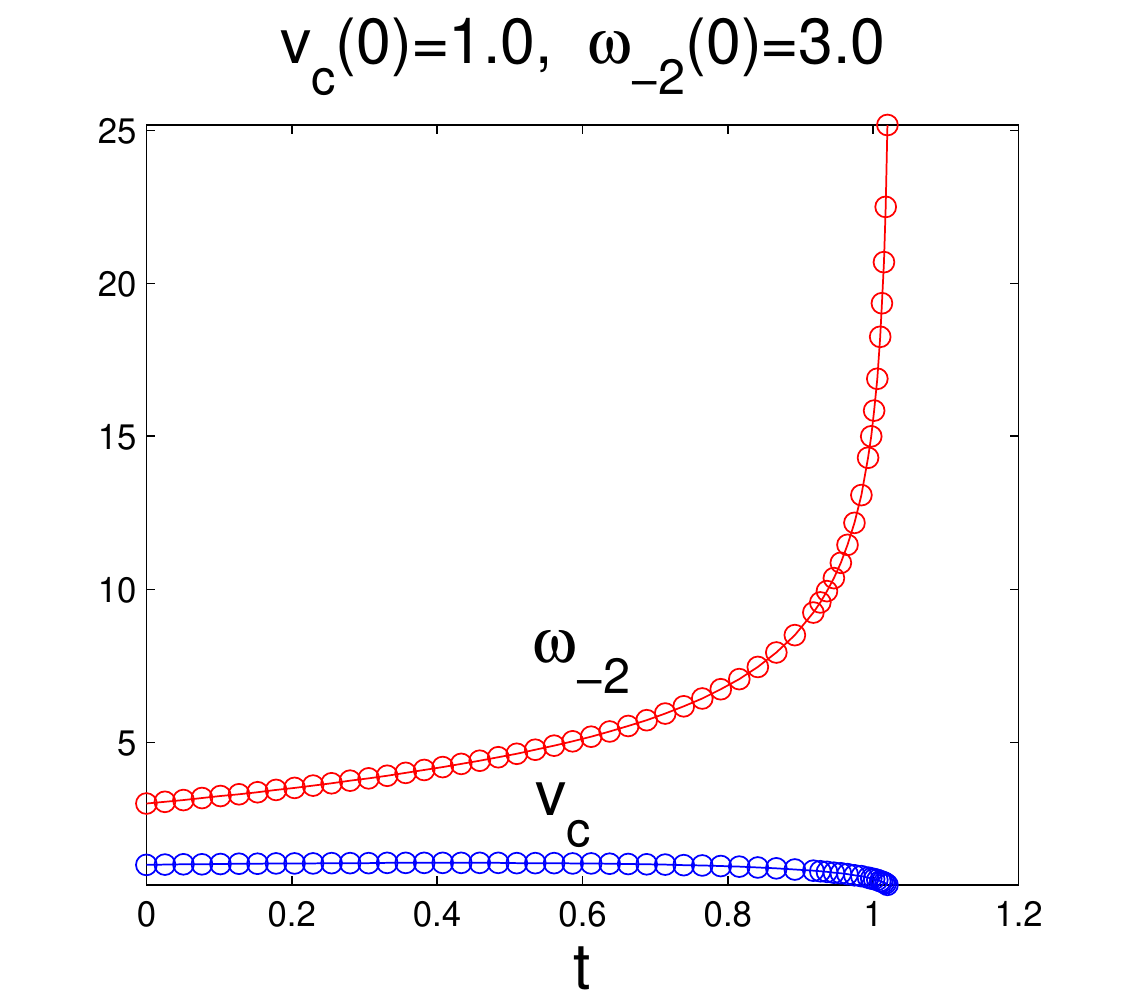}
\includegraphics[scale=0.60]{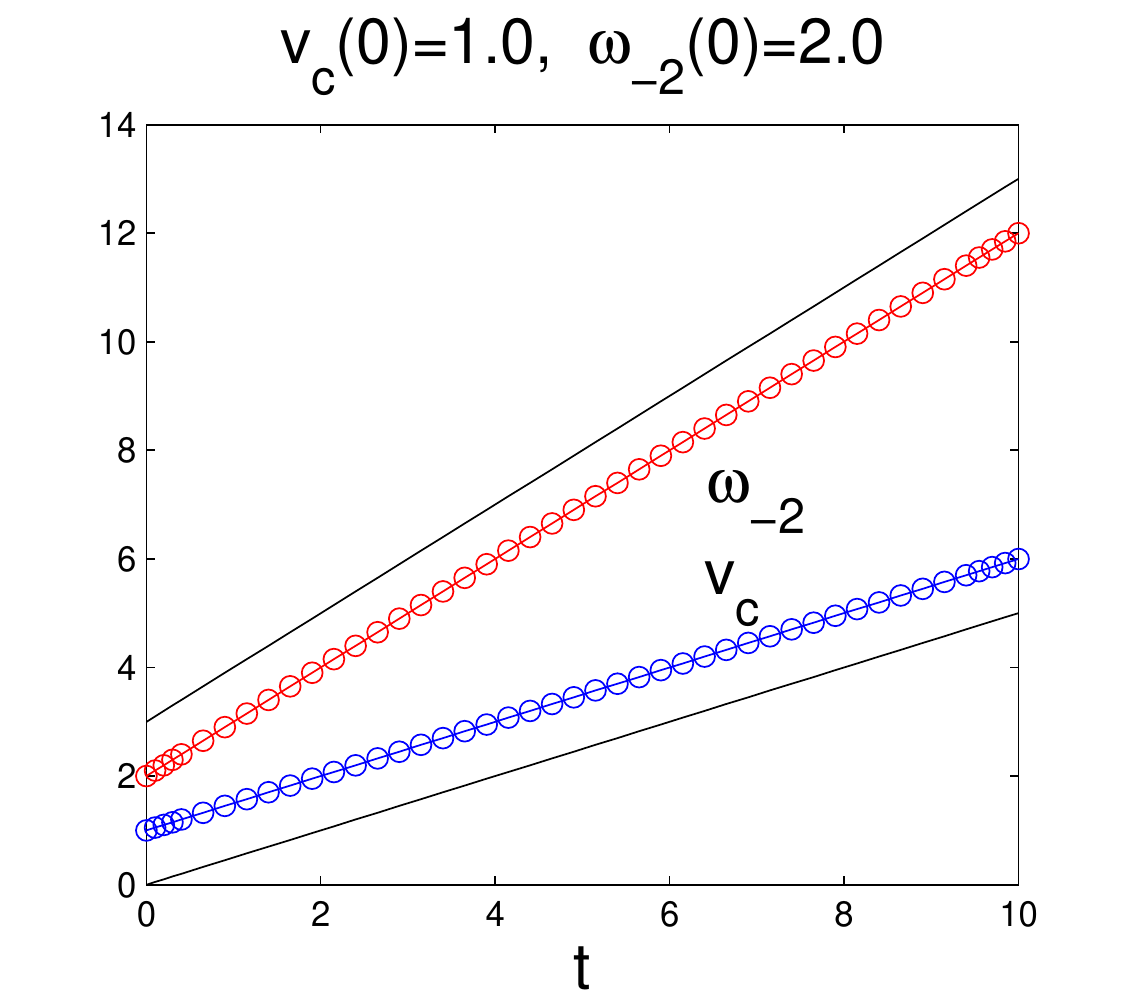}\\
\includegraphics[scale=0.60]{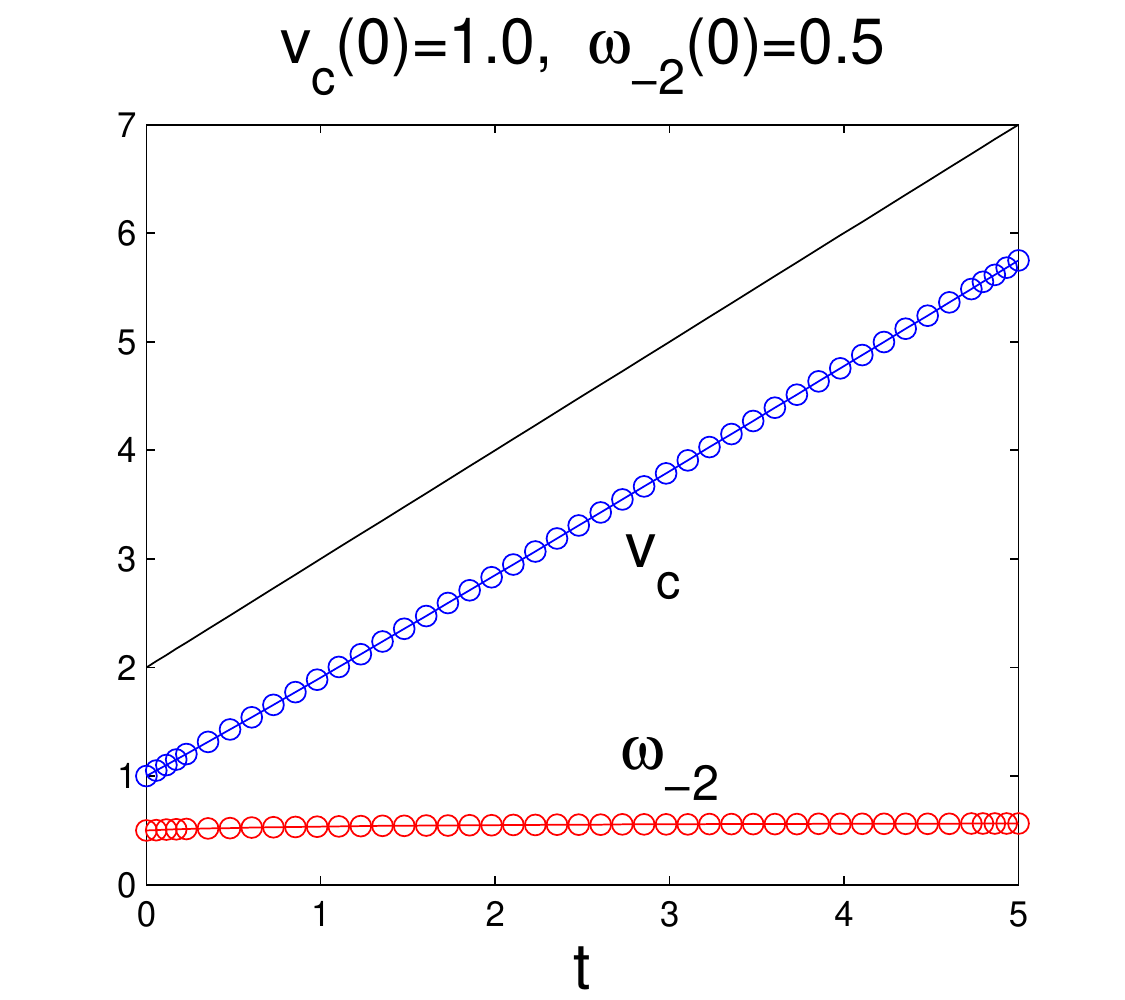}
\caption{Evolution of $\omt(t)$ and $v_c(t)$ for decreasing
$\Omega(0)=\omt(0)/v_c(0)$. (a) Finite time blow up for $\Omega(0)=3.0$,
(b) Global existence for $\Omega(0)=2.0$, solid lines have slope $1$ and
$1/2$, (c) Global existence for $\Omega(0)=0.5$, solid line has slope $1$.
\label{fig:figure1} }
\end{center}
\end{figure}

During blow up ($\Omega(0)>2$), the first term on the right-hand-side of
(\ref{v_c_ODE}) (or equivalently (\ref{AODE1})) grows rapidly over time,
and the solution asymptotically approaches \be \label{similarity_soln1}
v_c(t) = (t_c-t)^{1/3} \widetilde{v}_c,~~\omt(t)=\frac{4
\widetilde{v}_c^2}{3 (t_c-t)^{1/3}}, \ee in a space-time neighborhood of
the singularity $x \rightarrow x_0$ and $t \rightarrow t_c$, where
$\widetilde{v}_c>0$ and $t_c>0$ are two arbitrary real constants.
Numerical solutions illustrating this behavior are shown in Figure
\ref{fig:figure1}(a).   The vorticity (\ref{double_pole}) with $v_c(t)$
and $\omt(t)$ given by (\ref{similarity_soln1}) is an exact  solution of
the inviscid problem $\nu=0$  \cite{Lushnikov_Silantyev_Siegel}, and can
be written in the self-similar form \be \label{similarity_solution_omega}
\omega(x,t)=-\frac{1}{t_c-t}\frac{16 \tilde{v}_c^3 \xi}{3
(\xi^2+\tilde{v}_c^2)^2}, \ee where
\[
\xi=\frac{x-x_0}{(t_c-t)^{1/3}}.
\]

Numerical solutions to (\ref{v_c_ODE}), (\ref{omtil_ODE}) in the stable
case $\Omega(0) \leq 2$ are shown in  Figures
 \ref{fig:figure1}(b) and (c).  These are  plotted in the original variables $\omt$ and $v_c$. Figure  \ref{fig:figure1}(b) shows the solution
for $\Omega(0)=2$ (cf. (\ref{exact1})), while Figure  \ref{fig:figure1}(c)
shows it for $\Omega(0)=0.5$. In the latter, $\Omega(t) \sim 0$ and
$v_c(t)\sim t$ for $t \gg 1$, as is verified analytically from
(\ref{AODE1}), (\ref{AODE2}).  Here, the function $\omt(t)$ tends to a
constant as $t \rightarrow \infty$.

Finally, it is  noted  that we have been able to integrate (\ref{AODE1}),
(\ref{AODE2}) and obtain a solution in implicit form as
\begin{equation} \label{eq:implicit_soln_A}
\frac{\sqrt{\Omega-2}}{2 \Omega} +\frac{1}{2 \sqrt{2}} \tan^{-1}
\sqrt{\frac{\Omega-2}{2}} = c_1 t +c_2
\end{equation}
where $c_1$ and $c_2$ are constants. While it is not possible to obtain an
explicit solution for $\Omega$, the limit $\Omega \rightarrow \infty$ (or
equivalently $t \rightarrow t_c$) is easily computed with the result that
$\Omega(t) \sim (t_c-t)^{-2/3}$ in this limit. This gives from
(\ref{AODE2}) that  $v_c(t) \sim (t_c-t)^{1/3}$ and hence $\omega_{-2}(t)
\sim (t_c-t)^{-1/3}$ when $t \rightarrow t_c$. Thus, the similarity
scalings for the blow up solution  (\ref{similarity_soln1}) are recovered
from the implicit solution (\ref{eq:implicit_soln_A}).

In summary, the  analytical solutions derived here for $a=1/2$ and
$\sigma=1$ exist globally  in time and are smooth (analytic) for initial
data $\omega(x,0)$ of the form (\ref{double_pole}) with
$\Omega(0)=\omt(0)/v_c(0) \leq 2$.  When $\Omega(0)>2$, there is
finite-time blow up, and by taking $v_c(0)$ large, the blow-up can be made
to occur from arbitrarily small data  $\| \omega(\cdot, 0) \|_{L^2}$ or
$\| \omega(\cdot, 0) \|_{L^\infty}$.

\subsection{Exact solution  for $a=0$ and $\sigma=1$}
\label{sec:vCLMa0sigma1exact}

Another new solution can be found by the method of pole dynamics when
$a=0$ and $\sigma=1$. We look for a solution to (\ref{CLM}) in the form
of two poles as  %
\begin{equation}\label{omegagamma1}
\omega(x,t)=\frac{\omega_{-1}(t)}{x-\I v_c(t)}+\frac{\bar
\omega_{-1}(t)}{x+\I\bar  v_c(t)}
\end{equation}
This is easily seen to  result  in a  solution %
\begin{equation}\label{omeg1pvc0}
\omega_{-1}(t)=\omega_{-1}(0),  \qquad
v_c(t)=(\omega_{-1}(0)+\nu)t+v_c(0),
\end{equation}
where $\omega_{-1}(0)$ and $v_c(0)$ are arbitrary complex constants (in
contrast to Section \ref{a=1/2_semi_analytical}, in which the
corresponding constants are real) with $Re[v_c(0)]>0$. Note that

(a) If $\omega_{-1}(0)=-\nu$ then the solution \e{omegagamma1},
\e{omeg1pvc0} is time-independent.

(b) If $Re[\omega_{-1}(0)]>-\nu$ then the solution  \e{omegagamma1},
\e{omeg1pvc0}  exists for all $t>0$ because the poles moves away from the
real axis.

(c) If $Re[\omega_{-1}(0)]<-\nu$ then the solution  \e{omegagamma1},
\e{omeg1pvc0} exists until $t=t_c>0,$ where
\begin{equation}\label{tcdef1pole}
t_c=-\frac{Re[v_c(0)]}{Re[\omega_{-1}(0)]+\nu}
\end{equation}
is the collapse time (i.e. the time when the poles reach the real axis). Using   \e{omegagamma1}-\e{tcdef1pole} we find that at $t=t_c$, both poles impinge on the real axis with spatial location given by   %
\begin{equation}\label{x_c}
x=x_c:=\I v_c(t_c)=\frac{Im[\omega_{-1}(0)]Re[v_c(0)]}
{Re[\omega_{-1}(0)]+\nu}-Im[v_c(0)]
\end{equation}

  If  $\omega_{-1}(0)\ne-\nu$ then we can rewrite  \e{omegagamma1}-\e{x_c} in a self-similar form
\begin{equation}\label{omegaser1polesnu}
\omega(x,t)=\frac{\I}{t_c-t}\left(\frac{\xi_++\I\nu}{\xi-\xi_+}-\frac{\xi_--\I\nu}{\xi-\xi_-}\right
),
\end{equation}
where  %
\begin{equation}\nonumber
\xi_\pm=\frac{\pm\I v_c(0)-x_c}{t_c}
\end{equation}
are positions of poles in the complex plane of $\xi$ and%
\begin{equation}\label{xidef0}
\xi:=\frac{x-x_c}{t_c-t}
\end{equation}
is the self-similar variable.
 Equation \e{omegaser1polesnu} is the analog of equation (30) in  \cite{Lushnikov_Silantyev_Siegel}, which describes  self-similar blow up in the inviscid problem.  The solution \e{omegaser1polesnu} belongs to the general self-similar form (\ref{self-similar1})
%
with ${\alpha}=\beta=1.$

From (\ref{omegagamma1}) we can directly compute norms
 \begin{equation} \nonumber
 \| \omega(x,t) \|_{L^\infty} = \frac{|\omega_{-1}(t)|+|Im[\omega_{-1}(t)]|}{Re[v_c(t)]} ,
 \quad \| \omega(x,t) \|_{L^2} = \sqrt{2\pi}\frac{|\omega_{-1}(t)|}{\sqrt{Re[v_c(t)]}}.
 \end{equation}
 Both of the norms can be made arbitrarily small for the initial data of the collapsing solution \e{omegaser1polesnu} by choosing $Re[v_c(0)]$ large enough.

The solution \e{omegagamma1}, \e{omeg1pvc0} has  infinite kinetic energy
$E_K(t_c)$
  on the line $x\in \R$ for general values of the parameters $\omega_{-1}(0)$ and $v_c(0). $ An exception in which $E_K$ is finite occurs for $\omega_{-1}(0)=\bar \omega_{-1}(0)$ and $v_c(0)=\bar v_c(0)$, i.e. for purely real values of the residue $\omega_{-1}(0)$ and $v_c(0)$ in the solution \e{omegagamma1}.

We can also consider a solution with two pairs of poles as  %
\begin{equation}\label{omegagamma2pairssigma1}
\omega(x,t)=\frac{\omega_{-1,1}(t)}{x-\I v_{c,1}(t)}+ \frac{\bar
\omega_{-1,1}(t)}{x+\I\bar  v_{c,1}(t)}+\frac{\omega_{-1,2}(t)}{x-\I
v_{c,2}(t)}+ \frac{\bar \omega_{-1,2}(t)}{x+\I\bar  v_{c,2}(t)}
\end{equation}
in which the  poles are located at $x=\I v_{c,1}(t)$, $x=\I v_{c,2}(t)$
and their  complex conjugate points. Here we assume that
$Re[v_{c,1}(0)]>0$ and  $Re[v_{c,2}(0)]>0$.
Plugging \e{omegagamma2pairssigma1} into   \e{omegagamma1} and equating  the most singular terms (which are proportional to $(x-\I v_{c,1}(t))^{-2}$ and $(x-\I v_{c,2}(t))^{-2}$ at $x=\I v_{c,1}(t)$ and $x=\I v_{c,2}(t)$) results in %
\begin{equation}\label{vc1pcond}
\frac{d v_{c,1}(t)}{dt}=\nu+ \omega_{-1,1}(t)
\end{equation}
and %
\begin{equation}\label{vc2pcond}
\frac{d v_{c,2}(t)}{dt}=\nu+ \omega_{-1,2}(t).
\end{equation}

Collecting now the next most  singular terms which are proportional to  $(x-\I v_{c,1}(t))^{-1}$ and $(x-\I v_{c,2}(t))^{-1}$ at $x=\I v_{c,1}(t)$ and $x=\I v_{c,2}(t)$ results in  %
\begin{equation}\label{omm1pcond}
\frac{d \omega_{-1,1}(t)}{dt}=\frac{ 2\omega_{-1,1}(t)\omega_{-1,2}(t)}{
v_{c,1}(t)- v_{c,2}(t)}
\end{equation}
and %
\begin{equation}\label{omm1pcond2}
\frac{d\omega_{-1,2}(t)}{dt}=\frac{ 2\omega_{-1,1}(t)\omega_{-1,2}(t)}{
v_{c,2}(t)- v_{c,1}(t)}=-\frac{d \omega_{-1,1}(t)}{dt}.
\end{equation}
Substitution of  \e{omegagamma2pairssigma1}-\e{omm1pcond2}
 into  the governing equations (\ref{CLM}) reveals that they are identically  satisfied. A solution of the system \e{vc1pcond}-\e{omm1pcond2} follows from  the observation that $\frac{d\omega_{-1,2}(t)}{dt}=-\frac{d \omega_{-1,1}(t)}{dt} $  from \e{omm1pcond},\e{omm1pcond2}, so that%
\begin{equation}\label{c0def}
c_0:=\omega_{-1,1}(t)+\omega_{-1,2}(t)=\omega_{-1,1}(0)+\omega_{-1,2}(0).
\end{equation}
  Together with \e{vc1pcond},\e{vc2pcond}, this implies   that $\frac{d [v_{c,1}(t)+v_{c,2}(t)]}{dt}=2\nu+c_0, $ i.e. $v_{c,2}(t)=-v_{c,1}(t)+(2\nu+c_0)t+v_{c,1}(0)+v_{c,2}(0)$. Thus we reduce the system \e{vc1pcond}-\e{omm1pcond2} from four ordinary differential equations (ODEs) to two ODEs for $v_{c,1}(t)$ and $\omega_{-1,1}(t)$, which is easily solved. The  solution of the system \e{vc1pcond}-\e{omm1pcond2} for $c_0\neq0$ is %

\begin{equation}
\begin{split} \label{omega1c}
&\omega_{-1,1}(t)=\frac{c_0}{2}+\frac{1}{2}\frac{c_1+\frac{c_0^2
t}{v_{c,1}(0)-v_{c,2}(0)}}
{\sqrt{1 + \frac{2c_1 t}{v_{c,1}(0)-v_{c,2}(0)} + \frac{c_0^2t^2}{(v_{c,1}(0)-v_{c,2}(0))^2}}}, \\
&\omega_{-1,2}(t)=-\omega_{-1,1}(t)+c_0,\\
&v_{c,1}(t)=\frac{c_0}{2}t+\nu t+\frac{v_{c,1}(0)-v_{c,2}(0)}{2}\sqrt{1 + \frac{2c_1 t}{v_{c,1}(0)-v_{c,2}(0)} + \frac{c_0^2t^2}{(v_{c,1}(0)-v_{c,2}(0))^2}}+ \\ &+\frac{v_{c,1}(0)+v_{c,2}(0)}{2},  \\
&v_{c,2}(t)=-v_{c,1}(t)+(2\nu+c_0)t+v_{c,1}(0)+v_{c,2}(0),
\end{split}\nonumber
\end{equation}
where $c_0$ is given by Eq. \e{c0def},
$c_1:=\omega_{-1,1}(0)-\omega_{-1,2}(0)$, and we assumed a principle
branch of the square root.

For $c_0=0$, when $\omega_{-1,1}(0)=-\omega_{-1,2}(0)$, the  solution of
the system \e{vc1pcond}-\e{omm1pcond2} is
\begin{equation}\label{omega1c_c0}
\begin{split}
&\omega_{-1,1}(t)=-\omega_{-1,2}(t)=\frac{\omega_{-1,1}(0)}
{\sqrt{1+\frac{4\omega_{-1,1}(0)t}{v_{c,1}(0)-v_{c,2}(0)}}}, \\
&v_{c,1}(t)=\nu t+\frac{v_{c,1}(0)-v_{c,2}(0)}{2}\sqrt{1+\frac{4\omega_{-1,1}(0)t}{v_{c,1}(0)-v_{c,2}(0)}} + \frac{v_{c,1}(0)+v_{c,2}(0)}{2},  \\
&v_{c,2}(t)=-v_{c,1}(t)+2\nu t+v_{c,1}(0)+v_{c,2}(0).  \\
\end{split}
\end{equation}

The above solutions develop a finite-time singularity on the real line of
$x$  at $t=t_c$ provided   $Re[v_{c,1}(t_c)]=0$ or  $Re[v_{c,2}(t_c)]=0$.
By relabeling complex singularities if necessary, we can assume without
loss of generality that  $x=\I v_{c,1}(t)$  reaches the real line first
(ahead of    $x=\I v_{c,2}(t))$ thus resulting in collapse.
(It remains an open question whether it is possible to have
$v_{c,1}(t_c)= v_{c,2}(t_c)$, thus creating a higher order singularity at
$t=t_c.$)
Then for $t\to t_c-$ and in a small spatial neighborhood of $x=\I
v_{c,1}(t_c)$, the solution \e{omegagamma2pairssigma1} is dominated by
singularities at $x=\I v_{c,1}(t)$  and $x=-\I \bar v_{c,1}(t)$, so that
\e{omegagamma2pairssigma1} reduces to

\begin{equation}\label{omegagamma2pairssigma1short}
\omega(x,t)\simeq\frac{\omega_{-1,1}(t)}{x-\I v_{c,1}(t)}+ \frac{\bar
\omega_{-1,1}(t)}{x+\I\bar  v_{c,1}(t)}.
\end{equation}
Generically the singularity at $t=t_c$   located  at $x=\I v_{c,1}(t)$
hits the real line $Im(x)=0$ with a nonzero vertical velocity $\left
.\frac{d Re[v_{c,1}(t)]}{dt}\right|_{t=t_c}< 0. $ Then for $t\to t_c-$ the
solution \e{omegagamma2pairssigma1short} can be further reduced using the
Taylor series approximation  $v_{c,1}(t)=\I Im[v_{c,1}(t_c)] +
v_{c,1}'(t_c)(t-t_c)+O(t-t_c)^2$ (here $ v_{c,1}'(t_c):=\frac{d
v_{c,1}(t_c)}{dt}$) and neglecting the $O(t-t_c)^2$ term. We also assume
that $\omega_{-1,1}(t_c)\ne 0$ and  replace $\omega_{-1,1}(t)$ by
$\omega_{-1,1}(t_c)$  to obtain from \e{omegagamma2pairssigma1short}

\begin{equation}\nonumber
\omega(x,t)\simeq\frac{\omega_{-1,1}(t_c)}{x+ Im[ v_{c,1}(t_c)]-\I(t-t_c)
v_{c,1}'(t_c)}+ \frac{\bar \omega_{-1,1}(t_c)}{x+ Im[
v_{c,1}(t_c)]+\I(t-t_c)\bar  v_{c,1}'(t_c)}.
\end{equation}
This has  the self-similar form  (\ref{self-similar1}) with
$\alpha=\beta=1$.

A special situation occurs when  $Re[v_{c,1}'(t_c)]= 0. $ In that case $v_{c,1}'(t)=\I Im[v_{c,1}'(t_c)] + O(t_c-t)$, which corresponds to the pole singularity hitting  the real line of $x$ with  vanishing vertical velocity. In that case equation \e{omegagamma2pairssigma1short} turns into %

\begin{multline}\label{omegagamma2pairssigma1shortself-similar2}
\omega(x,t)\simeq\frac{\omega_{-1,1}(t_c)}{x+ Im[ v_{c,1}(t_c)]+Im[ v_{c,1}'(t_c)](t-t_c)-\frac{\I}{2}(t-t_c)^2 v_{c,1}''(t_c)}+\\
\frac{\bar \omega_{-1,1}(t_c)}{x+ Im[ v_{c,1}(t_c)]+Im[
v_{c,1}'(t_c)](t-t_c)+\frac{\I}{2}(t-t_c)^2\bar  v''_{c,1}(t_c)}.
\end{multline}
This occurs, for example, when $\omega_{-1,1}(0)=-\omega_{-1,2}(0)=K$,
where $K<0$ is a real number, and
\begin{equation}\nonumber
\frac{Re[v_{c,2}(0)]}{Re[v_{c,1}(0)]}=\left(\frac{K-\nu}{K+\nu}\right)^2,
Im[v_{c,1}(0)]=Im[v_{c,2}(0)]=0.
\end{equation}
In this case, $t_c = \frac{K^2-\nu^2}{4K\nu^2}Re[v_{c,2}(0)-v_{c,1}(0)]$.

Interestingly, the solution
\eqref{omegagamma2pairssigma1shortself-similar2} has a different
self-similar scaling than \e{omegaser1polesnu},\e{xidef0},
%
with  $\alpha=\beta=2.$
Direct numerical simulations in Section \ref{sec:num_a=0_sigma=1} show
this type of self-similar collapse is unstable to perturbations, as might
be expected.

The solution \e{omegagamma2pairssigma1},  similar to  \e{omegagamma1}, may
have arbitrarily small $L^{\infty}$ and $L^2$ norms at $t=0$ if we choose
$Re[v_{c,1}(0)]$ and $Re[v_{c,2}(0)]$ large enough.
These norms simplify in the case $Im[v_{c,1}(0)]=Im[v_{c,2}(0)]=0$, e.g.,
in which
 \begin{multline} \nonumber
 \| \omega(x,t) \|_{L^2} =\\ \sqrt{2\pi}\sqrt{\frac{|\omega_{-1,1}(t)|^2}{v_{c,1}(t)}+\frac{|\omega_{-1,2}(t)|^2}{v_{c,2}(t)} +\frac{4(Re[\omega_{-1,1}(t)] Re[\omega_{-1,2}(t)] + Im[\omega_{-1,1}(t)] Im[\omega_{-1,2}(t)]) }{v_{c,1}(t)+v_{c,2}(t)} }.
\end{multline}

\subsection{Exact solution  for $a=0$ and $\sigma=0$}
\label{sec:vCLMa0sigma0exact} Another solution can be found by the method
of pole dynamics when $a=0$ and $\sigma=0$. In this case
${\Lambda}^0(\omega)=\omega$.  We look for a solution to (\ref{CLM}) in
the form
of two simple poles as  %
\begin{equation}\label{omegagamma0}
\omega(x,t)=\frac{\omega_{-1}(t)}{x-\I v_c(t)}+\frac{\bar
\omega_{-1}(t)}{x+\I\bar  v_c(t)},
\end{equation}
where $\omega_{-1}(0)$ and $v_c(0)$ are arbitrary complex constants with
$Re[v_c(0)]>0$.

From (\ref{omegagamma0}) we directly compute norms
 \begin{equation} \nonumber
 \| \omega(x,t) \|_{L^\infty}=  \frac{|\omega_{-1}(t)|+|Im[\omega_{-1}(t)]|}{Re[v_c(t)]} ,
 \quad \| \omega(x,t) \|_{L^2} = \sqrt{2\pi}\frac{|\omega_{-1}(t)|}{\sqrt{Re[v_c(t)]}}.
 \end{equation}

Substituting \e{omegagamma0} into (\ref{CLM}) we get the following
equations:
\begin{equation}\nonumber
\frac{d\omega_{-1}(t)}{dt}=-\nu\omega_{-1}(t),  \qquad  \frac{dv_c(t)}{dt}
=\omega_{-1}(t),
\end{equation}
and their solution:
\begin{equation}\label{a0sigma0_solution}
\omega_{-1}(t)=\omega_{-1}(0) e^{-\nu t},  \qquad
v_c(t)=\frac{\omega_{-1}(0)}{\nu}(1-e^{-\nu t})+v_c(0),
\end{equation}
which for $\nu=0$ reduces to
\begin{equation}\label{a0sigma0_solutio_nu0}
\omega_{-1}(t)=\omega_{-1}(0),  \qquad v_c(t)=\omega_{-1}(0)t+v_c(0).
\end{equation}

For the case $\nu=0$ we always have a collapsing solution (even for
arbitrarily small data in $L^\infty$ and $L^2$ norms) if
$Re[\omega_{-1}(0)]<0$, since $Re[v_c(t_c)]=0$ at
$t_c=\frac{Re[v_c(0)]}{-Re[\omega_{-1}(0)]}$. This solution is equivalent
to
  (32) in  \cite{Lushnikov_Silantyev_Siegel}, which describes  self-similar blow up in the inviscid problem.

Note that for $\nu>0$ equations \e{a0sigma0_solution} indicate either
global existence of the solution \e{omegagamma0} or a collapsing solution
depending on initial values of $\omega_{-1}(0)$ and $v_c(0)$:

(a) If $Re[\omega_{-1}(0)]>-\nu Re[v_c(0)]$ then 
$Re[v_c(t)]>0$ for all $t>0$ and the solution \e{omegagamma0},
\e{a0sigma0_solution} exists for any $t>0$.

(b) If $Re[\omega_{-1}(0)]=-\nu Re[v_c(0)]$ then
$Re[v_c(t)]=Re[v_c(0)]e^{-\nu t} \rightarrow 0$ as $t\rightarrow \infty$,
and the solution  \e{omegagamma0}, \e{a0sigma0_solution}  exists for all
$t>0$ because the poles approach the real axis exponentially in time. They
approach the real line at the point
\begin{equation}\label{xcgamma0tinf}
x=x_c:=Re [iv_c(\infty)]=-Im[v_c(\infty)] =
-\frac{Im[\omega_{-1}(0)]}{\nu} - Im[v_c(0)],
\end{equation}
and $\| \omega(x,t) \|_{L^\infty}=\| \omega(x,0) \|_{L^\infty}= const>0$,
$\| \omega(x,t) \|_{L^2} \sim e^{-\nu t/2} \rightarrow 0$ as $t\rightarrow
\infty$.

(c) If $Re[\omega_{-1}(0)]<-\nu Re[v_c(0)]$ then the solution
\e{omegagamma0}, \e{a0sigma0_solution} exists until the collapse time
%
\begin{equation}\label{tcgamma0}
t_c=-\frac{1}{\nu}\ln \left(1+\frac{\nu
Re[v_c(0)]}{Re[\omega_{-1}(0)]}\right),
\end{equation}
at which the poles reach the real axis. Here $\| \omega(x,t)
\|_{L^\infty},\| \omega(x,t) \|_{L^2} \rightarrow \infty$ as $t\rightarrow
t_c$ since $\omega_{-1}(t_c)\neq0$ and $Re[v_c(t_c)]=0$. Crucially,
collapse occurs even when the   initial norm $\| \omega(x,0) \|_{L^2}$ is
made arbitrarily small by taking small $v_c(0)$. In contrast,  $\|
\omega(x,0) \|_{L^\infty}>\nu$, i.e., the $L^\infty$ norm is bounded from
below.

Using   \e{a0sigma0_solution}-\e{tcgamma0} we find that at $t=t_c$, both
poles cross the real axis at the location
\begin{equation}\nonumber
x=x_c:=\I v_c(t_c)=-Im[v_c(t_c)]=\frac{Im[\omega_{-1}(0)]}
{Re[\omega_{-1}(0)]}Re[v_c(0)]-Im[v_c(0)],
\end{equation}
with the complex velocity of the first pole being
\begin{equation}\nonumber
v_c'(t_c)=\omega_{-1}(t_c)=\omega_{-1}(0)e^{-\nu t_c}=\omega_{-1}(0)
\left(1+\frac{\nu Re[v_c(0)]}{Re[\omega_{-1}(0)]}\right).
\end{equation}
Since $v_c(t_c)=\I Im[v_c(t_c)]$,
we have that in a space-time neighborhood of the singularity $x
\rightarrow x_c$ and $t \rightarrow t_c$ the solution
\e{a0sigma0_solution} asymptotically approaches
\begin{equation} \nonumber
\omega_{-1}(t)=\omega_{-1}(t_c) + O(t_c-t), v_c(t) =\I Im[v_c(t_c)] -
\omega_{-1}(t_c) (t_c-t) + O(t_c-t)^2,
\end{equation}
and the solution \e{omegagamma0} can be written in a self-similar form
\begin{equation}\label{similarity_solution_omega_gamma0}
\omega(x,t)=\frac{\I}{t_c-t}\left(\frac{\xi_+}{\xi-\xi_+}-\frac{\xi_-}{\xi-\xi_-}\right
)+ O(1),
\end{equation}
where  %
\begin{equation}\nonumber
\xi_+=-\I\omega_{-1}(t_c),\quad \xi_-=\I\bar\omega_{-1}(t_c),
\end{equation}
are positions of poles in the complex plane of $\xi$ and%
\begin{equation} \label{xidef0_} \nonumber
\xi:=\frac{x-x_c}{t_c-t}
\end{equation}
is the self-similar variable.
 Equation \e{similarity_solution_omega_gamma0} is a viscous analog of equation (30) in  \cite{Lushnikov_Silantyev_Siegel}, which describes  self-similar blow up in the inviscid problem.  The solution \e{similarity_solution_omega_gamma0} belongs to the general self-similar form (\ref{self-similar1}) with $\alpha=\beta=1.$

The kinetic energy $E_K(t)=\int u^2(x,t)dx$ in (a), (b) scales like
$E_K(t) \sim e^{-2 \nu t}$ as $t\rightarrow \infty$, whereas in (c)
$E_K(t_c)$ is finite  for any complex values of the parameters
$\omega_{-1}(0)$ and $v_c(0)$, in contrast to the  $a=0, \sigma=1$ case.

\section{Exact solution to the periodic problem for $~~~$ $a=0$ and $\sigma=0$} \label{sec:periodic}


In this section we adapt the analysis of Section
\ref{sec:vCLMa0sigma0exact} to obtain an exact analytical  solution in the
periodic geometry. We take $a=0$, $\sigma=0$, $\nu>0$, in which case
(\ref{mainEquation2}) becomes
\begin{equation} \label{mainEquationa0sigma0_periodic}
\omega_{t}=\omega \mathcal{H}(\omega)-\nu \omega.
\end{equation}
We take initial data with zero mean value on $x \in [-\pi,\pi]$, which is
then preserved under the evolution.

Using the Hilbert transform representation \e{Hilbert} we can rewrite
\e{mainEquationa0sigma0_periodic} as:
\begin{equation} \label{mainEquationa0sigma0_periodic_lower}
\omega_{-t}=\I \omega_-^2 -\nu \omega_-,
\end{equation}
where $\omega_-$ is analytic in the lower half-plane  $\mathbb{C}^-$. We
look for a solution to \e{mainEquationa0sigma0_periodic_lower} in the form
of a single pole in $\tan(\frac{x}{2})$-space:
\begin{equation}\label{omegagamma0_periodic_lower}
\omega_-(x,t)=\omega_{-1}(t)\left[\frac{1}{\tan(\frac{x}{2})-\I v_c(t)}-
\frac{1}{-\I-\I v_c(t)}\right],
\end{equation}
where $\omega_{-1}(0)$ and $v_c(0)$ are arbitrary complex constants with
$Re[v_c(0)] > 0$. The term $\frac{1}{-\I-\I v_c(t)}$ is subtracted so that
$\omega_-(x,t)$ has zero mean value on $x \in [-\pi,\pi]$,
$\int_{-\pi}^\pi{\omega_-(x,t)dx}=0$. We supplement $\omega_-(x,t)$ from
\e{omegagamma0_periodic_lower} with
$\omega_+(x,t)=\overline{\omega_-(\bar{x},t)}$ to get a real-valued
solution $\omega=\omega_-+\omega_+$ of \e{mainEquationa0sigma0_periodic}.

From (\ref{omegagamma0_periodic_lower}) we  compute norms
 \begin{multline} \label{a0sigma0_periodic_solution_norms}
 \| \omega(x,t) \|_{L^\infty}=  \frac{|\omega_{-1}(t)|-Im[\omega_{-1}(t)]}{Re[v_c(t)]}+ 2\frac{ Im[\omega_{-1}(t)] (1 + Re[v_c(t)]) - Im[v_c(t)] Re[\omega_{-1}(t)]}{(Im[v_c(t)])^2 + (1 + Re[v_c(t)])^2},\\
 \| \omega(x,t) \|_{L^2} = 2\sqrt{\pi}\frac{|\omega_{-1}(t)|}{\sqrt{Re[v_c(t)]((Im[v_c(t)])^2 + (1 + Re[v_c(t)])^2)}},\\
 \| \omega(x,t) \|_{B_0} =\frac{|\omega_{-1}(t)|}{Re[v_c(t)]}\left(\sqrt{\frac{(Im[v_c(t)])^2 + (1 - Re[v_c(t)])^2}{(Im[v_c(t)])^2 + (1 + Re[v_c(t)])^2}} +1\right).
 \end{multline}

Substituting \e{omegagamma0_periodic_lower} to
\e{mainEquationa0sigma0_periodic_lower} we get the following equations:
\begin{equation}\nonumber
\frac{d\omega_{-1}(t)}{dt}=-\nu\omega_{-1}(t) +
\frac{2\omega_{-1}^2(t)}{1+v_c(t)},  \qquad  \frac{dv_c(t)}{dt}
=\omega_{-1}(t),
\end{equation}
and their solution:
\begin{equation}\label{a0sigma0_periodic_solution}
\omega_{-1}(t)=\omega_{-1}(0)e^{-\nu t} \left(\frac{1-e^{-\nu
t_0}}{1-e^{-\nu (t+t_0)}}\right)^2,  v_c(t)=\frac{\omega_{-1}(0)}{\nu}
\frac{(1-e^{-\nu t_0})(1-e^{-\nu t})}{1-e^{-\nu(t+ t_0)}}+v_c(0),
\end{equation}
where $e^{\nu t_0}=1-\frac{\nu(1+v_c(0))}{\omega_{-1}(0)}$  is complex
valued in general.

Equations \e{a0sigma0_periodic_solution} for $\nu=0$ reduce to:
\begin{equation}\label{a0sigma0_periodic_solution_nu0}
\omega_{-1}(t)=\frac{\omega_{-1}(0)}{(1-\frac{\omega_{-1}(0)}{1+v_c(0)}t)^2},
\qquad
v_c(t)=\frac{v_c(0)+\frac{\omega_{-1}(0)}{1+v_c(0)}t}{1-\frac{\omega_{-1}(0)}{1+v_c(0)}t}.
\end{equation}
For the case $\nu=0$ we always have a collapsing solution (even for
arbitrarily small data) at $t=t_c$ with collapse location
\begin{equation} \label{a0sigma0_xc_periodic}
x=x_c=2 \tan^{-1}(\I v_c(t_c))=2 \tan^{-1}(- Im[v_c(t_c)]),
\end{equation}
since $v_c(\infty)=-1<0$ for any $\omega_{-1}(0)$. 

We have $Re[v_c(t_c)]=0$ at time
\begin{equation} \label{a0sigma0_tc_periodic}
t_c=\frac{X + \sqrt{4|\omega_{-1}(0)|^2 Re[v_c(0)] (| v_c(0) | ^2 + 1 + 2
Re[v_c(0)]) + X^2}}{2|\omega_{-1}(0)|^2},
\end{equation}
$$X=Re[\omega_{-1}(0)](1 + Im[v_c(0)]^2 - Re[v_c(0)]^2) - 2 Im[\omega_{-1}(0)] Im[v_c(0)] Re[v_c(0)].$$

For purely real $v_c(0)$ and $\omega_{-1}(0)$ equations
\e{a0sigma0_tc_periodic} and \e{a0sigma0_xc_periodic} reduce to
\begin{align} \nonumber
t_c&=-\frac{v_c(0)(1+v_c(0))}{\omega_{-1}(0)}, &&x_c=0  &\mbox{for}\quad \omega_{-1}(0)<0,\\
t_c&=\frac{1+v_c(0)}{\omega_{-1}(0)}, &&x_c=\pm \pi  &\mbox{for} \quad
\omega_{-1}(0)>0.\nonumber
\end{align}

 The solution \e{omegagamma0_periodic_lower}, \e{a0sigma0_periodic_solution_nu0} is a periodic analog of equation (32) in \cite{Lushnikov_Silantyev_Siegel}, which describes  self-similar blow up in the inviscid problem. The solution  \e{omegagamma0_periodic_lower} with $\omega(x,t)=\omega_-(x,t)+\overline{\omega_-(\bar{x},t)}$ belongs to the general self-similar form (\ref{self-similar1}) with $\alpha=\beta=1.$

When $\nu>0$ the same analysis as (a)-(c) in Section
\ref{sec:vCLMa0sigma0exact} can be done. In this case
\e{a0sigma0_periodic_solution} gives either global existence of the
solution \e{omegagamma0_periodic_lower} or a collapsing solution depending
on initial values of $\omega_{-1}(0)$ and $v_c(0)$. For simplicity,  we
assume that  $\omega_{-1}(0)$ and $v_c(0)$ are purely real. Then
\e{a0sigma0_periodic_solution_norms} becomes
\begin{align} \label{a0sigma0_periodic_solution_norms_real}
 \| \omega(x,t) \|_{L^\infty} &=  \frac{|\omega_{-1}(t)|}{v_c(t)}, \qquad \; \; \,
\| \omega(x,t) \|_{L^2} = 2\sqrt{\pi}\frac{|\omega_{-1}(t)|}{\sqrt{v_c(t)}(1 + v_c(t))}, \nonumber \\
\| \omega(x,t) \|_{B_0} &=\frac{|\omega_{-1}(t)|}{v_c(t)}\left(\frac{|1 -
v_c(t)|}{1 + v_c(t)} +1\right).
\end{align}

Rewriting the second equation of \e{a0sigma0_periodic_solution} we get
\begin{equation}\nonumber
v_c(t)=\frac{\omega_{-1}(0)(1-e^{-\nu t})+\nu
v_c(0)(1+v_c(0))}{-\omega_{-1}(0)(1-e^{-\nu t})+\nu (1+v_c(0))},
\end{equation}
from which we can conclude that:

(a) If $-\nu v_c(0)(1+v_c(0))< \omega_{-1}(0) < \nu (1+v_c(0))$ then
$0<v_c(t)<\infty$ for all $t>0$ and the solution
\e{omegagamma0_periodic_lower}, \e{a0sigma0_periodic_solution} exists for
any $t>0$.

(b1) If $\omega_{-1}(0)=-\nu v_c(0)(1+v_c(0))$ then $\omega_{-1}(t),v_c(t)
\sim e^{-\nu t} \rightarrow 0$ as $t\rightarrow \infty$, and the solution
\e{omegagamma0_periodic_lower}, \e{a0sigma0_periodic_solution} exists for
all $t>0$ because the poles approach the real line at $x_c=0$
exponentially in time.

(b2)  If $\omega_{-1}(0)=-\nu (1+v_c(0))$ then $\omega_{-1}(t),v_c(t) \sim
e^{\nu t} \rightarrow \infty$ as $t\rightarrow \infty$, and the solution
\e{omegagamma0_periodic_lower}, \e{a0sigma0_periodic_solution} exists for
all $t>0$ because the poles approach the real line at $x_c=2 \tan^{-1}(\I
\infty)=\pm\pi$ exponentially in time.

In both (b1) and (b2) cases $\| \omega(x,t) \|_{L^\infty}, \| \omega(x,t)
\|_{B_0}/2 \rightarrow \nu>0$, $\| \omega(x,t) \|_{L^2} \sim e^{-\nu t/2}
\rightarrow 0$ as $t\rightarrow \infty$.

(c1) If $\omega_{-1}(0) <- \nu v_c(0)(1+v_c(0))$ then the solution
\e{omegagamma0_periodic_lower}, \e{a0sigma0_periodic_solution} exists
until the collapse time $t_c$ (when the poles reach the real axis at
$v_c(t_c)=0$, $x_c=0$), where
\begin{equation}\nonumber
t_c=\frac{1}{\nu}\ln \left(\frac{\omega_{-1}(0)}{\omega_{-1}(0)+\nu
v_c(0)(1+v_c(0))}\right).
\end{equation}
Using \e{a0sigma0_periodic_solution_norms_real}, we get that blow up
occurs for any initial data satisfying
\begin{align} \label{a0sigma0_periodic_solution_norms_real_c1}
 &\| \omega(x,0) \|_{L^\infty}> \nu (1+v_c(0))>\nu, &&
 \| \omega(x,0) \|_{L^2} > 2\nu \sqrt{\pi v_c(0)},
 \nonumber
 \\
 &\| \omega(x,0) \|_{B_0} > 2\nu, \mbox{ if } v_c(0)<1, &&
 \| \omega(x,0) \|_{B_0} > 2\nu v_c(0), \mbox{ if } v_c(0)\geq1. \nonumber
\end{align}
We therefore see that $\| \omega(x,0) \|_{L^2}$ can be made arbitrarily
small by choosing small enough $v_c(0)$, but $\| \omega(x,0)
\|_{L^\infty}$ and $\| \omega(x,0) \|_{B_0}$ cannot be made similarly
small.

(c2) If $\omega_{-1}(0) > \nu (1+v_c(0))$ then the solution
\e{omegagamma0_periodic_lower}, \e{a0sigma0_periodic_solution} exists
until the collapse time $t_c$ (when the poles reach the real axis at
$v_c(t_c)=\infty$, $x_c=\pm\pi$),  where
\begin{equation}\nonumber
t_c=\frac{1}{\nu}\ln
\left(\frac{\omega_{-1}(0)}{\omega_{-1}(0)-\nu(1+v_c(0))}\right).
\end{equation}

Using \e{a0sigma0_periodic_solution_norms_real}, we get that blow up
occurs for any initial data satisfying
\begin{align} \label{a0sigma0_periodic_solution_norms_real_c2}
 &\| \omega(x,0) \|_{L^\infty}> \nu \left(1+\frac{1}{v_c(0)}\right)>\nu, &&
\| \omega(x,0) \|_{L^2} > 2\nu \sqrt{ \frac{\pi}{v_c(0)}}, \nonumber
\\
 &\| \omega(x,0) \|_{B_0} > 2\frac{\nu}{v_c(0)},\mbox{ if } v_c(0)<1, &&\| \omega(x,0) \|_{B_0} > 2\nu ,\mbox{ if } v_c(0)\geq1, \nonumber
\end{align}
We again see that $\| \omega(x,0) \|_{L^2}$ (but not  $\| \omega(x,0)
\|_{L^\infty}$ and $\| \omega(x,0) \|_{B_0}$)   can be made arbitrarily
small by choosing large enough $v_c(0)$.

In both (c1) and (c2), the collapse is self-similar and the solution
\e{omegagamma0_periodic_lower} together with
$\omega(x,t)=\omega_-(x,t)+\overline{\omega_-(\bar{x},t)}$ belongs to the
general self-similar form (\ref{self-similar1}) with $\alpha=\beta=1$ and
$\| \omega(x,t) \|_{L^\infty},\| \omega(x,t) \|_{L^2} \rightarrow \infty$
as $t\rightarrow t_c$, since in (c1) $\omega_{-1}(t_c)\neq0$ and $v_c\sim
(t_c-t)$, in (c2) $\omega_{-1} \sim (t_c-t)^{-2}$ and $v_c\sim
(t_c-t)^{-1}$.

Similarly to the real line solution, the kinetic energy $E_K(t)=\int
u^2(x,t)dx$ in (a), (b1), (b2) scales like $E_K(t) \sim e^{-2 \nu t}$ as
$t\rightarrow \infty$, whereas in (c1), (c2) $E_K(t_c)$ is finite for any
complex values of the parameters $\omega_{-1}(0)$ and $v_c(0)$.

Similar exact solutions can be derived with one pair of simple poles in
$\tan(\frac{x}{2})$-space  for $a=0,~ \sigma=1$,  and one pair of double
poles for $a=1/2,~ \sigma=0,1$, as periodic analogues of exact solutions
on the real line. Details are left for future work.

\section{Numerical results} \label{sec:numerical_results}

 We present the results of direct time-dependent numerical simulations of (\ref{CLM}) in both the periodic and real-line geometries. The numerical results are consistent with the analytical theory on global existence for small data in the periodic setting, and further indicate that finite-time singularities can form for sufficiently large data. They also are in quantitative agreement with   with  the exact solutions presented in Sections \ref{sec:real_line} and \ref{sec:periodic}, and give information on the stability of those solutions.

\subsection{Numerical method} \label{sec:numerical_method}

We provide a brief description of the numerical method  and the procedure
for tracking complex singularities. More details are given in
\cite{Lushnikov_Silantyev_Siegel}. In the periodic case,  (\ref{CLM}) is
numerically solved   for $x \in  \mathbb{S}=[\-\pi,\pi]$ using a
pseudo-spectral Fourier method based on the representation \be \nonumber
\omega(x,t) = \sum_{k=-N}^{N-1} \hat{\omega}_k(t) e^{ikx} \ee in terms of
$2N$ Fourier modes. Derivatives along with the periodic Hilbert transform
and the dissipation term  are computed by wavenumber multiplication  in
Fourier space. Time stepping is performed using an 11-stage explicit
Runge-Kutta method of $8th$ order \cite{RK8_CooperVerner72} with adaptive
time step determined by the condition $\Delta t = CFL \cdot  \min[\Delta
x/(a \max_x |u(x,t)|), 1/\max_x |u_x(x,t)|, (\Delta x)^\sigma/\nu]$, where
$\Delta x=\pi/N$ and the numerical constant $CFL$ is chosen as $1/16,
1/32$ or $1/64$.  This condition ensures numerical stability and that the
error in time-stepping is near round-off.

The decay of the Fourier spectrum is checked at the end of every time
step, and if $|\omega_k(t)|$ is larger than numerical round-off at $|k|
\sim N$,  the simulation is ``rewound''  one time step backward, $N$
increased by a factor of $2$ via zero padding (i.e., Fourier
interpolation),  and the time step is adjusted before time-stepping is
continued. Rewinding helps avoid accumulation of error from the tails of
the spectrum not being fully resolved.

To compute on the infinite domain, we make a change of variable \be
\label{transform} x = \tan \left( \frac{q}{2} \right), \ee which maps
$(-\pi,\pi)$ in $q$ to $(-\infty,\infty)$ in $x$.   The transformed
equations are \cite{Lushnikov_Silantyev_Siegel}
\begin{equation}\label{CLM_transformed}
\begin{split}
&\omega_t=-a (1+\cos{q}) u\omega_q +\omega[ {\mathcal H}^{q}\omega+ C_\omega^q] - \nu[(1+\cos q) \partial_q \mathcal{H}^{q}]^\sigma \omega,  \\
&(1+\cos{q})u_q=[ {\mathcal H}^{q}\omega+ C_\omega^q], ~~q\in (-\pi,\pi),
\end{split}
\end{equation}
where $\mathcal{H}^q$ is the periodic Hilbert transform in $q$
\[
{\mathcal{H}}^q f(q) = \frac{1}{2 \pi} PV \int_{-\pi}^{\pi}  f(q') \cot
\left(\frac{q-q'}{2} \right) \  dq',
\]
 and the constant $ C_\omega^q$ is determined by%
\begin{equation}\nonumber
 C_\omega^q=-\frac{1}{2\pi}\int^{\pi}_{-\pi}\omega(q')\tan(\frac{q'}{2})\D q',
\end{equation}
so that ${\mathcal H}^{q}\omega(\pm \pi)+ C_\omega^q=0$.

A pseudo-spectral method similar  to that used for the periodic case is
then employed to solve (\ref{CLM_transformed}), using the Fourier
representation in $q$ space \be \nonumber \omega(q,t) = \sum_{k=-N}^{N-1}
\hat{\omega}_k(t) e^{ikq} \ee and a modified adaptive time-step condition
$\Delta t = CFL \cdot  \min[\Delta q/(a \max_q |(1+\cos q) u(q,t)|),
1/\max_q |(1+\cos q)u_q(q,t)|, \max_q |\omega(q,t)|/\max_q |\nu[(1+\cos q)
\partial_q \mathcal{H}]^\sigma \omega(q,t)|]$.  We only consider cases in
which  $\sigma$ is a non-negative integer, so that the dissipation term
can be easily computed by wavenumber multiplication in Fourier space.

Two complementary methods are employed to detect singularities in the
complex plane. The first method uses a least squares fit of the asymptotic
Fourier decay \be \label{FCdecay} |\hat{\omega}_k(t)| \approx C(t)
\frac{e^{-\delta(t) |k|}}{|k|^{p(t)}} \ee for $|k|\gg 1$
\cite{CarrierKrookPearson1966}, where $C(t), \delta(t)>0$ and $p(t)$ are
fitting parameters. The value of $\delta(t)$ gives the distance at time
$t$ of the (single) closest complex singularity in $\omega$  to the real
line, and $p(t)$ is related to the type or power of singularity.
 If the closest singularity to the real line has the power law form $(q- q_c)^{-\gamma}$, then $\delta=|Im(q_c)|$  and $p=1-\gamma$.
 On the infinite domain with the  additional transform \eqref{transform}, the fitting \eqref{FCdecay} provides the distance of the closest complex singularity of $\omega$ to the real line in $q$-space. To find the the distance to the closest singularity in $x$-space, we use $\delta_x=\tanh(\delta/2)$, when $|Re(q_c)|=0$ and $\delta_x=\coth(\delta/2)$, when $|Re(q_c)|=\pm\pi$.

This type of  Fourier fitting procedure for  tracking complex
singularities was originally proposed by Sulem et al.
\cite{SulemSulemFrischJCompPhys1983} and extended in
\cite{BakerCaflischSiegelJFM1993}, \cite{Pugh}. For more details about the
version employed here, see \cite{Lushnikov_Silantyev_Siegel}.

The second method for detecting complex singularities makes use of
analytical continuation based on rational interpolants. Specifically, we
employ a modified version of the AAA algorithm originally due to
Nakatsukasa et al. \cite{TrefethenAAA}. This has the advantage of
providing a structure of complex singularities beyond the one closest to
the real line, although $|Im(q_c)|$ and $\gamma$ can be determined more
accurately by Fourier fitting than via the  AAA algorithm.
See Section 10 of \cite{Lushnikov_Silantyev_Siegel} for more details on
the AAA algorithm employed here.

\subsection{Periodic problem} \label{sec:periodic_numerics}

\begin{figure}[ht!]
\begin{center}
\includegraphics[width=0.495\textwidth]{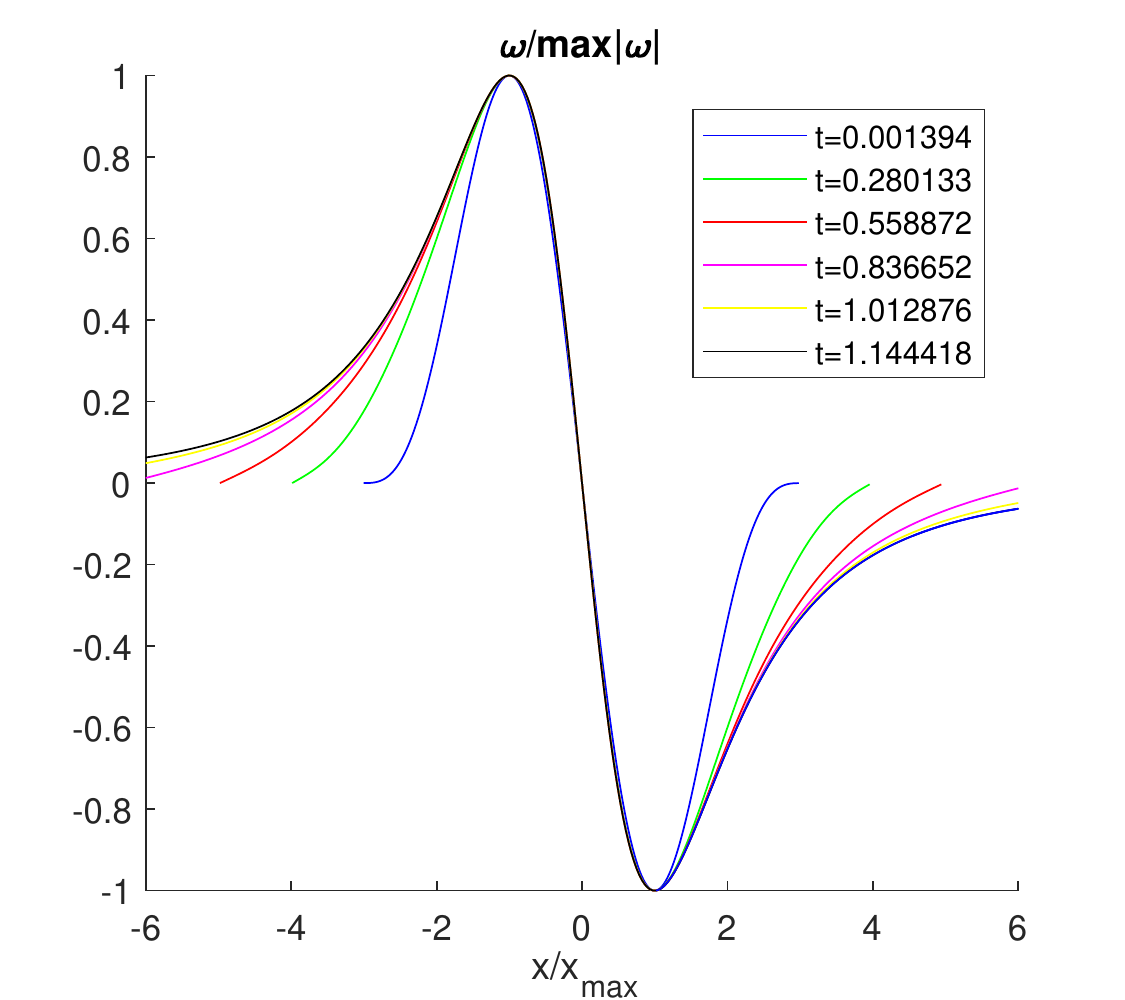}
\includegraphics[width=0.495\textwidth]{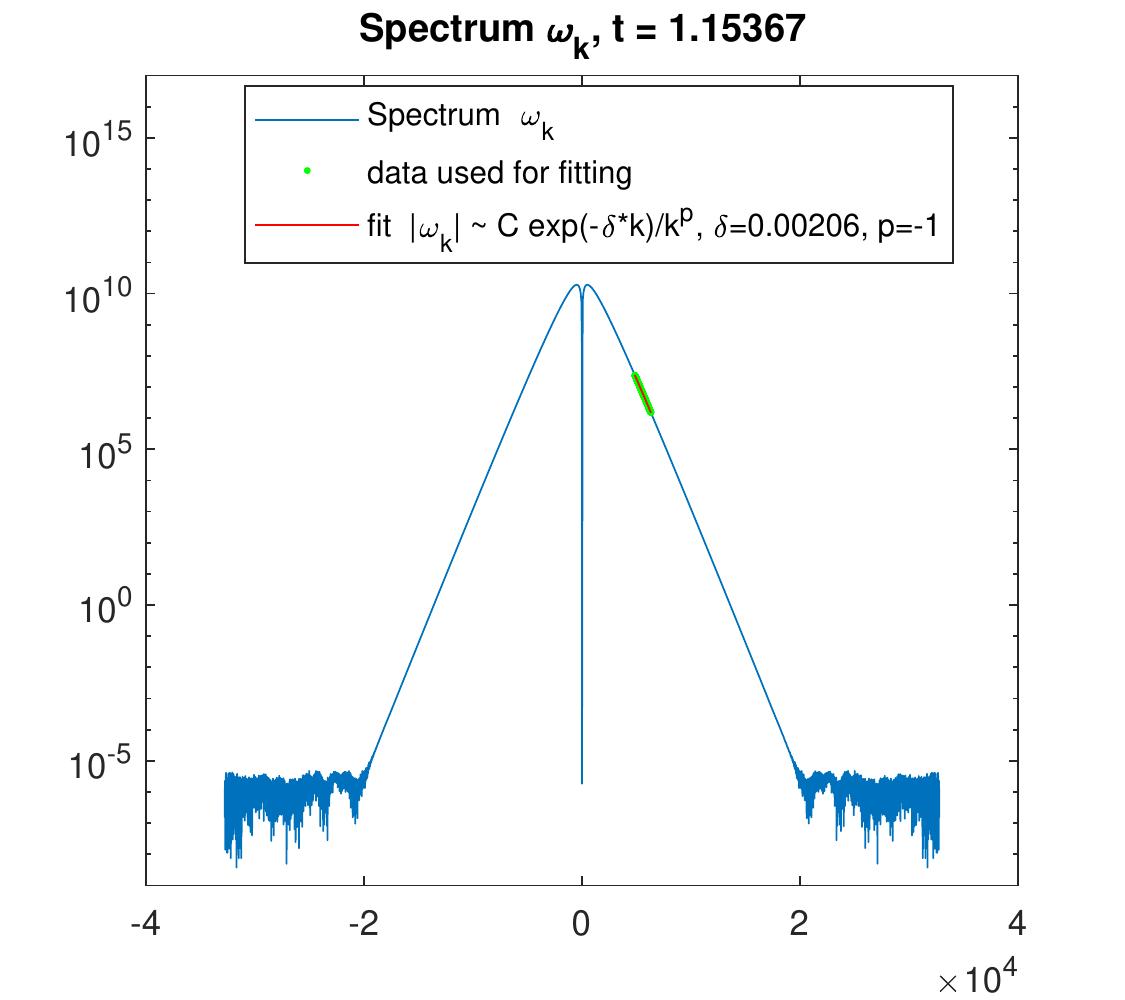}\\
\includegraphics[width=0.495\textwidth]{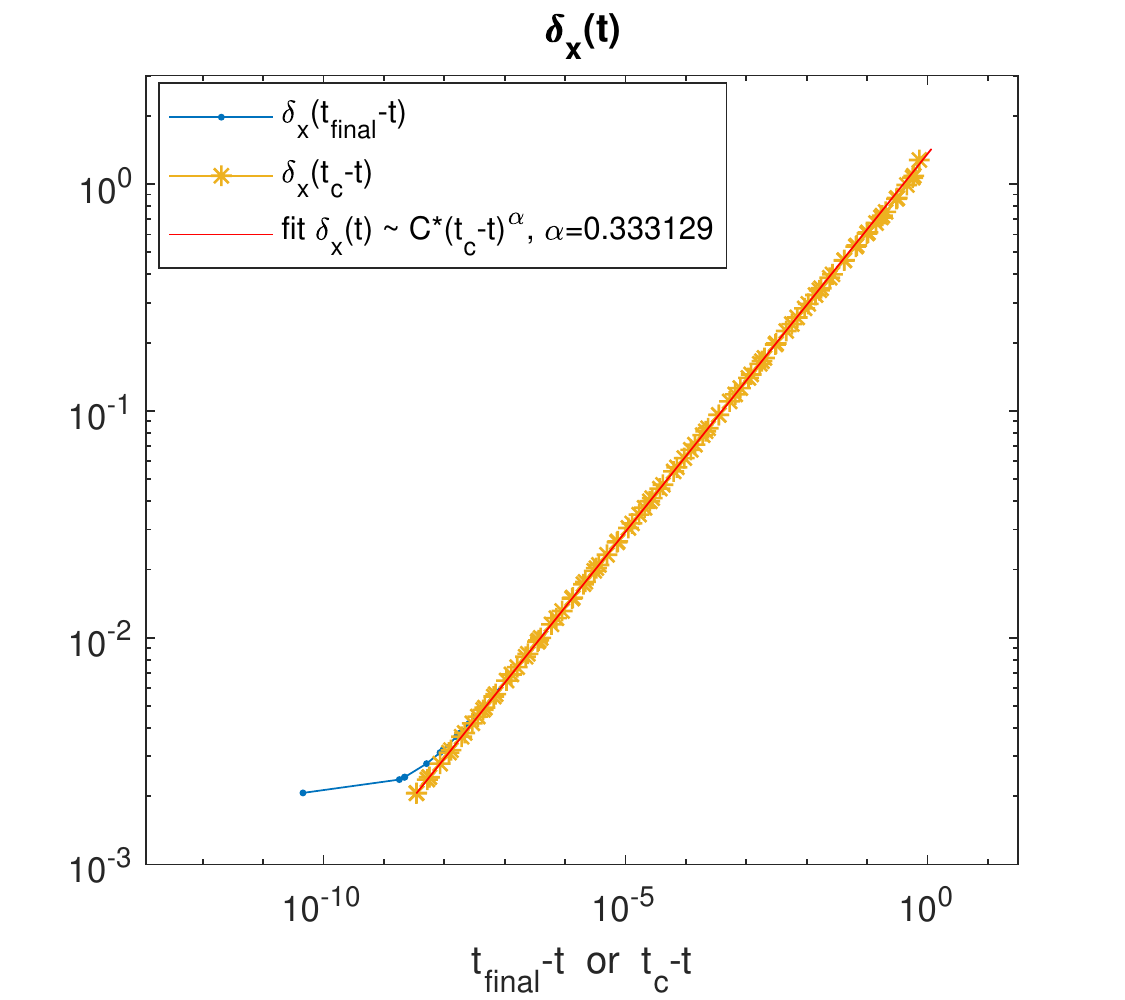}
\includegraphics[width=0.495\textwidth]{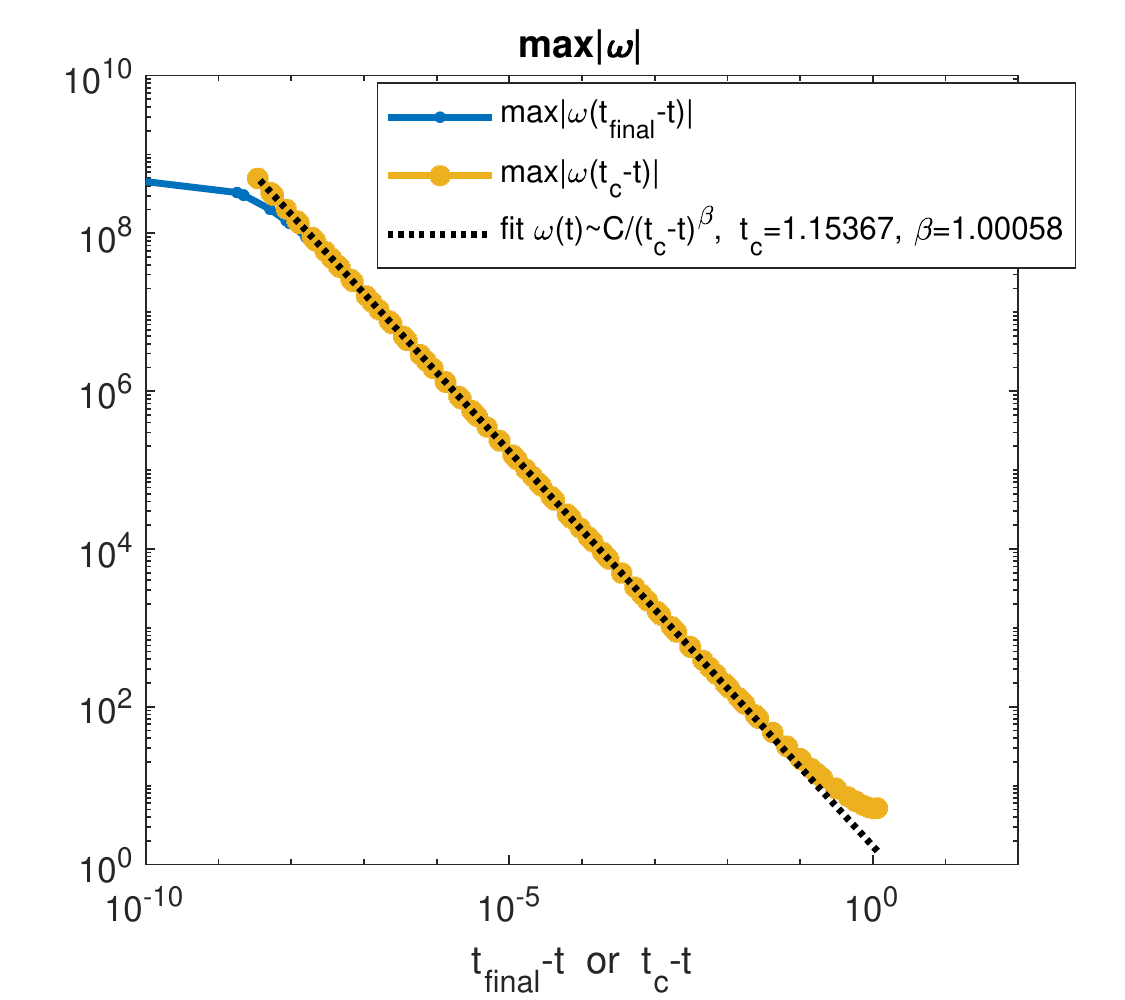}
\caption{Evolution of the collapsing periodic solution with parameters
$a=1/2$, $\sigma=1$, $\nu=1$ and two-mode initial data
(\ref{initial_data}) with $A=4$. Top-left: Scaled  solution
$\omega(x/|x_{max}|)/\max_x |\omega|$  versus $x/|x_{max}|$ (where
$x_{max}>0$ is the  location of $\max_x |\omega|$)  at different times in
the evolution. Top-right: Spectrum $\log |{\hat \omega}_k|$ versus $k$ at
$t=1.1536657$ and fit by \eqref{FCdecay}. Bottom-left: Log-log plot of
$\delta_x(t)=\delta(t)$, the distance of the closest singularity to the
real line, versus $t_{final}-t$ (in blue) and  $t_c-t$ (in yellow) and fit
to $\delta_x(t) \sim C(t_c-t)^\alpha, \alpha\approx1/3$. Bottom-right:
Log-log plot of $\max_x |\omega|$ versus $t_c-t$ and the fit $\max_x
|\omega|(t) \sim C/(t_c-t)^\beta, t_c\approx1.15367 ,\beta\approx1$.
\label{fig:a0.5sigma1_A4} }
\end{center}
\end{figure}

Figure \ref{fig:a0.5sigma1_A4} gives an illustrative example of
finite-time collapse in the case of the periodic problem with parameter
values $a=1/2$, $\sigma=1$ and $\nu=1$, for two-mode   initial data \be
\label{initial_data} \omega_0(x)=\I A \left(
\frac{1}{(\tan(\frac{x}{2})-\I)^2} - \frac{1}{(\tan(\frac{x}{2})+\I)^2}
\right) = -A (\sin x +\sin(2x)/2). \ee The top-left panel  plots a scaled
solution
 $\omega(x/|x_{max}|)/\max_x |\omega|$  versus $x/|x_{max}|$ at different times in the evolution.  The solution curves approach a universal self-similar profile $f(\xi)$  in a space-time neighborhood of the collapse point. This verifies the self-similar nature of the collapse.
 The top-right panel show the spectrum $\hat{\omega}_k$ of the solution at $t=1.1536657$ and fit by \eqref{FCdecay}. We choose an interval of $k$ somewhere between $1/4$ and $1/3$ of the full length of the spectrum to obtain the best balance between   numerical precision  and asymptotic or  large-$k$ behavior in the data. The fit $p \approx -1$ indicates the presence of  a persistent double-pole in $\mathbb{C}$ for this periodic geometry, similar to the  exact solution in the infinite geometry (cf. \S \ref{a=1/2_semi_analytical}).
 The bottom-left panel presents a log-log plot of  $\delta_x(t)=\delta(t)$, the distance of the closest singularity to the real line,  versus $t_c-t$ ($\delta_x(t)$ versus raw time $t_{final}-t$ is also shown, where  $t_{final}<t_c$ is the final simulation time).
 The linear behavior in this log-log plot indicates an algebraic approach of the singularity toward the real line when $t$ is near $t_c$, and a least squares fit to $\delta_x(t) \sim C(t_c-t)^\alpha$ gives  the similarity parameter $\alpha \approx 1/3$.
 The bottom-right panel shows a log-log plot of $\max_x |\omega|$ versus $t_c-t$, which shows $\max_x |\omega|(t) \sim C/(t_c-t)^\beta$  behavior near the singularity time.
 To estimate $t_c$ and $\beta$, we found it most accurate and reliable to fit to  $\max_x |\omega|(t) \sim C/(t_c-t)^\beta$  using the last quarter of $k$-space data for $|\omega|(t)$. This fit gives $t_c \approx 1.15$ and
 similarity parameter $\beta \approx 1$.  Note that the maximum value  $\max\limits_q|\omega(q,t)|$  of the numerical solution increases from an initial value $\sim 10$ up to $\sim 10^9$ at the final simulation time $t_{final}$.

\begin{figure}[ht!]
\begin{center}
\includegraphics[width=0.495\textwidth]{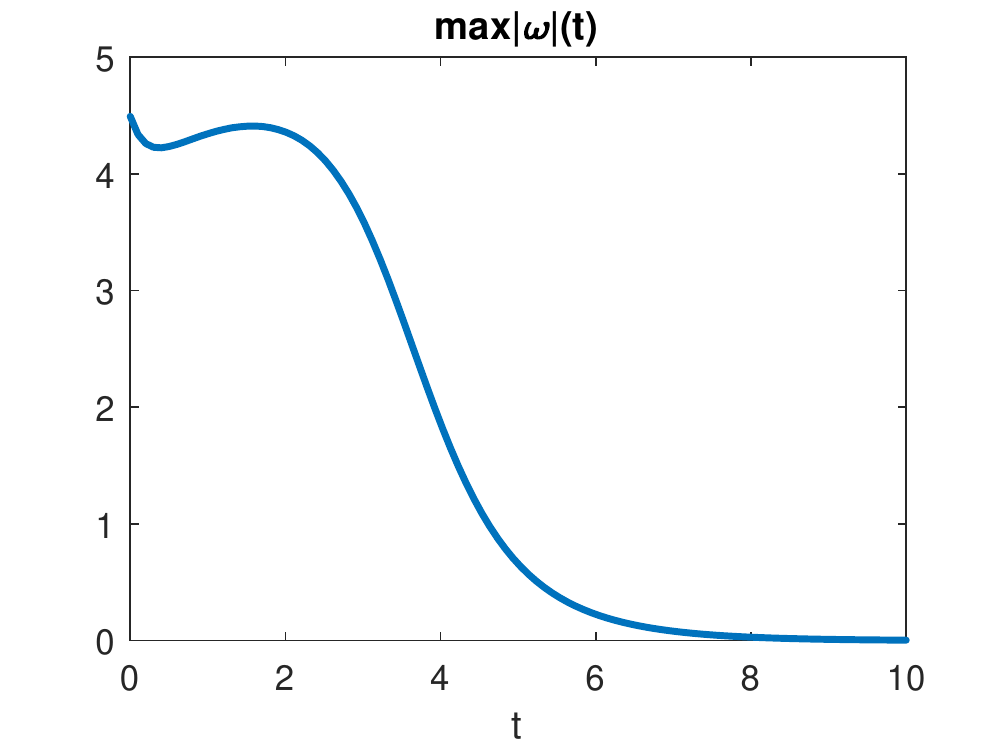}
\includegraphics[width=0.495\textwidth]{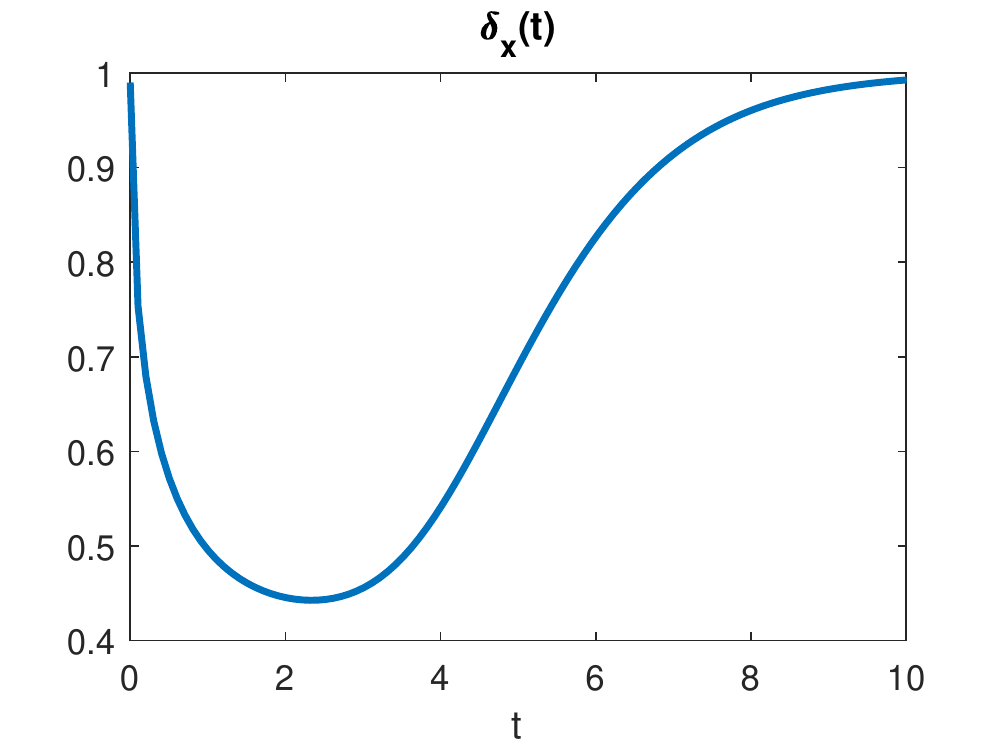}
\caption{Evolution of $\max_x |\omega|(t)$ and of the distance to the
closest singularity $\delta_x(t)$ from the real line in $x$-space for the
periodic solution with parameters $a=1/2$, $\sigma=1$, $\nu=1$ and
two-mode initial data (\ref{initial_data}) with
$A=3.46$.\label{fig:a0.5sigma1_A3.46}}
\end{center}
\end{figure}

 The fitted  values of $\alpha$ and $\beta$ are the same as for
 the exact solution on the real line (\ref{similarity_soln1}, \ref{similarity_solution_omega}). This is expected, since the local form of a collapsing similarity solution does not depend on the far-field boundary conditions, i.e., whether they are posed on $x \in \mathbb{R}$ or $\mathbb{S}$.  Notably, collapse  is only observed for $A \geq 3.47$, and the numerics suggest that there is global existence when $A \leq 3.46$, as illustrated in Figure \ref{fig:a0.5sigma1_A3.46}.

 The singularity structure in $\mathbb{C}$, as determined by the AAA algorithm,  consists of two double poles at $x=\pm i \delta$  and two branch cuts coming out of them vertically.
 There are also two more branch points at $x=\pm \pi \pm i \delta_2$ with $\delta_2>\delta$.  This singularity structure, as well as the similarity  exponents $\alpha$ and $\beta$,  are the same as in the problem without dissipation \cite{Lushnikov_Silantyev_Siegel}.
 However, the collapse takes longer to develop  when there is dissipation, e.g., $t_c=1.15367$ compared to $t_c=0.491637$ in the inviscid problem when $A=4$.
A more important distinction is that  collapse can occur in the inviscid
problem for any amplitude $A$, with the collapse time found to scale like
$t_c \sim 1/A$.

 \begin{table}[ht!]
 \centering
 \begin{tabular}{|c |c| c| c|}
 \hline
 $a$ & $\sigma$ & no blow up&  blow up \\
 \hline
        & $2$    & $A \leq 18.4$  & $18.5  \leq A$ \\
        & $1$    & $A \leq 4.53$  & $4.54  \leq A$ \\
 $0$    & $1/2$  & $A \leq 2.35$  & $2.35  \leq A$ \\
        & $0$    & $A \leq 1.33$  & $1.34   \leq A$ \\
        & $-1/2$ & $A \leq 0.81$  & $0.82 \leq A$ \\
 \hline
        & $2$    & $A \leq 5.49$  & $5.50  \leq A$ \\
        & $1$    & $A \leq 3.46$  & $3.47  \leq A$ \\
 $1/2$  & $1/2$  & $A \leq 2.64$  & $2.65  \leq A$ \\
        & $0$    & $A \leq 1.96$  & $1.97  \leq A$ \\
        & $-1/2$ & $A \leq 1.42$  & $1.43  \leq A$ \\
 \hline
        & $2$    & $A \leq 6.66$  & $6.67  \leq A$ \\
        & $1$    & $A \leq 4.73$  & $4.74  \leq A$ \\
 $0.8$  & $1/2$  & $A \leq 4.01$  & $4.02  \leq A$ \\
        & $0$    & $A \leq 3.42$  & $3.43  \leq A$ \\
        & $-1/2$ & $A \leq 2.93$  & $2.94  \leq A$ \\
 \hline
 \end{tabular}
 \caption{Critical  amplitude $A$ for blow up starting from  initial data \eqref{initial_data} in the  periodic problem with $\nu=1$. The blow up is of collapsing type (i.e., with $\alpha>0$) for $a=0,~1/2$, and neither collapsing nor expanding type (i.e., with $\alpha=0$) for $a=0.8$.
 }
 \label{table1}
 \end{table}

 The dependence of the critical initial amplitude $A$ for blow up on the dissipation exponent  $\sigma$, starting from initial  data (\ref{initial_data}),   is shown in Table \ref{table1}. The critical amplitude decreases with $\sigma$, as expected, but most importantly for all values of $\sigma$ in the table we find that there is no blow up for sufficiently small data.
 This differs from the inviscid problem, in which blow up can occur for arbitrarily small amplitude. Of course, the absence of a blow up for $\sigma<1$ and sufficiently small data could be the consequence of restricting to the particular class of initial conditions \e{initial_data}. In fact, for $a=0,~ \sigma=0,~ \nu>0$ we have established in Section \ref{sec:periodic} that  blow up occurs for arbitrarily small data of type \e{omegagamma0_periodic_lower} in the  $L^2$ norm (but not in $L^\infty$ or $B_0$ norms). That data  contains a pair of simple poles in the finite complex plane.
 Additional numerical simulations have been performed with the data
 \e{omegagamma0_periodic_lower} which  (1)  validates the analytical solution described in Section \ref{sec:periodic} in cases (a)-(c2) and  confirms the formulas for $t_c$, $x_c$ and $\| \omega(x,t) \|_{L^\infty}$, $\| \omega(x,t) \|_{L^2}$, $\| \omega(x,t) \|_{B_0}$,  $E_K(t)$; and (2) shows that blow up for this data  does not occur when  $\sigma \geq 1$ and the data is  sufficiently small in the  $L^2,~ L^\infty$ and $B_0$ norms. This is consistent with the analytical theory.

 Examination of the solution at other values of $a$  and $\sigma\geq1$ gives results that are consistent with Theorems \ref{thm_ambrose} and \ref{thm_main}, namely, that finite-time singularity formation in the periodic problem does not occur for sufficiently small data  when using initial conditions of type \e{omegagamma0_periodic_lower} or \e{initial_data}.

 \subsection{Problem on the real line} \label{sec:real_line_numerics}

 In contrast to the periodic case, the problem on the real line can exhibit finite-time blow up for arbitrarily small data.

 \subsubsection{Schochet's solution for $a=0$, $\sigma=2$}

 We have numerically computed the solution to the initial value problem   (\ref{CLM}) on $x \in \mathbb{R}$ using Schochet's initial condition.  The initial singularity locations are taken on the negative imaginary axis in $\mathbb{C}$, e.g.,  $x_1(0)=-i, x_2(0)=-2i$,  and data for  $\omega=2 Re (\omega_+)$ is specified as in (\ref{Schochet_soln}), with $A(0),B(0),C(0),D(0)$ given by  (\ref{A_eqn}), (\ref{B_eqn}). We use the corrected values  $K_\pm=24(3 \pm  \sqrt{6})$.  In all cases we observe  singularity motion exactly as given by
 (\ref{x1_eqn}), (\ref{x2_eqn}). We also
observe  self-similar collapse which scales precisely as predicted by
(\ref{eqn:selfsim}), with collapse time given by (\ref{TcSchochet}). This
verifies the corrected form of Schochet's solution, and shows that it is
stable to discretization and round-off errors. Crucially, this solution
develops finite time singularities from arbitrarily small data.

Perturbations of Schochet's initial data, for example by slightly altering
some of the coefficients $K_\pm$ or $A(0)$ through $D(0)$, leads to the
formation of additional branch points/cuts in the complex singularity
structure. In particular, we observe the formation of a branch cut between
the initial  two double poles at $x_1(t)$ and $x_2(t)$ in each of
$\omega_+$ and $\omega_-$ when $t>0$. Despite the change in the complex
singularity structure, the solution exhibits the same self-similar blow-up
as described by (\ref{eqn:selfsim}), with similarity exponents $\alpha=1$
and $\beta=2$. For small perturbations in the data, the value of $t_c$ is
only slightly perturbed from (\ref{TcSchochet}).

\subsubsection{Solutions for $a=1/2$, $\sigma=1$}

Numerical computations of the initial value problem for $a=1/2$ and
$\sigma=1$ using double pole data of the form (\ref{double_pole}) have
also  been performed. These give results that are in complete quantitative
agreement with the analytical solution described in Section
\ref{a=1/2_semi_analytical}. In particular, we find that the complex
singularity pattern for $t>0$ consists of two double poles as described by
(\ref{double_pole}). We also find that there is blow up with the local
self-similar form (\ref{similarity_solution_omega})
when $\Omega(0)=\omt(0)/v_c(0) \geq 2$, and global existence with $\omt(t)
\rightarrow 0$ and $v_c(t)  \rightarrow t+c$ when $\Omega(0)<2$, identical
to  Figure \ref{fig_a=1/2_sigma=1_exact}. This numerically validates the
analysis of Section \ref{a=1/2_semi_analytical},
 and further shows that the analytical solution derived  there  is stable to discretization and round-off errors. We have additionally  verified that this solution develops finite time singularities from arbitrarily small data as measured by the $L^2$ or $L^\infty$ norms of $\omega$, by taking  the imaginary singularity location $v_c(0)$ large while retaining $\Omega(0) >2$.

 Small perturbations of the initial data (\ref{double_pole}) have no effect on the complex singularity pattern of the self-similar part of the collapsing solution near $t \rightarrow t_c$. Blow up continues to follow the self-similar form  (\ref{similarity_solution_omega}) with $\alpha=1/3$ and $\beta=1$ and two double poles being the closest singularity to the real line.

 \subsubsection{Solutions for $a=0$, $\sigma=1$} \label{sec:num_a=0_sigma=1}
 We performed numerical computations for the initial value problem with $a=0$ and $\sigma=1$ and initial data of the form (\ref{omegagamma1}) containing a pair of simple poles.
 The computations validate (i.e., agree quantitatively) with the analytical solution \eqref{omeg1pvc0} described in Section \ref{sec:vCLMa0sigma1exact} by showing (1)  global existence if $Re(\omega_{-1}(0))>-\nu$, (2) stationarity if $\omega_{-1}(0)=-\nu$ and (3) collapse  with self-similar form (\ref{self-similar1}) and $\alpha=\beta=1$ if $Re(\omega_{-1}(0))<-\nu$.

\begin{figure}[ht!]
\begin{center}\label{fig:a0sigma1_real_line}
\includegraphics[width=0.41\textwidth]{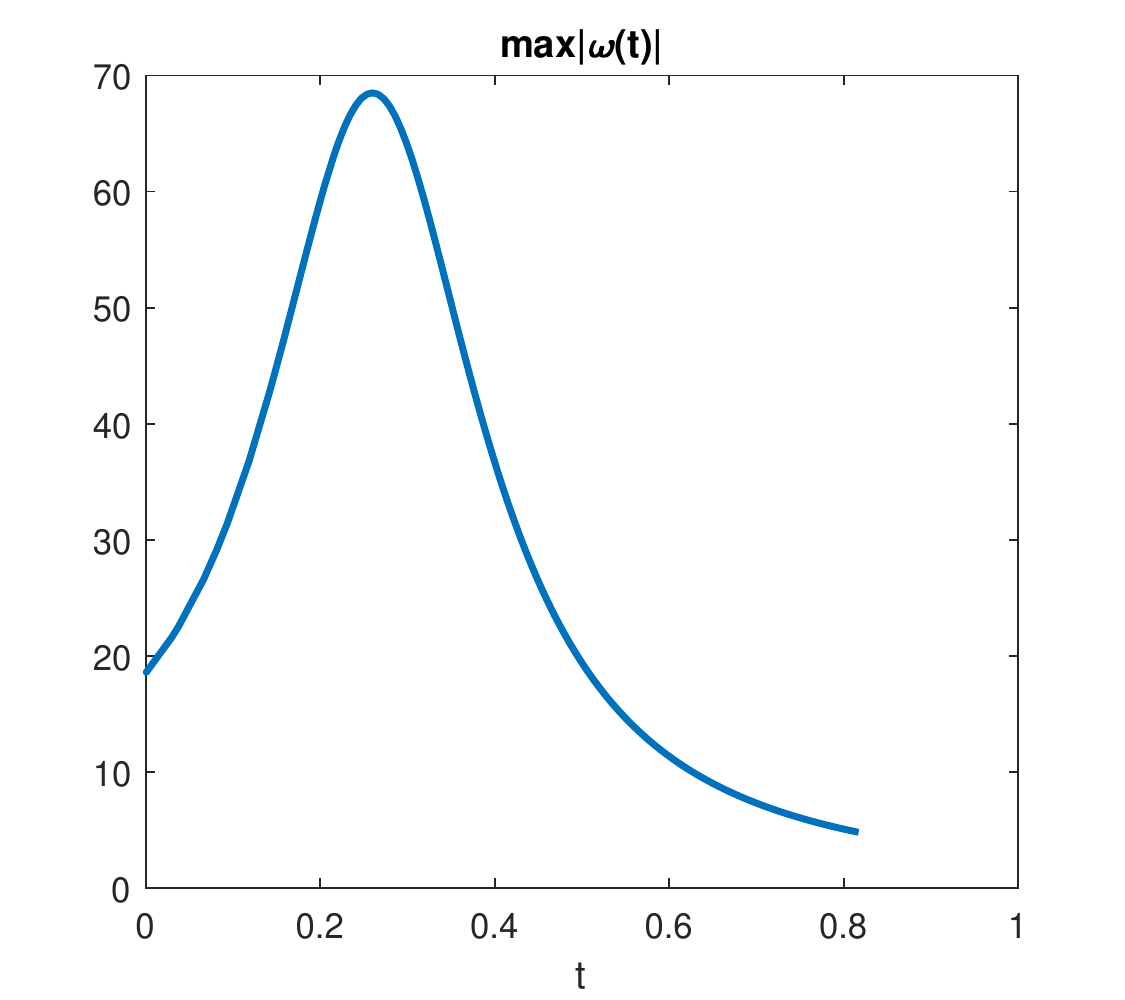}
\includegraphics[width=0.41\textwidth]{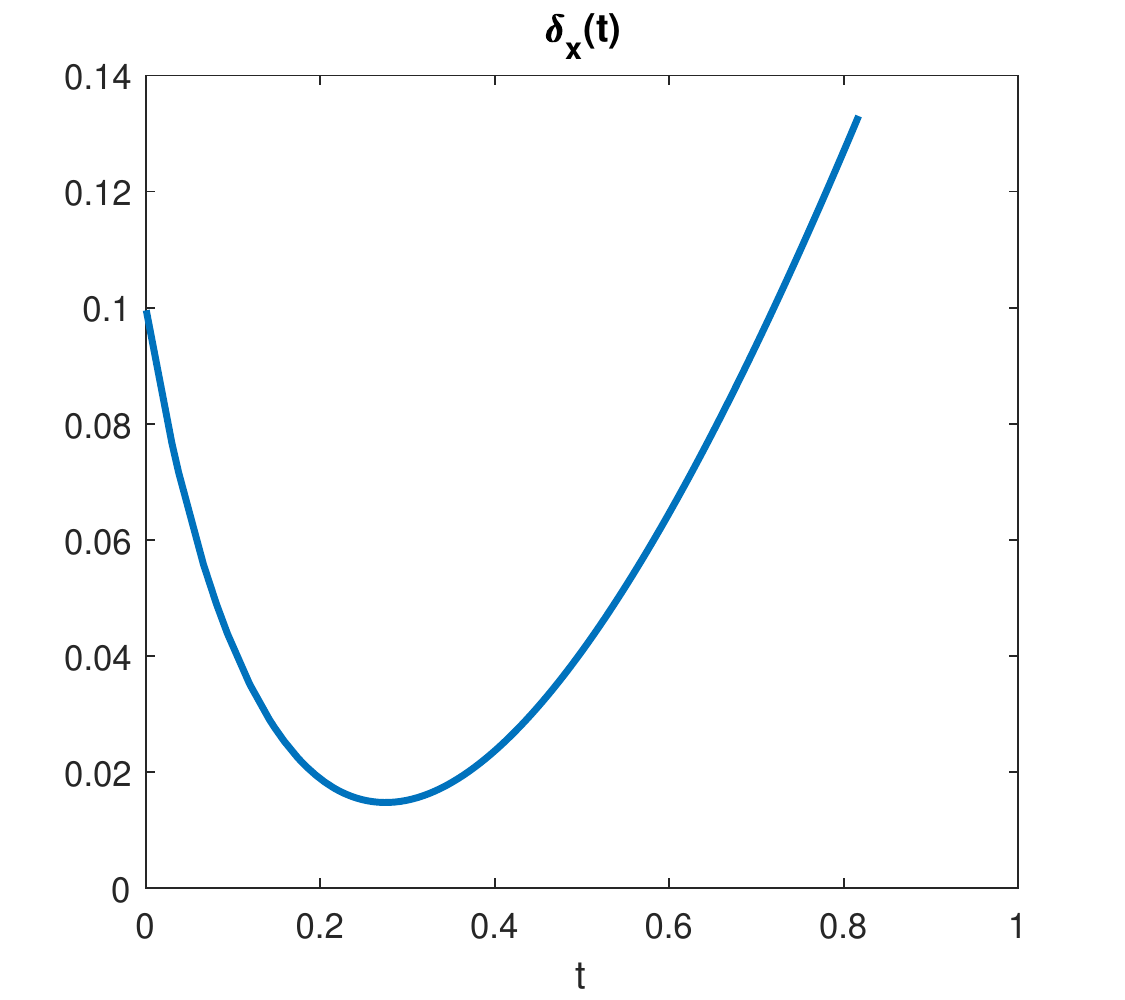}\\
\includegraphics[width=0.41\textwidth]{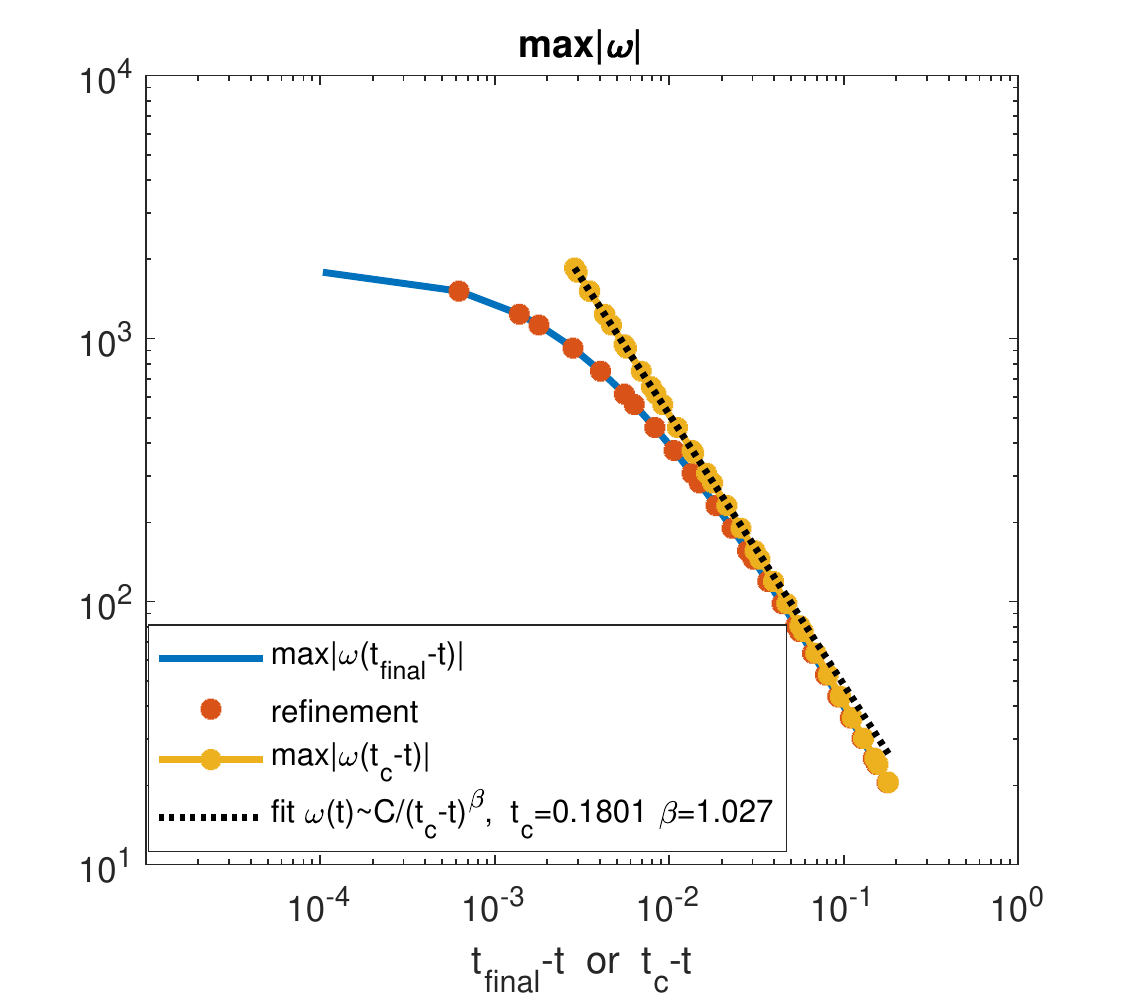}
\includegraphics[width=0.41\textwidth]{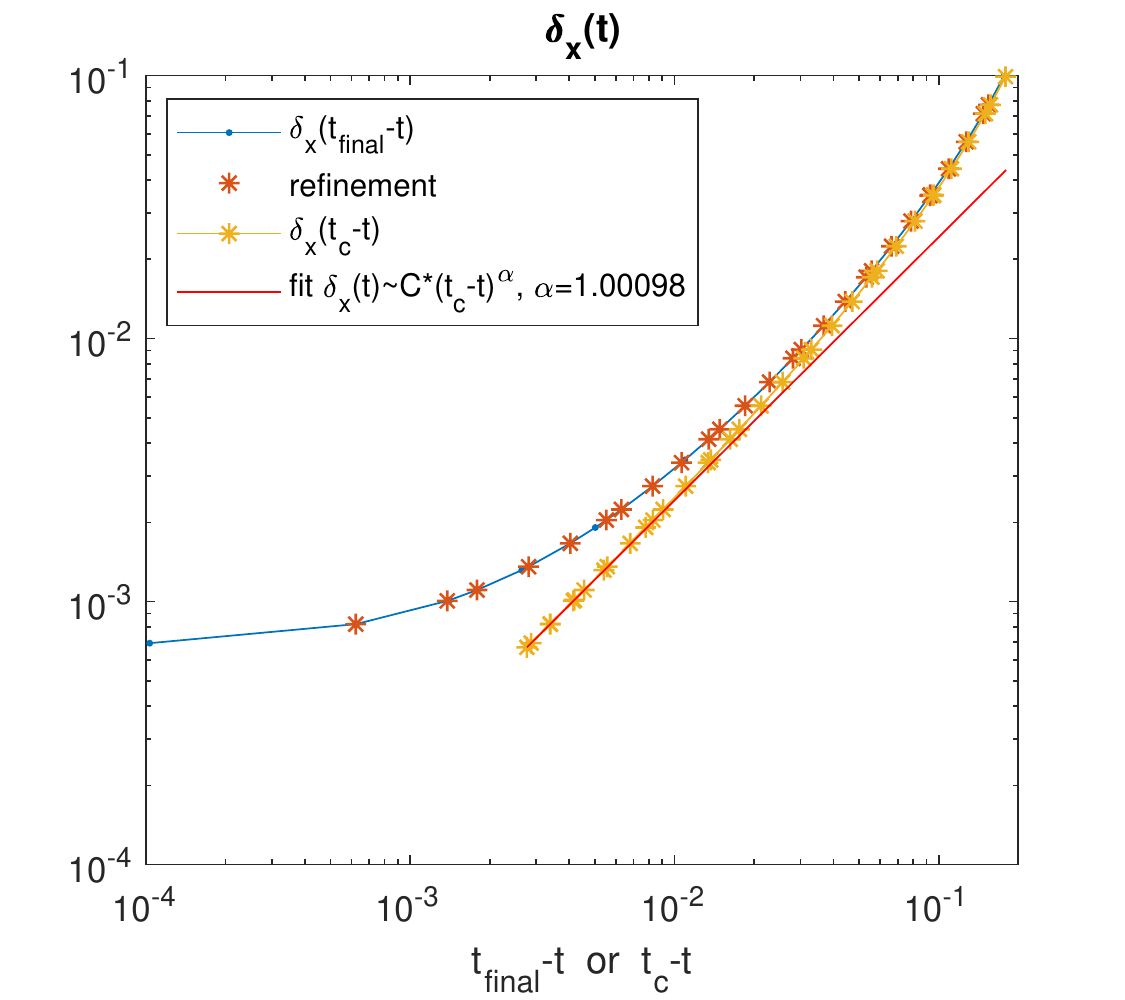}\\
\includegraphics[width=0.41\textwidth]{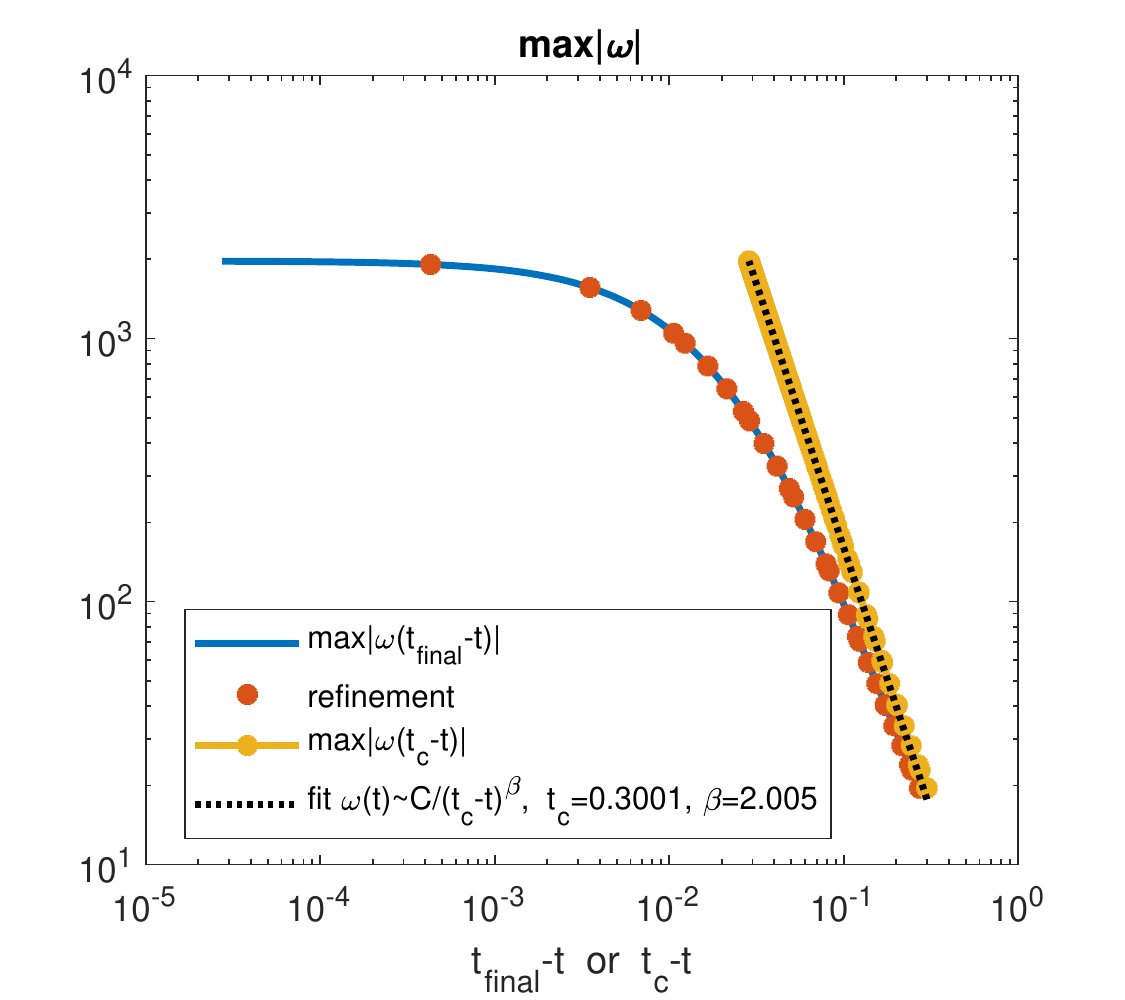}
\includegraphics[width=0.41\textwidth]{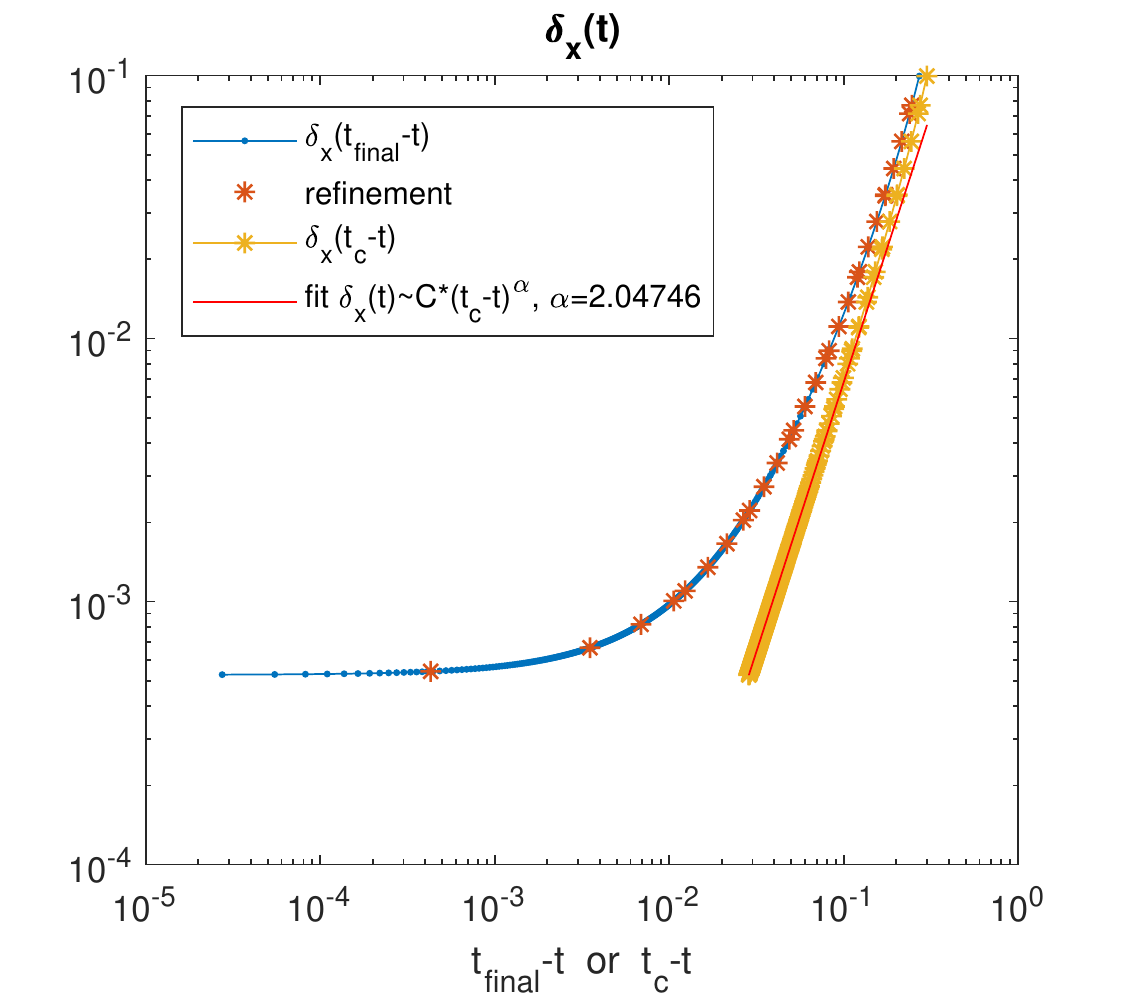}
\caption{Evolution of $\max_x |\omega|(t)$ and
$\delta_x(t)$
during numerical simulations with  $a=0$, $\sigma=1$, $\nu=1$ and initial
data (\ref{omegagamma2pairssigma1}),  which contains two pairs of simple
poles. The initial pole positions and amplitudes are $v_{c,1}(0)=0.1,
v_{c,2}(0)=0.9$ and  $\omega_{-1,1}(0)=-\omega_{-1,2}(0)=K$. Top:
$K=-1.9$, global existence of the solution. Middle: $K=-2.1$, self-similar
collapse in the form (\ref{self-similar1}) with $\alpha=\beta=1$. Bottom:
$K=-2$, self-similar collapse in the form (\ref{self-similar1}) with
$\alpha=\beta=2$. Note that in the case of collapse,  data for
$\mbox{max}_x |\omega|$ and $\delta_x(t)$  display linear behavior (with
slope $\beta$ and $\alpha$, respectively) when  plotted versus $t_c-t$. }
\end{center}
\end{figure}

We also performed numerical computations for the initial value problem
with $a=0$ and $\sigma=1$ using initial data of the form
(\ref{omegagamma2pairssigma1}) which contains two pairs of simple poles.
We considered the two cases $\omega_{-1,1}(0)+\omega_{-1,2}(0)\neq0$ and
$\omega_{-1,1}(0)+\omega_{-1,2}(0)=0$. The computations agree
quantitatively with  the analytical solutions described in Section
\ref{sec:vCLMa0sigma1exact}.  We are able to observe self-similar blow up
in the form (\ref{self-similar1}) with  $\alpha=\beta=1$ and with
$\alpha=\beta=2$ (see Figure \ref{fig:a0sigma1_real_line}).
 However, the latter type of blow up is unstable, in the sense that arbitrarily small perturbations of the initial data transform it  into blow up of type (\ref{self-similar1}) with $\alpha=\beta=1$ or lead to no collapse at all (see the top of Figure \ref{fig:a0sigma1_real_line}).

 Also we have found and checked numerically that the initial condition $w_0(x)=2 \I A(1/(x - \I V_c)^2 - 1/(x+\I V_c)^2)$  with two double poles, where $A,V_c>0$ are real numbers, leads to similar solutions. Each of the double poles splits into two single poles at $t=0$, one of which initially moves toward the real line while the other one moves away.
 This initial condition could be obtained from (\ref{omegagamma2pairssigma1}) as a limit $\epsilon \rightarrow 0$  of
 $$\omega_{-1,1}(t=0)=-\omega_{-1,2}(t=0)\sim - \frac{A}{\epsilon}, \quad v_{c,1}(t=0)\sim V_c-\epsilon, \quad
 v_{c,2}(t=0)\sim V_c+\epsilon.$$

 Substitution  of this initial data 
 into the solution \e{omega1c_c0} gives
 $$\omega_{-1,1}(t)=-\omega_{-1,2}(t)= - \sqrt{\frac{A}{2t}},
 v_{c,1}(t)=\nu t - \sqrt{2At} + V_c,
 v_{c,2}(t)=\nu t + \sqrt{2At} + V_c. $$

 If $A>2\nu V_c$, this solution leads to  self-similar blow up in the form (\ref{self-similar1}) with $\alpha=\beta=1$ at time $t_c=(A-\nu V_c-\sqrt{(A-\nu V_c)^2-(\nu V_c)^2})/\nu^2$. At the time $t=t_c$ the lower  poles cross the real axis at $x=x_c=0$, so that $v_{c,1}(t_c)=0$, with nonzero velocity  $v'_{c,1}(t_c)<0$.

 If $A=2\nu V_c$, this solution provides a self-similar blow up
 with $\alpha=\beta=2$ at the time $t_c=A/2=\nu V_c$, when the lower  pole approaches the real axis at $x=x_c=0$ with zero velocity $v'_{c,1}(t_c)=0$.

 For $A<2\nu V_c$, the solution exists for all $t>0$ since the lower pole doesn't reach the real axis.

 \subsubsection{Solutions for $a=0$, $\sigma=0$} \label{sec:num_a=0_sigma=0}
 We performed numerical computations for the initial value problem with $a=0$ and $\sigma=0$ and initial data of the form (\ref{omegagamma0}) containing a pair of simple poles.
 The computations validate (i.e., agree quantitatively) with the analytical solution \eqref{a0sigma0_solution} described in Section \ref{sec:vCLMa0sigma0exact} which exhibits global existence if $Re[\omega_{-1}(0)]>-\nu Re[v_c(0)]$), steady states if $Re[\omega_{-1}(0)]=-\nu Re[v_c(0)]$, and self-similar collapse  with $\alpha=\beta=1$ if $Re[\omega_{-1}(0)]<-\nu Re[v_c(0)]$.  We have also numerically confirmed all other formulas and claims made in Section \ref{sec:vCLMa0sigma0exact} regarding $t_c$, $x_c$ and $\| \omega(x,t) \|_{L^\infty}$, $\| \omega(x,t) \|_{L^2}$, $\| \omega(x,t) \|_{B_0}$, $E_K(t)$.

 \section{Conclusion} \label{sec:conclusion}

 We have shown global-in-time existence of solutions to the generalized Constantin-Lax-Majda equation with dissipation, in the case of  small data  in the periodic geometry, for $\sigma \geq 1$ and  any  $a$.
 This extends
 previous results on global existence theory
 from a subset of the range $a \leq -1$   to all $a$.   Our analysis is by two complementary approaches.  The first result, Theorem \ref{thm_ambrose}, proves that the solution exists globally in time for $\sigma \geq 1$ and  sufficiently small data as measured by the Wiener norm
 $\| \omega \|_{B_0}$.  Furthermore, the solution is analytic in a strip in $\mathbb{C}$ containing the real line for any $t>0$.  The theorem also gives a lower bound  on the critical initial  magnitude of vorticity (in the Wiener norm)  for global existence.

Our second main result, Theorem \ref{thm_main}, shows global-in-time
existence for  small periodic data in  $L^2$ when $\sigma>1$. The proof
shows the solution  at any time $t>0$ exists in $H^\gamma$ for all
$1/2<\gamma<\min [1,\sigma-1/2]$. Following the approach of
\cite{GrujicKukavica},  this solution is also expected to be analytic in a
strip in the complex plane for $t>0$.

The analytical theory is complemented by numerical computations for
different $a$ and $\sigma$. The numerics are able to track the formation
and motion of singularities in the complex plane. Computations in the
periodic geometry for $\sigma \geq 1$ are always found to  indicate global
existence of solutions when the initial vorticity is below a critical
amplitude.  This is in agreement with the analytical theory. On the other
hand, the numerics shows that   finite-time blow up can occur for
sufficiently large amplitude data.  We  derive a new exact analytical
solution in the periodic geometry for $a=0$ and $\sigma=0$ which forms
finite-time singularities for arbitrarily small $L^2$ data (but not for
arbitrarily small data in $L^\infty$ or the Wiener space $B_0$). This
result is suggestive of the existence in the periodic geometry of a
critical value of $\sigma$, below which there can be singularity formation
for arbitrarily small data.

In contrast, the problem on the real line can exhibit finite-time
singularity formation for arbitrarily small data as measured by  the $L^2$
or $L^\infty$  norm of $\omega$, at least  for $\sigma=0,~1,$ and $2$ at
various $a$.  This is established by the derivation of new exact
analytical solutions for $a=0$ and $1/2$. The new solutions exhibit
interesting dynamics, which are further explored by numerical simulation.
We also revisit an analytical solution derived by Schochet \cite{Schochet}
for $a=0$ and $\sigma=2$, which leads to finite-time singularity formation
for arbitrarily small data.   A minor correction is made to the solution
(after which the analytical results agree with numerical computations) and
the solution is reinterpreted from the standpoint of self-similar blow up.

In future work, we will provide a comprehensive numerical investigation of
finite-time singularity formation for a wide range of $a$  in both the
periodic and real-line problems.  Of particular interest is the effect of
the dissipation on the critical parameter $a_c$ which separates
self-similar collapsing solutions from  expanding and `neither collapsing
nor expanding' solutions observed in the problem without dissipation
\cite{Lushnikov_Silantyev_Siegel}. Another  interesting question  is
whether $\sigma=1$ is the optimal lower bound for which  global existence
for small data can be guaranteed. A related question for the problem on
the real line is whether there exists a value  of $\sigma$ greater than
$2$ for which one can guarantee global existence for small data.  These
questions are left for future work.

\section{Acknowledgements}
D.M.A. was supported by National Science Foundation Grant DMS-1907684.
M.S. was supported by National Science Foundation Grant DMS-1909407. The
work of P.M.L. was supported by the National Science Foundation, Grant
DMS-1814619.  Simulations were performed at the Texas Advanced Computing
Center using the Extreme Science and Engineering Discovery Environment
(XSEDE), supported by NSF Grant ACI-1053575.

\appendix

\section{Appendix} \label{sec:appendix}

\subsection{Proof of the inequality \eqref{desiredBound}}\label{inequalityAppendix} \label{sec:appendix1}

Given $\theta>0,$ we want to find $C\geq0$ such that for all
$(j,k)\in\mathbb{Z}^{2},$ we have
\begin{equation}\nonumber
|k|^{\theta}\leq C(|k-j|^{\theta}+|j|^{\theta}).
\end{equation}

First, note that if $k=0,$ the inequality is satisfied for any $C\geq0.$
We now focus on the case $k\neq0.$ Notice that if $k\neq0$ then also
$|k-j|^{\theta}+|j|^{\theta}\neq0.$ The value $C$ may then be taken to be
the maximum of
\begin{equation}\nonumber
\displaystyle\frac{|k|^{\theta}}{|k-j|^{\theta}+|j|^{\theta}}=\frac{1}{\left|1-\frac{j}{k}\right|^{\theta}+\left|\frac{j}{k}\right|^{\theta}}.
\end{equation}
We define $z=j/k,$ and we seek to find the maximum value for
$z\in\mathbb{R}$ of the function
\begin{equation}\nonumber
f(z)=\frac{1}{|1-z|^{\theta}+|z|^{\theta}}.
\end{equation}

We consider this in three regions.  First, if $z\geq1,$ then
\begin{equation}\nonumber
f(z)=\frac{1}{(z-1)^{\theta}+z^{\theta}}.
\end{equation}
We note the values $f(1)=1$ and
$\displaystyle\lim_{z\rightarrow\infty}f(z)=0.$ We compute
\begin{equation}\nonumber
f'(z)=-\theta\frac{(z-1)^{\theta-1}+z^{\theta-1}}{\left((z-1)^{\theta}+z^{\theta}\right)^{2}}.
\end{equation}
If $f'(z)=0,$ then $(z-1)^{\theta-1}=-z^{\theta-1},$ and this equation has
no solutions on $[1,\infty).$  Therefore the maximum of $f$ for
$z\in[1,\infty)$ is attained at $z=1,$ and is $f(1)=1.$

Next we consider $0\leq z \leq 1.$  On this domain, the function $f$
becomes
\begin{equation}\nonumber
f(z)=\frac{1}{(1-z)^{\theta}+z^{\theta}}.
\end{equation}
The boundary values on this domain are $f(0)=1$ and $f(1)=1.$  We take the
derivative, finding
\begin{equation}\nonumber
f'(z)=-\theta\frac{-(1-z)^{\theta-1}+z^{\theta-1}}{\left((1-z)^{\theta}+z^{\theta}\right)^{2}}.
\end{equation}
Setting $f'(z)=0,$ we find a critical point at $z=1/2.$  At this point, we
have the function value $f(1/2)=2^{\theta-1}.$  So, for $z\in[0,1],$ we
have $f(z)\leq\max\{1,2^{\theta-1}\}.$

Finally, we let $z\in(-\infty,0];$ on this domain, $f$ is given by
\begin{equation}\nonumber
f(z)=\frac{1}{(1-z)^{\theta}+(-z)^{\theta}}.
\end{equation}
On this domain the boundary values are $f(0)=1$ and
$\displaystyle\lim_{z\rightarrow-\infty}f(z)=0.$ The derivative of $f$ is
\begin{equation}\nonumber
f'(z)=-\theta\frac{-(1-z)^{\theta-1}-(-z)^{\theta-1}}{\left((1-z)^{\theta}+(-z)^{\theta}\right)^{2}}.
\end{equation}
Setting $f'(z)=0,$ we find the equation
$(1-z)^{\theta-1}=-(-z)^{\theta-1}.$  There are no solutions of this, so
the maximum of $f$ on the present domain is attained at $z=0,$ and is
$f(0)=1.$

Overall, we have demonstrated that \eqref{desiredBound} holds with
$C=\max\{1,2^{\theta-1}\}.$

\subsection{Proof of Lemma \ref{lemma:int_est}} \label{sec:appendix2}

Denote the integrand in (\ref{time_int}) as $J(t,\tau)$, and decompose the
integral as
\begin{align}
\int_0^t J(t,\tau) \ d \tau &= \left( \int_0^1 + \int_1^t \right) J(t,\tau) \ d \tau \nonumber \\
&= I_1 + I_2\nonumber
\end{align}
For $t \leq 2$ (say), the integral  can be bounded by a constant that is
independent of $t$.  This follows by using $|e^{-{\hat{q} (t - \tau)}}|
\leq 1$ and making the change of variable $\theta=\tau/t$ which gives
\[
\int_0^t J(t,\tau) \ d \tau \leq t^{1 - {\hat{\alpha}} - {\hat{\beta}}}
\int_0^1 \frac{1}{(1-\theta)^{\hat{\alpha}}
\theta^{{\hat{\beta}}+{\hat{\delta}}}} \ d\theta,
\]
which is bounded for $0  \leq {\hat{\alpha}} <1$, $0 \leq {\hat{\alpha}} +
{\hat{\beta}} \leq 1$, $0 \leq {\hat{\beta}} + {\hat{\delta}}<1$ and $0
\leq t \leq 2$.

 Therefore, w.l.o.g. assume $t > 2$.  To bound $I_1$,  note that $e^{-{\hat{q}}(t-\tau)} \leq  e^{-{\hat{q}}(t-1)}$ on the integration interval  and make the change of variable $\theta = \tau/t$ to obtain
\begin{align}
I_1 &\leq  e^{-{\hat{q}}(t-1)} t^{1-{\hat{\alpha}}-{\hat{\beta}}}
\int_0^{1/t} \frac{1}{(1-\theta)^{\hat{\alpha}} \theta^{{\hat{\beta}}+{\hat{\delta}}}} \ d \theta \nonumber \\
&\leq 2^{\hat{\alpha}}  e^{-{\hat{q}}(t-1)}
t^{1-{\hat{\alpha}}-{\hat{\beta}}}
\int_0^{1/t} \frac{1}{\theta^{{\hat{\beta}}+{\hat{\delta}}}} \ d \theta \nonumber \\
&\leq \frac{2^{\hat{\alpha}}}{1 - {\hat{\beta}} - {\hat{\delta}}}
e^{-{\hat{q}}(t-1)} t^{{\hat{\delta}}-{\hat{\alpha}}},\nonumber
\end{align}
where we have used $1/(1-\theta)^{\hat{\alpha}} \leq 2^{\hat{\alpha}}$ on
$\theta \in [0,1/t]$ when $t \geq 2$. Hence $I_1<C$ for $t>2$. To bound
$I_2$,  note that in the integration interval $\tau^{\hat{\beta}}>1$ and
make the change of variable $v=t-\tau$ to obtain
\begin{align}
I_2 &\leq \int_0^{t-1}  \frac{e^{-{\hat{q}} v}}{v^{\hat{\alpha}}}  \frac{1}{(1-\frac{v}{t})^{\hat{\delta}}} \ d \tau \nonumber \\
&=
\left( \int_0^{(t-1)/2} + \int_{(t-1)/2}^{t-1} \right) \frac{e^{-{\hat{q}} v}}{v^{\hat{\alpha}}}  \frac{1}{(1-\frac{v}{t})^{\hat{\delta}}} \ dv \nonumber \\
&=J_1+J_2.\nonumber
\end{align}
$J_1$ is bounded as
\begin{align}
J_1 &\leq  \left( \frac{1}{2} + \frac{1}{2t} \right)^{- {\hat{\delta}}}
\int_0^\infty
\frac{e^{-{\hat{q}} v}}{v^{\hat{\alpha}}}   \ dv \nonumber \\
&\leq C,\nonumber
\end{align}
while $J_2$ satisfies the estimate
\begin{align}
J_2 &\leq \left( \frac{t-1}{2} \right)^{-{\hat{\alpha}}} e^{-{\hat{q}}
\left( \frac{t-1}{2} \right)} \int_{(t-1)/2}^{t-1}   \frac{dv}{(1-
\frac{v}{t})^{\hat{\delta}}}
\nonumber \\
&\leq 2^{\hat{\alpha}}  e^{-{\hat{q}} \left( \frac{t-1}{2} \right)}
\frac{t}{1 - {\hat{\delta}}}  \nonumber \\
& \leq   C.\nonumber
\end{align}
Hence $I_2<C$ for $t>2$, and the result follows.

\subsection{Proof of Lemma \ref{lemma:sumest}}

Let $H_t(k)=e^{-2 t |k|^\sigma}$, and assume $\sigma > 0$.   Then
\begin{align}
\sum_{k \in \mathbb{Z}} H_t(k)  &=1 + 2 \sum_{k=1}^\infty H_t(k) \nonumber \\
&\leq 1 + 2 \left( H_t(1) + \int_1^\infty H_t(k) \ dk \right).
\label{sum_est1}
\end{align}
Substitute $u=2 t  k^{\sigma}$  into   the integral in (\ref{sum_est1})
and estimate it as
\begin{align}
\int_1^\infty H_t(k) \ dk &= \frac{(2 t )^{-1/\sigma}}{\sigma}  \int_{2 t} ^\infty e^{-u}  u^{1/\sigma-1} \ du \nonumber \\
&\leq C  t^{-1/\sigma} e^{-t}.  \label{int_est}
\end{align}
The second inequality follows from
 elementary estimates.

Substitute (\ref{int_est}) and  $H_t(1)=e^{-2t}$ into (\ref{sum_est1}) to
obtain \be \label{sum_est} \sum_{k \in \mathbb{Z}} e^{-2 t |k|^\sigma}
\leq 1 + 2 \left( e^{-2t} + C e^{-t} t^{-1/\sigma} \right) \ee which can
be easily simplified to obtain the final form  (\ref{l2_est}). Finally,
note that if
the $k=0$ term in the norm
 $\| e^{- t \rho(\cdot)} \|_{l_2}$ is omitted, which in turn implies that the first $1$ in (\ref{sum_est1}), (\ref{sum_est}), and (\ref{l2_est}) can be omitted.


\end{document}